\documentclass[a4paper,oneside,11pt]{article}
\usepackage{amsmath}
\usepackage{amsmath}
\usepackage{amsfonts}
\usepackage{amssymb}
\usepackage{titlesec}
\usepackage{graphicx}
\usepackage{subfigure}
\usepackage{soul}
\usepackage{fancyhdr}
\usepackage{float}
\usepackage{pdfsync}
\usepackage{latexsym}
\usepackage{mathrsfs}
\usepackage{mathrsfs}
\usepackage{mathtools}
\usepackage{amsthm}
\usepackage{wasysym}
\usepackage{cite}
\usepackage{color}
\usepackage{geometry}
\usepackage{extarrows}
\usepackage{enumerate}
\usepackage{dsfont}
\usepackage{bm}
\usepackage[marginal]{footmisc}
\usepackage[dvips]{epsfig}
\usepackage{amscd}
\usepackage[marginal]{footmisc}
\usepackage[bookmarks,colorlinks,citecolor=blue,linkcolor=red]{hyperref}
\allowdisplaybreaks
\def\QEDopen{{\setlength{\fboxsep}{0pt}\setlength{\fboxrule}{0.2pt}\fbox{\rule[0pt]{0pt}{1.3ex}\rule[0pt]{1.3ex}{0pt}}}} 
\def\QED{\QEDopen} 

\def\endproof{\hspace*{\fill}~\QED\par\endtrivlist\unskip}

\numberwithin{equation}{section}

\newtheorem{theorem}{Theorem}[section]
\newtheorem{lemma}{Lemma}[section]
\newtheorem{remark}{Remark}[section]
\newtheorem{proposition}{Proposition}[section]
\renewcommand{\thefootnote}{}

\title{Large Time Behavior and Sharp Interface Limit of Compressible Navier-Stokes/Allen-Cahn System for Interacting Shock Waves}
\date{  }
\author{Yazhou C{\small HEN}$^1$, Qiaolin H{\small E}$^2$, Xiaoding S{\small HI}$^{1*}$, Xiaoping W{\small ANG}$^3$ \\[3mm]
\scriptsize$^{1}$ {College of Mathematics and Physics, Beijing University of
Chemical Technology, Beijing 100029, China}\\
\scriptsize$^{2}$ {School of Mathematics, Sichuan University, Chengdu, Sichuan, 610065,  China} \\
\scriptsize$^3$ {The Chinese University of Hong Kong, Shenzhen, Guangdong, 518055, China  }
}

\begin{document}
\maketitle
\renewcommand{\thefootnote}{\fnsymbol{footnote}}
\footnotetext[1]{{Corresponding author. }\\
{Email: chenyz@mail.buct.edu.cn (Y.Chen), qlhejenny@scu.edu.cn (Q.He), shixd@mail.buct.edu.cn (X.Shi), wangxiaoping@cuhk.edu.cn, (X.Wang)}}

\begin{abstract}
In this paper, we study the large time behavior and sharp interface limit of the Cauchy problem for compressible Navier-Stokes/Allen-Cahn system with interaction shock waves in the same family. This system is an important mathematical model for describing the motion of immiscible two-phase flow. The results show that, if the initial density and velocity are near the superposition of two shock waves in the same family,  then there exists a unique global solution to the compressible Navier-Stokes/Allen-Cahn system, and this solution asymptotically converges to the superposition of the viscous shock wave and rarefaction wave which moving in opposite directions. Moreover, this global-in-time solution converges to the entropy solution of $p$-system in $L^\infty$-norm as the thickness of the diffusion interface tends to zero.
\end{abstract} 
\noindent{\bf Keywords:} Large time behavior, Sharp interface limit,  Compressible Navier-Stokes equations, Allen-Cahn equation,  Shock waves,  Rarefaction waves

\

\noindent{\bf AMS subject classifications:} 35Q35; 35B65; 76N10; 35M10; 35B40; 35C20; 76T30

\

\section{Introduction and Main Results}

\setcounter{equation}{0}

\qquad
It is well known that the immiscible two-phase flow with diffusion interface is usually described by the Navier-Stokes/Allen-Cahn system(called as NSAC system) as following
\begin{equation}\label{original NSAC}
\left\{\begin{array}{llll}
\displaystyle \rho_{t}+\textrm{div}(\rho \mathbf{u})=0,\\
\displaystyle (\rho \mathbf{u})_{t}+\mathrm{div}\big(\rho \mathbf{u}\otimes\mathbf{u}\big)-\nu(\epsilon)\Delta\mathbf{u}-\big(\nu(\epsilon)+\lambda(\epsilon)\big)\nabla\mathrm{div}\mathbf{u}+\nabla p(\rho)\\
\displaystyle \qquad\qquad\qquad\qquad\qquad \qquad\qquad=\mathrm{div}\big(\frac{\eta(\epsilon)}{2}\nabla\chi-\epsilon^2\nabla\chi\otimes\nabla\chi\big),\\
\displaystyle(\rho\chi)_{t}+\mathrm{div}\big(\rho\chi \mathbf{u}\big)=-L_d(\epsilon)\mu,\\
\displaystyle\rho\mu=\rho(\chi^3-\chi)-\eta(\epsilon)\Delta\chi,
\end{array}\right.
\end{equation}
where $\mathbf{u}$ is the fluid velocity, $\rho$  the density and $\chi$  the concentration difference of the immiscible two-phase flow. $\nabla$, div, and $\Delta$ represent the gradient operator,  the divergence operator and the Laplace operator respectively, $t$ the time and $\cdot_t=\frac{d}{dt}$. The pressure $p$ is assumed to be a function of $\rho$ given by
\begin{equation}\label{p}
p= \rho^{\gamma},\  \  \gamma>1\ \mathrm{ the \ adiabatic \ exponent}.
\end{equation}
 $\mu$ is the chemical potential.  $\nu(\epsilon)$ and $\lambda(\epsilon)$ represents the viscosity coefficient for the immiscible two-phase flow, $\eta(\epsilon)$  represents the gradient energy coefficient related to the interfacial width, and   $L_d(\epsilon)$ represents the phenomenological mobility coefficient related to the speed at which the system  approaches an equilibrium configuration, $\epsilon>0$ is the parameter. Throughout this paper, without loss of generality, we assume that these   coefficients satisfy the following relation:
\begin{equation}\label{hypothesis-3-1}
 \nu(\epsilon)\sim\epsilon,\quad\lambda(\epsilon)\sim\epsilon,\quad L_d(\epsilon)\sim\frac1\epsilon,\quad \eta(\epsilon)\sim\epsilon^2.
\end{equation}

The physical derivation of the NSAC system \eqref{original NSAC} can be referred to Blesgen \cite{B1999}, Heida-Malek-Rajagopal \cite{HMR2012} and the references therein. Intuitively, the function $\chi$ can be roughly thought as an indicator function to distinguish between two fluids. The region with $\chi=1$ refers to the region where the first fluid is located,
the region with $\chi=-1$ refers to the region where the second fluid is located, and the region of $ -1<\chi<1$ refers to the region where diffusion interface of immiscible two-phase flow is located.
On the other hand, from a purely mathematical point of view, this model can be regarded as an approximation of the two-phase flow problem with sharp interface. In this paper, we are concerned with the singular limit problem for \eqref{original NSAC} where the thickness of the diffusion interface tends to zero. Considering the complexity and challenge of this singular limit problem, this paper starts from the one-dimensional Cauchy problem of NSAC systsem \eqref{original NSAC}. Taking into account the hypothesis \eqref{hypothesis-3-1}, without loss of generality, we assume the following
\begin{equation}\label{hypothesis 3}
 2\nu(\epsilon)+\lambda(\epsilon)=\epsilon,\quad L_d(\epsilon)=\frac1\epsilon,\quad \eta(\epsilon)=\epsilon^2.
\end{equation}
Therefore, in the one-dimensional Lagrange coordinate system (see Smoller \cite{S94} and the references therein), the  problem \eqref{original NSAC} is simplified to the following mathematical model
 \begin{equation}
\left\{\begin{array}{llll}\label{nsac-lagrange}
\displaystyle v_{t}-u_{x}=0,\\
\displaystyle u_{t}+p(v)_{x}=\epsilon\Big(\frac{u_{x}}{v}\Big)_{x}-\frac{\epsilon^2}{2}\Big(\frac{\chi_{x}^2}{v^2}\Big)_{x}, \\
\displaystyle \chi_{t}=-\frac{1}{\epsilon}v\mu,\\
\displaystyle \mu=\chi^3-\chi-\epsilon^2 \Big(\frac{\chi_{x}}{v}\Big)_{x},
\end{array}\right.
\end{equation}
with the initial condition
\begin{equation}\label{initial data}
 (v,u,\chi)(x,0)=(v_0,u_0,\chi_0)(x)\xrightarrow {x\rightarrow\pm\infty} (v_{\pm}, u_\pm, \pm1),\ \ x\in\mathbb{R},
\end{equation}
where $v$ is the specific volume defined by
\begin{equation}\label{p-v}
  v=\frac1\rho,\qquad p(v)=v^{-\gamma},
  \end{equation}
and $v_\pm>0, u_\pm$  in \eqref{initial data} are the given constants, satisfy the following conditions
  \begin{equation}\label{Condition at infinity}
   u_- >u_+, \qquad \ 0<v_-<v_+.
 \end{equation}
Moreover, we assume that initial phase field $\chi_0$ satisfies the following physical assumption
\begin{equation}\label{physical constraint}
  -1\leq\chi_0\leq1.
\end{equation}

 As the parameter $\epsilon\rightarrow0^+$, that is, the diffuse interface thickness tends to zero, the system (\ref{nsac-lagrange})  is formally transformed into the following free boundary problem of hyperbolic conservation law systems
\begin{equation}\label{NSAC-Euler-equation}
\left\{\begin{array}{llll}
\displaystyle v^{\pm}_{t}-u^{\pm}_{x}=0,   &\text{in}\ \Omega^{\pm}(t),\\
\displaystyle u^{\pm}_{t}+p_{x}(v^{\pm})=0, &\text{in}\ \Omega^{\pm}(t),\\
\displaystyle \chi^{\pm}=\pm 1, &\text{in}\ \Omega^{\pm}(t),
\end{array}\right.\end{equation}
where for fixed $t\geq0$, $\Omega^{\mp}({t})$ are defined as
\begin{equation*}
 \Omega^{-}({t})\xlongequal{\mathrm{def}} \big\{{x}\in\mathbb{R}\big|\chi({x},t)=-1\big\},\qquad\Omega^+({t})\xlongequal{\mathrm{def}} \big\{x\in\mathbb{R}\big|\chi({x},{t})=1\big\},
 \end{equation*}
and
\begin{equation*}
  \Gamma({t})\xlongequal{\mathrm{def}}  \mathbb{R}\backslash\big\{\Omega^-(t)\cup\Omega^+(t)\big\},\qquad \mathrm{meas}\Gamma(t)=0.
\end{equation*}

\vskip 0.3cm
The sharp interface limit of the immiscible two-phase flow has been an important and challenging problem in recent years, and it is becoming increasingly visible, see Wang-Wang \cite{WW-2007}, Xu-Di-Yu \cite{XDY-2018}, Witterstein \cite{Witterstein2010} and the references therein. It is well known that, when the interface thickness tends to zero, in addition to the discontinuity of the phase field on both sides of the two-phase flow interface, the continuity of density and velocity of the flow may also be disrupted, which prevent us to solve the problem of singular limit by using the traditional analytical tools, and new methods are needed to deal with this open problem. 

Before analyzing this singular limit problem, let us give a briefly review of the relevant works on the well-posedness and asymptotic stability for the system \eqref{original NSAC}. The global existence of renormalized finite energy weak solutions of compressible NSAC system \eqref{original NSAC} in 3-D with no-slip boundary condition is established  by Feireisl-Petzeltov$\acute{a}$-Rocca-Schimperna \cite{FPRS2010} for the adiabatic exponent of pressure $\gamma>6$
, and this existence result was then generalized  to $\gamma>2$ by Chen-Wen-Zhu \cite{cwz2019}. The existence of the local strong solutions in high dimensional for NSAC system \eqref{original NSAC} is established  by Kotschote \cite{Kotschote2012},  Chen-He-Huang-Shi \cite{CHHS2022-1}, etc. For the global existence of the strong solutions in 1-dimensional, refer to Ding-Li-Luo \cite{DLL2013} for the Cauchy problem, Ding-Li-Tang \cite{dlt2019} for the free boundary problem, Chen-Guo \cite{CG2017} and Li-Yan-Ding-Chen\cite{lydc2023} for  the problem of the initial density containing vacuum, Chen-He-Huang-Shi \cite{CHHS2021,CHHS2022} for the problem  that the coefficient of heat conduction depends on temperature, and Yan-Ding-Li \cite{ydl2022} for the problem that the  viscosity depends on phase field. Moreover, Chen-Hong-Shi \cite{CHS-2021,CHS-2023}, Chen-Li-Tang\cite{clt2021-non-nsac} and Zhao \cite{zhao2019} presented the global existence and uniqueness of  strong solution for the Cauchy problem with small initial perturbations in 3-dimensional, and also established the asymptotic stability and decay estimates of the solutions. When the initial conditions are near the rarefaction wave, contact wave or stationary solution, the stabilities of the solutions for the Cauchy problem of compressible NSAC system with some certain restrictions on initial phase field are investigated by Yin-Zhu \cite{yz2019}, Luo-Yin-Zhu \cite{lyz2018}, Luo-Yin-Zhu \cite{lyz2020}, etc.

However, very little work has been done on the sharp interface limit for compressible NSAC system \eqref{original NSAC}.  The results that have been obtained so far are, specifically, for incompressible immiscible two-phase flow, the sharp interface limit is formally presented in Wang-Wang \cite{WW-2007}, Xu-Di-Yu \cite{XDY-2018} by means of the method of matched asymptotic expansion, they showed that the leading order problem satisfies the incompressible Navier-Stokes equations with the interface being a free boundary. Abels-Fei \cite{af2022} and Jiang-Su-Xie\cite{jsx2022} proved convergence of the solutions of the incompressible NSAC system in a two-dimensional bounded domain to solutions of a sharp interface model respectively. Hensel-Liu \cite{hl2022} showed the sharp interface limit of the incompressible NSAC system in a bounded region in $\mathbb{R}^2$ and $\mathbb{R}^3$. Recently, for the nonhomogeneous incompressible NSAC system, Jiang-Su-Xie\cite{jsx2023} investigated the sharp interface limit for inhomogeneous incompressible NSAC system in a bounded domain via a relative energy method. As for compressible immiscible two-phase
flow, the density, velocity and phase field are coupled together nonlinearly,  this makes the asymptotic expansion matching method, which is effective for the incompressible fluids, difficult to apply to the rigorous analysis of the compressible fluids.  Furthermore, unlike the incompressible fluids, considering that various basic waves, such as shock waves, rarefaction waves and contact discontinuity, will appear in the compressible immiscible two-phase flow,  the interactions between these waves and the interface,  and the interactions between the interface can be seen everywhere and are complex and important and extremely difficult to analyze. 
Witterstein \cite{Witterstein2010} formally generalized the results of Wang-Wang \cite{WW-2007} to compressible NSAC system \eqref{original NSAC},   but the arguments are only formal and lack rigorous analysis.

From a mathematical point of view, the diffused interface model is equivalent to an approximation of the ideal free boundary model, which avoids the jump of the phase field near the interface, especially makes it easy to construct a suitable numerical computation format for numerical simulation. Nevertheless, in terms of physical observations, the corresponding physical scale of the diffusion interface is actually  so small that even for the numerical computation format,  it exceeds the computing capacity of the current computers, so that the scale has to be artificially enlarged to obtain the approximate numerical simulation results, and whether it is really approximate, it needs to be strictly explained on the analysis.  Therefore, the in-depth analysis of sharp interface limit is extremely important both from the point of view of numerical simulation and mathematical physics.

In view of the above analysis, this paper starts with the analysis of interface interaction. As for the interaction of the interface, our first concern at present is the interaction of two codirectional interfaces chasing each other at the velocity of the shock wave, that is, when the front and back interfaces moving in the same direction with different shock velocity interact with each other, there's reason to guess, a forward shock wave and a backward rarefaction wave will be generated, and the newly generated rarefaction wave tail will produce a complicated wave-shift effect on the forward shock wave. Further, when the interface thickness approaches zero, these so-called viscous shock waves, smooth approximate rarefaction waves, and smooth phase fields caused by diffusion interfaces $L^\infty(\mathbb{R})$-norm converge to corresponding entropy solutions and jump discontinuous phase fields caused by phase separation interfaces, respectively.

From the physical intuition, we believe that the physical state of the sharp interface limit for immiscible two-phase flow, should be the same as the physical state after a long enough time, this is the fundamental reason for our above speculation, so we think that when the thickness of the two-phase flow interface approaches zero, the velocity field and density should converge to the corresponding entropy solution, and  the phase field converges to $\pm1$. Therefore, the following corresponding isentropic Navier-Stokes system in the Lagrangian coordinates is considered:
 \begin{equation}\label{NS-lagrange}
\left\{\begin{array}{llll}
\displaystyle v_{t}-u_{x}=0,\\
\displaystyle u_{t}+p(v)_{x}=\epsilon\Big(\frac{u_{x}}{v}\Big)_{x},
\end{array}\right.
\end{equation}
and the inviscid Navier-Stokes equations associated with the system \eqref{NS-lagrange} is as follows
 \begin{equation}\label{Riemann problem for Euler equation}
\left\{\begin{array}{llll}
\displaystyle  V_{t}-U_{x}=0,\\
\displaystyle  U_{t}+p(V)_x=0, \\
\end{array}\right.
\end{equation}
 the eigenvalues of  system (\ref{Riemann problem for Euler equation})  are
\begin{eqnarray}
\lambda_1(V)=-\sqrt{-p'(V)}=-\sqrt{\gamma \frac{p(V)}{V}}<0,\quad \lambda_2(V)=\sqrt{-p'(V)}=\sqrt{\gamma \frac{p(V)}{V}}>0.
\end{eqnarray}

\begin{figure}
\centering
\subfigure[entropy solution $\mathcal{V}$, $\mathcal{U}$]{
\begin{minipage}[t]{0.45\textwidth}
\centering
\includegraphics[width=2.8 in,height=2.2 in]{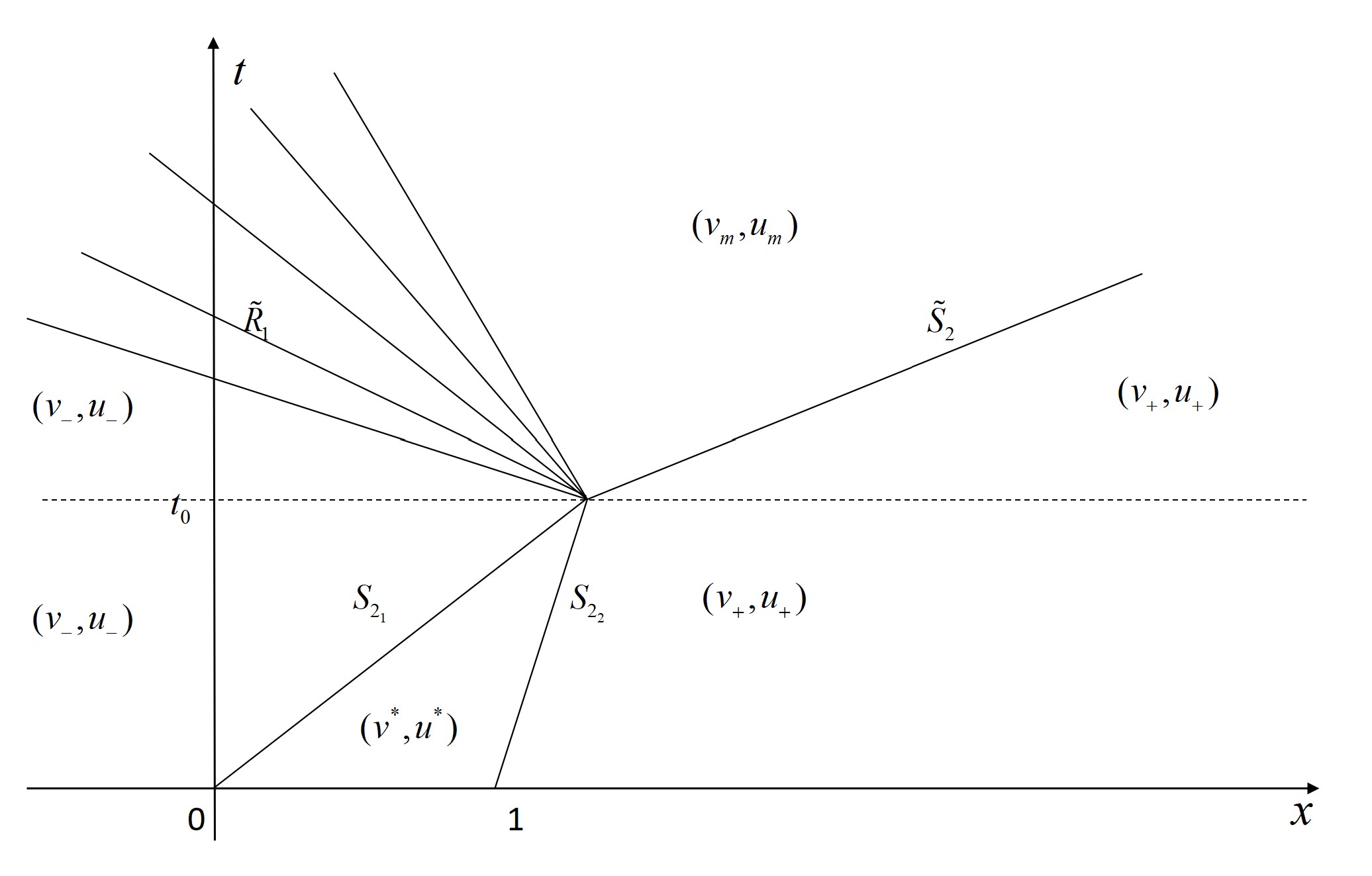}
\end{minipage}}
\hfill
\subfigure[phase plane  of  $\mathcal{V}$, $\mathcal{U}$]{
\begin{minipage}[t]{0.45\textwidth}
\centering
\includegraphics[width=2.8 in,height=2.2 in]{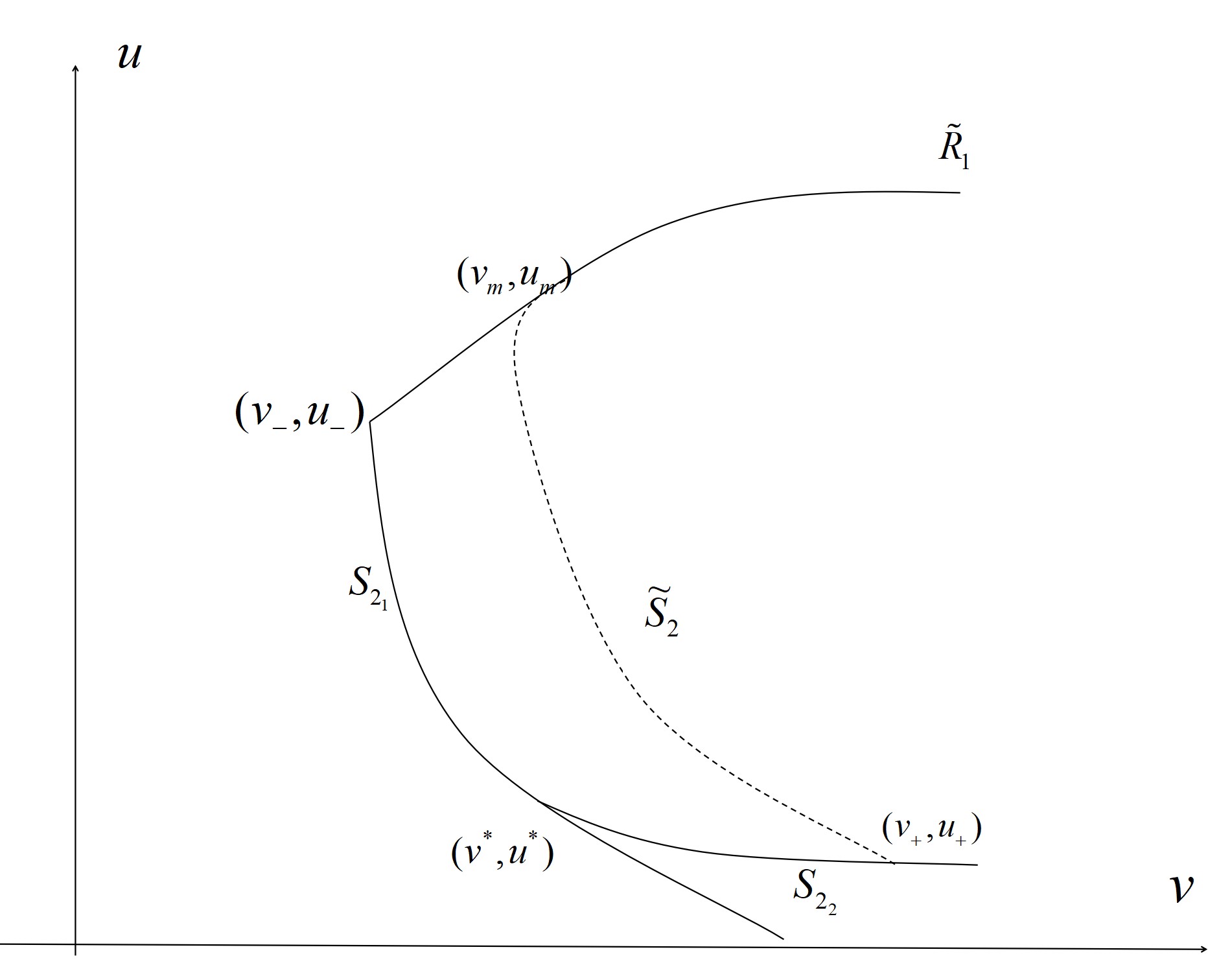}
\end{minipage}}
\centering
\caption{Figure of the interaction of shock waves in the same group}
\label{model2figure1}
\end{figure}

Given two incoming shocks denoted by $S_{2_1},S_{2_2}$ which all belongs to the second family shock waves. $S_{2_2}$ is the 2-shock wave which connects $(v^*,u^*)$ as the left state
and $(v_+,u_+)$ as the right state with the speed $s_{2_2}>0$, while $S_{2_1}$ is the another 2-shock which connects $(v_-,u_-)$ as the left state
 and $(v^*,u^*)$ as the right state with the speed $s_{2_1}>0$, and $s_{2_1}>s_{2_2}$. Because the shock wave $S_{2_1}$ propagates faster than shock wave $S_{2_2}$,
the two shocks must interact at some point $Q(x_0,t_0)$ with $t_0>0$. After the interaction, two nonlinear waves will be generated in opposite directions, that is,  the outgoing rarefaction wave $\tilde{R}_1$ and the outgoing shock wave $\tilde{S}_2$ are generated. In this case, the intermediate state
between $\tilde{R}_1$ and $\tilde{S}_2$ is $(v_m,u_m)$.  It is important to note that the  intermediate states $(v^*,u^*)$ and $(v_m,u_m)$ mentioned here are determined by \eqref{R-H condition of s21}-\eqref{R-H condition of S22} and \eqref{23}, \eqref{solution after t0}.
Without loss of generality, we assume the initial data of the Euler system (\ref{Riemann problem for Euler equation}) is
\begin{eqnarray} \label{initial for Euler}
(V,U)(x,0)=
\left\{ \begin{aligned}
& (v_-,u_-),\,\,x<0,  \\
& (v^*,u^*), \,\,0<x<1, \\
& (v_+,u_+), \,\,x>1,
\end{aligned} \right.
\end{eqnarray}
where $(v^*,u^*)$  satisfies
  \begin{equation}\label{v* and u*}
   u_- >u^*>u_+, \ \ v_-<v^*<v_+,
 \end{equation}
and is determined by the following R-H conditions \eqref{R-H condition of s21}-\eqref{R-H condition of S22}
\begin{eqnarray}\label{R-H condition of s21}
\left\{ \begin{aligned}
& -s_{2_1}(v^*-v_-)-(u^*-u_-)=0,  \\
& -s_{2_1}(u^*-u_-)+\big(p(v^*)-p(v_-)\big)=0,
\end{aligned} \right.
\end{eqnarray}
\begin{eqnarray}\label{R-H condition of S22}
\left\{ \begin{aligned}
& -s_{2_2}(v_+-v^*)-(u_+-u^*)=0,  \\
& -s_{2_2}(u_+-u^*)+\big(p(v_+)-p(v^*)\big)=0,
\end{aligned} \right.
\end{eqnarray}
where $s_{2_1}$ and $s_{2_2}$ represent the corresponding shock velocity, and satisfy the entropy condition
\begin{eqnarray}
0<\lambda_2(v_+,u_+)<s_{2_2}<\lambda_2(v^*,u^*)<s_{2_1}<\lambda_2(v_-,u_-).
\end{eqnarray}
So far, the 2-shock wave curves  $S_{2_1}$  in a suitable neighbourhood $O_{(v_{-}, u_{-})}$ of $(v_{-}, u_{-})$ and $S_{2_2}$ in a suitable neighbourhood $O_{(v^*,u^*)}$ of $(v^*,u^*)$ are defined as
\begin{eqnarray}\label{S}
&&S_{2_1}(v_-,u_-) \xlongequal{\mathrm{def}} \Big\{(v,u)\in O_{(v_-,u_-)}\Big| u=u_--s_{2_1}(v-v_-),v>v_-\Big\},\\
&& S_{2_2}(v^*,u^*) \xlongequal{\mathrm{def}} \Big\{(v,u)\in O_{(v^*,u^*)}\Big| u=u^*-s_{2_2}(v-v^*),v>v^*\Big\}.
\end{eqnarray}
By a simple calculation, the two incoming shocks hit at the point
\begin{equation}\label{Time of collision}
 (x_0,t_0)=\big(\frac{s_{2_1}}{s_{2_1}-s_{2_2}},\frac1{s_{2_1}-s_{2_2}}\big),
\end{equation}
and thus, the unique entropy solution $(\mathcal{V},\mathcal{U})$ of (\ref{Riemann problem for Euler equation}),(\ref{initial for Euler}) before the interaction time $t_0$ can be given as follows
\begin{eqnarray}\label{solution before t0}
(\mathcal{V},\mathcal{U})(x,t) \xlongequal{\mathrm{def}} 
\left\{ \begin{aligned}
& (v_-,u_-),\,\,x<s_{2_1}t, \,\,t\leq t_0 , \\
& (v^*,u^*), \,\,s_{2_1}t<x<s_{2_2}t+1, \,\,t\leq t_0, \\
& (v_+,u_+), \,\,x>s_{2_2}t+1, \,\,t\leq t_0.
\end{aligned} \right.
\end{eqnarray}
Further, after the interacting time $t_0$, it is equivalent to resolve the Riemann problem of system (\ref{Riemann problem for Euler equation}) again with the initial data at $t=t_0$ below
\begin{eqnarray}
(V,U)(x,t_0)=
\left\{ \begin{aligned}
& (v_-,u_-),\,\,x<x_0,  \\
& (v_+,u_+), \,\,x>x_0.
\end{aligned} \right.
\end{eqnarray}
The outgoing  1-rarefaction wave curve $\tilde{R}_1\left(v_{-}, u_{-}\right)$ in a suitable neighbourhood $O_{(v_{-}, u_{-})}$ of $(v_{-}, u_{-})$ is defined as
\begin{equation*}
\tilde{R}_{1}\big(v_{-}, u_{-}\big) \xlongequal{\mathrm{def}} \Big\{(v, u) \in O_{(v_-, u_-)} \Big| \ u=u_{-}-\int_{v_{-}}^{v} \lambda_1(s) d s, \ v>v_{-}\Big\},
\end{equation*}
and the corresponding Lipschitz continuous 1-rarefaction wave $(v^r, u^r)(x,t)=(v^r, u^r)(\frac{x-x_0}{t-t_0})$ of the Riemann problem \eqref{Riemann problem for Euler equation} is given explicitly by
\begin{equation}
 \lambda_1\big(v^r)=w^r, \qquad
  u^r=u_{-}-\int_{v_{-}}^{v^r} \lambda_1(s)ds,  \label{23}
\end{equation}
where $w^r(x, t)$ has the following expression
\begin{equation}
w^r(x, t)=\left\{\begin{array}{ll}
\displaystyle \lambda_1(v_{-}), & x-x_0<w_{-} (t-t_0), \\
\displaystyle \frac{x-x_0}{t-t_0}, & w_{-} (t-t_0) \leq x-x_0 \leq w_m (t-t_0), \\
\displaystyle \lambda_1(v_m), & x-x_0>w_m (t-t_0).
\end{array}\right.
\end{equation}
$(v_m,u_m)$ is the intermediate state between $\tilde{R}_1$ and $\tilde{S}_2$ can be determined by \eqref{23}, \eqref{R-H condition for shock wave after interation}.
To be clear here, $w^r(x, t)$ is actually satisfies the following Riemann problem of inviscid Burgers equation
\begin{equation}
\left\{\begin{array}{l}
w_t+w w_x=0, \\
w(x,t_0)=w_0^r(x)= \begin{cases}w_{-}=\lambda_1(v_{-}), & x<x_0, \\
w_m=\lambda_1(v_m), & x>x_0.\label{Burgers}\end{cases}
\end{array}\right.
\end{equation}
 On the other hand, the outgoing  2-shock wave  curve $\tilde{S}_{2}$  in a suitable neighbourhood $O_{(v_{m}, u_{m})}$ of $(v_{m}, u_{m})$ is defined as
\begin{equation}\label{tilde S}
\tilde{S}_{2}(v_m,u_m) \xlongequal{\mathrm{def}} \Big\{(v,u)\in O_{(v_m,u_m)}\Big| u=u_m-\tilde{s}_{2}(v-v_m),\lambda_2(v,u)<\lambda_2(v_m,u_m)\Big\},
\end{equation}
where the shock speed $\tilde{s}_2$ satisfies the following R-H condition
\begin{eqnarray}\label{R-H condition for shock wave after interation}
\left\{ \begin{aligned}
& -\tilde{s}_2(v_+-v_m)-(u_+-u_m)=0,  \\
& -\tilde{s}_2(u_+-u_m)+\big(p(v_+)-p(v_m)\big)=0,  \\
\end{aligned} \right.
\end{eqnarray}
and the entropy condition
\begin{eqnarray}\label{lax entropy condition after interation}
0<\lambda_2(v_+,u_+)<\tilde{s}_2<\lambda_2(v_m,u_m),
\end{eqnarray}
$\tilde{s}_2$ is uniquely determined by (\ref{R-H condition for shock wave after interation}) and satisfies
\begin{eqnarray} \label{the value of a and b}
\tilde{s}_2=\sqrt{\frac{p(v_m)-p(v_+)}{v_+-v_m}},\ b_1=p(v_+)+\tilde{s}_2^2v_+=p(v_m)+\tilde{s}_2^2v_m,\ b_2=-(\tilde{s}_2v_{+}+u_+).\end{eqnarray}
Therefore,
the entropy solution $(\mathcal{V},\mathcal{U})$ after the interaction time $t=t_0$ can be given as follows:
\begin{equation} \label{solution after t0}
(\mathcal{V},\mathcal{U})(x,t) \xlongequal{\mathrm{def}} 
\left\{\begin{array}{ll}
\displaystyle (v_-,u_-), & x-x_0\leq\lambda_1(v_-)(t-t_0),   \\
\displaystyle \Big(\lambda_1^{-1}(w^{r}),u_--\int_{v_-}^{\lambda_1^{-1}(w^{r})}\lambda_1(s)ds\Big), &
\displaystyle \lambda_1(v_-)<\frac{x-x_0}{t-t_0}\leq\lambda_1(v_m),  \\
\displaystyle (v_m,u_m), &\displaystyle \lambda_1(v_m)<\frac{x-x_0}{t-t_0}\leq\tilde{s}_2,\\
\displaystyle (v_+,u_+), &\displaystyle x-x_0>\tilde{s}_2(t-t_0),
\end{array} \right.
\end{equation}
 Setting
\begin{equation}\label{initial wave strength}
\delta_1=\max\big\{|v^*-v_-|,|u^*-u_-|\big\},\qquad\delta_2=\max\big\{|v_+-v^*|,|u_+-u^*|\big\} ,
\end{equation}
be the strengths for the two incoming shocks of the same family, and letting
\begin{equation}
\tilde{\delta}_1=\max\big\{|v_m-v_-|,|u_m-u_-|\big\},\qquad  \tilde{\delta}_2=\max\big\{|v^+-v_m|,|u^+-u_m|\big\},\end{equation}
be the wave strengths of the outgoing waves,   and the relations of the wave strengths between the incoming and the outgoing waves can be provided by the following formula (see Smoller \cite{S94} Proposition 19.1)
\begin{equation}\label{Strength after interaction}
\tilde{\delta}_1=\delta_1+O(1)\delta_1\delta_2 ,\ \ \tilde{\delta}_2=\delta_2+O(1)\delta_1\delta_2.
\end{equation}

\

\noindent\textbf{\normalsize Notations.} We denote by $C$ and $c$ the positive generic constants  without confusion throughout this paper,
  $L^2(\mathbb{R})$ denotes the space of Lebesgue measurable functions on $\mathbb{R}$ which are square integrable, with the norm $\|f\|=(\int_{\mathbb{R}}|f|^2)^{\frac{1}{2}}$.
 $H^l(\mathbb{R})(l\geq0)$ denotes the Sobolev space of $L^2$-functions $f$ on $\mathbb{R}$ whose derivatives $\partial^j_x f,  j=1,\cdots$ are $L^2$
 functions too, with the norm
$ \|f\|_{H^l(\mathbb{R})}=(\sum_{j=0}^l\|\partial^j_x f\|^2)^{\frac{1}{2}}$.

\

So far, we will give the main results about the large time behavior of the solution and the singular limit of diffusion interface thickness respectively in the following theorems.

\begin{theorem}(Asymptotic stability)
Suppose that $(v^*, u^*)\in S_{2_1}(v_-, u_-)$, $(v_+, u_+)\in S_{2_2}(v^{*}, u^{*})$, $(v_m, u_m)\in \tilde{R}_1(v_-,u_-)$ and $(v_+, u_+)\in \tilde{S}_{2}(v_m, u_m)$.
Let $v_0, u_0, \chi_{0}$ be the initial data of the compressible NSAC system \eqref{nsac-lagrange} satisfies the well-prepared initial data (\ref{initial data of N-S}) in Section 2, $\chi_0\in H^2(\mathbb{R})$,  $(\bar V,\bar U)$ is the superposition of two incoming viscous shock waves given by \eqref{the composite wave profile before interaction},  $(v^r,u^r)$ and $(\tilde V^{\tilde{S}_2}, \tilde U^{\tilde{S}_2})$ are 
the outgoing rarefaction wave and viscous shock wave given by \eqref{23} and  \eqref{smooth-wave} respectively,  for fixed $\epsilon$, there exist  positive constants $\delta_0$ and $M_0$ such that
 \begin{equation}
\delta_1+\delta_2\leq\delta_0,\qquad\| \chi_{0x}\|_{{2}}+\|\chi_0^2-1\|\leq M_0, \end{equation}
 then the Cauchy problem of NSAC system \eqref{nsac-lagrange}-\eqref{initial data} admits a unique global-in-time solution $(v, u,\chi)$ satisfying
 $$\chi^2-1 \in C\left([0,+\infty ); H^2(\mathbb{R})\right)\cap L^2\left([0,+\infty) ; H^3(\mathbb{R})\right),\ \ -1\leq\chi\leq1,$$
and before the interaction time $t_0$, 
\begin{eqnarray}
&& (v-\bar{V},u-\bar{U}) \in C\big([0,t_0) ; H^2(\mathbb{R})\big),\notag\\
&&(v-\bar V) \in L^2\left([0,t_0); H^2(\mathbb{R})\right),\ \ (u-\bar U) \in L^2\left([0,t_0); H^3(\mathbb{R})\right),\notag
\end{eqnarray}
after the interaction time  $t_0$,
\begin{eqnarray}
&& v(x,t)-\Big( v^{r_1}(\frac{x-x_0}{t-t_0})+\tilde V^{\tilde{S}_2}\big(x-x_0-\tilde{s}_2(t-t_0)-\mathbf{X}(t-t_0)\big)-v_m\Big) \in C\big([t_0,+\infty) ; H^1(\mathbb{R})\big), \notag\\
&& u(x,t)-\Big( u^{r_1}(\frac{x-x_0}{t-t_0})+\tilde U^{\tilde{S}_2}\big(x-x_0-\tilde{s}_2(t-t_0)-\mathbf{X}(t-t_0)\big)-u_m\Big) \in C\big([t_0,+\infty) ; H^1(\mathbb{R})\big), \notag\\
&& u_{xx}-\tilde U^{\tilde{S}_2}_{xx}\big(x-x_0-\tilde{s}_2(t-t_0)-\mathbf{X}(t-t_0)\big) \in L^2\left([t_0,+\infty) ; L^2(\mathbb{R})\right),\notag
\end{eqnarray}
where $\mathbf{X}(t)$ is  determined by \eqref{Shift X}, In addition, it holds that
\begin{eqnarray}
&&    \lim_{t\rightarrow +\infty}\sup _{x \in \mathbb{R}}\Big|v(x,t)-\Big( v^{r_1}(\frac{x-x_0}{t-t_0})+\tilde V^{\tilde{S}_2}\big(x-x_0-\tilde{s}_2(t-t_0)-\mathbf{X}(t-t_0)\big)-v_m\Big)\Big| = 0,\notag\\
&&   \lim_{t\rightarrow +\infty}\sup _{x \in \mathbb{R}}\Big|u(x,t)-\Big( u^{r_1}(\frac{x-x_0}{t-t_0})+\tilde U^{\tilde{S}_2}\big(x-x_0-\tilde{s}_2(t-t_0)-\mathbf{X}(t-t_0)\big)-u_m\Big)\Big| = 0,\notag\\
&& \lim_{t\rightarrow+\infty}\sup _{x \in \mathbb{R}}\big|\chi^2-1\big|=0,\ \ \lim_{t\rightarrow+\infty}\big|\frac{d}{dt}\mathbf{X}(t)\big|=0.\notag
\end{eqnarray}
\end{theorem}

\begin{remark} Theorem 1.1 implies that, when two interfaces move in the same direction with different shock wave velocities, if the speed of the rear interface is faster than the speed of the front interface, the rear interface will catch up with the front interface and interact with each other in a finite time, so that it will evolve into a forward interface moving at the shock wave velocity and a backward interface moving at the rarefaction wave velocity.
\end{remark}

\begin{remark}It is not necessary to take the initial conditions directly as the composite viscous waves \eqref{the composite viscous wave profile before interaction}  and the smoothed approximate rarefaction wave constructed in Section 2. In fact, Theorem 1.1 also holds for small perturbations of the initial condition. 
\end{remark}

\begin{theorem} (Singular limit of diffusion interface thickness) \label{main thm}Suppose that $(v_{\pm}, u_{\pm})$, $(v^*, u^*)$ and $(v_m, u_m)$ satisfy $(v^*, u^*)\in S_{2_1}(v_-, u_-)$, $(v_+, u_+)\in S_{2_2}(v^{*}, u^{*})$, $(v_m, u_m)\in \tilde{R}_1(v_-,u_-)$ and $(v_+, u_+)\in \tilde{S}_{2}(v_m, u_m)$.
Let $v_0, u_0, \chi_{0}$ be the initial data of the compressible NSAC system \eqref{nsac-lagrange} satisfies the well-prepared initial data (\ref{initial data of N-S}) in Section 2, $\chi_0\in H^2(\mathbb{R})$,  $(\mathcal{V},\mathcal{U})$ be the entropy solution of system (\ref{Riemann problem for Euler equation}) given by (\ref{solution before t0}) and (\ref{solution after t0}),  there exist  positive constants $\delta_0$ and $M_0$ such that
\begin{equation}
  \delta_1+\delta_2\leq\delta_0,\qquad \|\chi_0^2-1\|\leq M_0,
\end{equation}
then the system (\ref{nsac-lagrange})-\eqref{initial data} admits a family
smooth solution $(v^{\epsilon},u^{\epsilon},\chi^{\epsilon})$ for any $\epsilon>0$. Moreover,  before the interaction time $t=t_0$, it holds that
\begin{equation}\label{sharp interface limit before interaction}
\lim_{\epsilon\rightarrow0^+}\big\|v^{\epsilon}-\mathcal{V},u^{\epsilon}-\mathcal{U},{(\chi^\epsilon)}^{2}-1\big\|_{L^\infty(\Sigma_\varsigma)}=0,
\end{equation}
where $\Sigma_\varsigma=\big\{(x,t)\big|  |x-s_{2_1}t|\geq \varsigma, \, |x-s_{2_2}t-1|\geq \varsigma,\, 0\leq t\leq t_0-\varsigma\big\}$ for any positive constant $\varsigma>0$.
After the interaction time $t_0$, it holds that
\begin{equation}\label{sharp interface limit after interfaction}
\lim_{\epsilon\rightarrow0^+}\big\|v^{\epsilon}-\mathcal{V},u^{\epsilon}-\mathcal{U},{(\chi^\epsilon)}^{2}-1\big\|_{L^\infty(\widetilde{\Sigma}_{\tilde{\varsigma}})}=0,\end{equation}
where
$\tilde{\Sigma}_{\tilde{\varsigma}}=\big\{(x,t)\big| |(x-x_0)-\tilde{s}_2(t-t_0)|\geq \tilde \varsigma,\  \tilde {\varsigma}\leq t-t_0\big\}$
for any positive constant $\tilde \varsigma>0$.
\end{theorem}

\begin{remark} The
Navier-Stokes equations can be regarded as a special case of the NSAC system, and the conclusion of Theorem 1.2 is also applicable for the problem of vanishing viscosity for isentropic Navier-Stokes equations, and in this sense,  Theorem 1.2 answers and generalizes the open problem raised in Huang-Wang-Wang-Yang \cite{HWWY-2015} about the single-phase flow. 
\end{remark}
\begin{remark} The conclusion of Theorem 1.2 can be seen as a special case of sharp interface limit for Navier-Stokes/Allen-Cahn system. Under the condition of Theorem 1.2, when the thickness of the diffusion interface tends to zero, the NSAC system \eqref{nsac-lagrange} tends to the $p$-system \eqref{NSAC-Euler-equation}, and the diffusion interface evolves into a free boundary.
\end{remark}

\begin{remark}  Theorem 1.2 is concerned with the $L^{\infty}$-norm which gives a detailed description on the sharp interface limit, the range of validity of this $L^{\infty}$-norm sharp interface limit requires the removal of the neighborhood near the interface and the interaction, and this constraint is necessary, although the scale $\varsigma,\tilde \varsigma$ can be small enough.
  \end{remark}

\begin{remark}
The assumption  $\displaystyle \lim_{x\rightarrow\pm\infty}\chi_0=\pm1$ in \eqref{initial data} means that the initial phase field of the immiscible two-phase flow is in a different phase state at a distance, this assumption can also be set more generally as  $\displaystyle \lim_{x\rightarrow\pm\infty}\chi^2_0=1$. 
\end{remark}

Here we make some comments on the key steps in the proof. First, scaling transformation is introduced to transformed the sharp interface limit  into the problem of  large time behavior for the solution of NSAC system.
Second, before the interacting shock waves of the same family, the anti-derivative technique is used to study the stability of viscous shocks. However, due to the effect of phase field in compressible immiscible two-phase flow, it is not possible to use anti-derivative technique directly on NSAC system. Therefore,  the  anti-derivative method is modified here, that is,  the antiderivative only be used for the mass equation and momentum equation, while keeping the phase field equation unchanged, see \eqref{antiferivative technique} and \eqref{equation of composite shock wave before interaction}. In addition, in order to get more information about the density before interaction,  the weighted energy estimation method is used (see \eqref{the lower estimate} and \eqref{quadratic form}).
Finally,  after interaction, as the shock waves of same family catch up with each other and interacting, a viscous shock wave in the same direction and an inverted rarefaction wave are generated. In order to continue to be able to use the energy estimation method, we need to derive the estimates separately before and after the interaction time and connect the incoming and outgoing on these nonlinear waves, see the energy inequality at $t=t_0$ \eqref{key Lemma of the estimate for initial data}. Moreover, since the traditional energy method for dealing with the stability of rarefaction waves does not match the anti-derivative method, to overcome this difficulty,  the time-dependent shift $\mathbf{X}(t)$ is introduced to the viscous shock wave (see \eqref{Shift X}), associated to the BD entropy (see \eqref{effective velocity}), and combined with the weight function technique (see \eqref{412}). This  new method is used by Kang-Vasseur-Wang \cite{kvw2021}  to discuss the asymptotic stability of composite waves, and is called as  \textit{a}-contraction method.

The outline of this paper is organized as follows. In Section 2, The scaling transformation is applied to the NSAC system \eqref{nsac-lagrange},  the approximate solutions (such as viscous shock waves, smooth approximation of rarefaction waves) of entropy solution are constructed, some basic properties (such as decay estimates etc.) of these approximate nonlinear waves are given, and the definition of relative quantities  and the related properties are also given here.  In Section 3,  the properties for the solutions of NSAC system \eqref{nsac-lagrange} is established before the interaction of the shocks. In Section 4, the desired energy estimates for the solutions of NSAC system \eqref{nsac-lagrange} are provided after the interaction time.

\section{Preliminaries}
\indent\qquad
This section is concerned with the  estimates for the viscous shock waves and rarefaction waves.  Without loss of generality,
we rewrite the  system (\ref{nsac-lagrange}) by  scaling method, letting
\begin{equation}\label{scaling method}
y=\frac{x-x_{0}}{\epsilon},\ \ \tau=\frac{t-t_0}{\epsilon},
\end{equation}
and introducing the following new variable
\begin{equation}\label{Variable substitution 1}
 \omega=\chi^2,
\end{equation}
then the original system (\ref{nsac-lagrange}) is reformulated into the form
\begin{equation}\label{NSAC-scaling}
\left\{\begin{array}{llll}
\displaystyle v_{\tau}-u_{y}=0,\\
\displaystyle u_{\tau}+p(v)_{y}=\Big(\frac{u_{y}}{v}\Big)_{y}-\frac{1}{8}\Big(\frac{\omega_{y}^{2}}{\omega v^{2}}\Big)_{y}, \\
\displaystyle \omega_{\tau}=-2v(\omega-1)\omega+v\Big(\frac{\omega_{y}}{v}\Big)_{y}-\frac{\omega_{y}^{2}}{2\omega}.
\end{array}\right.
\end{equation}
Now we begin to construct the viscous shock waves before $\tau=0$, that is,  considering the Navier-Stokes equations \eqref{NS-lagrange}. The left 2-viscous incoming shock wave $(V^{S_{2_1}},U^{S_{2_1}})(y-s_{2_1}\tau)$  of the system (\ref{NS-lagrange})  which connects $(v_-,u_-)$ on the left and $(v^*,u^*)$ on the right exists uniquely up to a shift and satisfies
\begin{eqnarray} \label{21 shock equation}
\left\{ \begin{aligned}
& -s_{2_1}V_\xi^{S_{2_1}}-U_\xi^{S_{2_1}}=0,  \\
& -s_{2_1}U_\xi^{S_{2_1}}+p_\xi(V^{S_{2_1}})=\Big(\frac{U_\xi^{S_{2_1}}}{V^{S_{2_1}}}\Big)_\xi, \\
&\lim_{\xi\rightarrow-\infty}(V^{S_{2_1}},U^{S_{2_1}})=(v_-,u_-),\ \ \lim_{\xi\rightarrow+\infty}(V^{S_{2_1}},U^{S_{2_1}})=(v^{*},u^{*}),
\end{aligned} \right.
\end{eqnarray}
where $\cdot_\xi=\frac d{d\xi}$, and $\xi=y-s_{2_1}\tau$. To fix the viscous shock wave, setting
\begin{equation}
V^{S_{2_1}}(0)=\frac{v_-+v^*}2,
\end{equation}
then by using $(\ref{21 shock equation})_1$ and (\ref{R-H condition of s21}), one has
\begin{eqnarray}
 U^{S_{2_1}}(0)=u_--s_{2_1}(V^{S_{2_1}}(0)-v_-)=u_--\frac{1}{2}s_{2_1}(v^*-v_-)=\frac{u_-+u^*}2.
\end{eqnarray}
Similarly, the right 2-viscous incoming shock wave  which connects $(v^*,u^*)$ on the left and $(v_+,u_+)$ on the right satisfies
\begin{eqnarray} \label{22 shock equation}
\left\{ \begin{aligned}
& -s_{2_2}V_\xi^{S_{2_2}}-U_\xi^{S_{2_2}}=0,  \\
& -s_{2_2}U_\xi^{S_{2_2}}+p(V^{S_{2_2}})_\xi=\Big(\frac{U_\xi^{S_{2_2}}}{V^{S_{2_2}}}\Big)_\xi, \\
&\lim_{\xi\rightarrow-\infty}(V^{S_{2_2}},U^{S_{2_2}})=(v^*,u^*),\ \ \lim_{\xi\rightarrow+\infty}(V^{S_{2_2}},U^{S_{2_2}})=(v_{+},u_{+}).
\end{aligned} \right.
\end{eqnarray}
In order to determine the above viscous shock wave, the following constraint need to be given
\begin{equation}
(V^{S_{2_2}},U^{S_{2_2}})(0)=\big(\frac{v_++v^*}2,\frac{u_++u^*}2\big),
\end{equation}
where $\cdot_\xi=\frac d{d\xi}$, and $\xi=y-s_{2_2}\tau$.

Further, after the interaction of these two shock waves of the same family,  we will construct the corresponding smooth rarefaction wave and viscous shock wave, that is, for  the  1-approximate outgoing rarefaction wave $(\tilde{V}^{\tilde{R}_{1}},\tilde{U}^{\tilde{R}_{1}})(y,\tau)$  to the problem (\ref{nsac-lagrange}) which connecting $(v_-,u_-)$ on the left and $(v_m,u_m)$ on the right, it can be expressed by the following formula
 \begin{equation} \label{1-rarefaction wave}
(\tilde{V}^{\tilde{R}_1},\tilde{U}^{\tilde{R}_1})(y,\tau)=\left\{ \begin{array}{ll}
(v_-,u_-),\ & y\leq \lambda_1(v_-)\tau,   \\
\Big(
\lambda_1^{-1}(w_{\ell}^{r}),u_--\int_{v_-}^{\lambda_1^{-1}(w^r_{\ell})}\lambda_1(s)ds\Big),&\lambda_1(v_-)\tau<y\leq\lambda_1(v_m)\tau,  \\
(v_m,u_m),&\lambda_1(v_m)\tau<y\leq\tilde{s}_2\tau,
\end{array} \right.
\end{equation}
 where $w_{\ell}^{r}(y,\tau)$ is the  smooth solution of the following  Burgers equation
 \begin{eqnarray} \label{Burgers equation}
\left\{ \begin{aligned}
& w_{t}+ww_{y}=0,  \\
& w(y,0)=w_{\ell}(y)=w\big(\frac{y}{\ell}\big)=\frac{w_m+w_-}{2}+\frac{w_m-w_-}{2}\tanh\frac{y}{\ell},
\end{aligned} \right.
\end{eqnarray}
with $  w_{-}=\lambda_1(v_-)<w_m=\lambda_1(v_m)<0$,
 $\ell>0$ is a small parameter to be determined. As is well-known, 
the approximate outgoing rarefaction wave $(\tilde{V}^{\tilde{R}_1},\tilde{U}^{\tilde{R}_1})$ has the following properties (see Matsumura-Nishihara \cite{MN86}):
\begin{lemma}\label{the properties of rarefaction wave}
The approximate rarefaction waves $(\tilde{V}^{\tilde{R}_1},\tilde{U}^{\tilde{R}_1})(y,\tau)$ which be constructed in  (\ref{1-rarefaction wave}) have the following properties:

\noindent i) $\tilde{U}_y^{\tilde{R}_1}>0$, for $y\in\mathbb{R}, \tau>0$.

\noindent ii) For any $1\leq q\leq+\infty$, the following estimates holds,
\begin{eqnarray}
&&\big\|(\tilde{V}^{\tilde{R}_1},\tilde{U}^{\tilde{R}_1})_y\big\|_{L^q(\mathbb{R})}\leq C \min\big\{\tilde\delta_1\ell^{-1+\frac{1}{q}},\tilde\delta_1^{\frac 1q}
\tau^{-1+\frac{1}{q}}\big\},\notag\\
&&\big\|(\tilde{V}^{\tilde{R}_1},\tilde{U}^{\tilde{R}_1})_{yy}\big\|_{L^q(\mathbb{R})}\leq C \min\big\{\tilde\delta_1\ell^{-2+\frac{1}{q}},\ell^{-1+\frac{1}{q}}
\tau^{-1}\big\},\notag
\end{eqnarray}
where $\tilde\delta_1=|v_m-v_-|$ is the rarefaction wave strength and the positive constant $C$ depending only on $q$.

\noindent iii) If $y\geq \lambda_1(v_m)\tau$, then
\begin{eqnarray}
&&\big|(\tilde{V}^{\tilde{R}_1},\tilde{U}^{\tilde{R}_1})-(v_m,u_m)\big|\leq C\delta^{\tilde R_1} e^{-\frac{2|y-\lambda_1(v_m)\tau|}{\ell}},\notag\\
&&\big|\partial^k_y(\tilde{V}^{\tilde{R}_1},\tilde{U}^{\tilde{R}_1})-(v_m,u_m)\big|\leq \frac{C}{\ell^k}\delta^{\tilde R_1} e^{-\frac{2|y-\lambda_1(v_m)\tau|}{\ell}}, k=1,2.\notag
\end{eqnarray}
\noindent iv) There exist positive constants C and $\ell_0$ such that for $\ell\in(0,\ell_{0})$ and $\tau>0$
  \begin{eqnarray}
 &&\sup_{y\in \mathbb{R}}\big|(\tilde{V}^{\tilde{R}_1},\tilde{U}^{\tilde{R}_1})(y,\tau)-(v^{r_1},u^{r_1})(\frac{y}{\tau})\big|\leq\frac{C}{\tau}\big[\ell\ln(1+\tau)+\sigma|\ln\ell|\big],\notag\\
  &&\big\|(\tilde{V}^{\tilde{R}_1},\tilde{U}^{\tilde{R}_1})(y,\tau)-(v^{r_1},u^{r_1})(\frac{y}{\tau})\big\|\leq C\tilde\delta_1\ell.\notag
\end{eqnarray}
\end{lemma}

On the other hand, the 2-viscous outgoing shock wave which connecting $(v_m,u_m)$
on the left and $(v^+,u^+)$ on the right satisfies
\begin{equation} \label{2 shock after}
\left\{
\begin{array}{l}
\displaystyle -\tilde{s}_{2}\tilde{V}_\xi^{\tilde{S}_{2}}-\tilde{U}_\xi^{\tilde{S}_{2}}=0,  \\
\displaystyle-\tilde{s}_{2}\tilde{U}_\xi^{\tilde{S}_{2}}+p_\xi(\tilde{V}^{\tilde{S}_{2}})=\Big(\frac{\tilde{U}_\xi^{\tilde{S}_{2}}}{\tilde{V}^{\tilde{S}_{2}}}\Big)_\xi, \\
\displaystyle\lim_{\xi\rightarrow-\infty}(\tilde{V}^{\tilde{S}_{2}},\tilde{U}^{\tilde{S}_{2}})=(v_m,u_m),\quad\lim_{\xi\rightarrow+\infty}(\tilde{V}^{\tilde{S}_{2}},\tilde{U}^{\tilde{S}_{2}})=(v_{+},u_{+}),
\end{array} \right.
\end{equation}
with  the following constraint:
\begin{equation}
(\tilde{V}^{\tilde{S}_{2}},\tilde{U}^{\tilde{S}_{2}})(0)=(\frac{v_++v_m}2,\frac{u_++u_m}2),
\end{equation}
where $\cdot_\xi=\frac d{d\xi}$, and $\xi=y-\tilde{s}_{2}\tau$.
And now, the properties of the viscous shock waves which were constructed in the system \eqref{21 shock equation}, \eqref{22 shock equation}
 and  \eqref{2 shock after} are presented as following(see Matsumura-Nishihara \cite{MN85}):
\begin{lemma}\label{estimate of viscous shock}
There are positive constant $C_0$, $C$ and $c_0$ such that
\begin{eqnarray}
&&\big|(V^{S_{2_1}}-v^*,U^{S_{2_1}}-u^*)\big|<C_0\delta_1e^{-c_0\delta_1|y-s_{2_1}\tau|},\quad y> s_{2_1}\tau,\ 0>\tau\geq-t_{0}/\epsilon, \notag \\
&&\big|(V^{S_{2_1}}-v_-,U^{S_{2_1}}-u_-)\big|<C_0\delta_1e^{-c_0\delta_1|y-s_{2_1}\tau|},\ y< s_{2_1}\tau,\ 0>\tau\geq-t_{0}/\epsilon, \notag \\
&&\big|(V^{S_{2_2}}-v^*,U^{S_{2_2}}-u^*)\big|<C_0\delta_2e^{-c_0\delta_2|y-s_{2_2}\tau|},\quad y< s_{2_2}\tau,\ 0>\tau\geq-t_{0}/\epsilon, \notag \\
&&\big|(V^{S_{2_2}}-v_+,U^{S_{2_2}}-u_+)\big|<C_0\delta_2e^{-c_0\delta_2|y-s_{2_2}\tau|},\ y> s_{2_2}\tau,\ 0>\tau\geq-t_{0}/\epsilon, \notag \\
&&\big|(\tilde{V}^{\tilde{S}_{2}}-v_m,\tilde{U}^{\tilde{S}_{2}}-u_m)\big|<C_0\tilde{\delta}_2e^{-c_0\tilde{\delta}_2|y-\tilde{s}_{2}\tau|},\quad  y< \tilde{s}_{2}\tau,\ \tau\geq 0, \notag\\
&&\big|(\tilde{V}^{\tilde{S}_{2}}-v_+,\tilde{U}^{\tilde{S}_{2}}-u_+)\big|<C_0\tilde{\delta}_2e^{-c_0\tilde{\delta}_2|y-\tilde{s}_2\tau|},\quad y> \tilde{s}_{2}\tau,
\,\,\tau\geq 0, \notag \\
&&\big|(V^{S_{2_i}}_{y},U^{S_{2_i}}_{y})\big|\leq C\delta_i^2e^{-c_0\delta_1|y-s_{2_i}\tau|},i=1,2,
\ \big|(\tilde{V}^{\tilde{S}_{2}}_{y},\tilde{U}^{\tilde{S}_{2}}_{y})\big|\leq C\tilde{\delta}_2^2e^{-c_0\tilde{\delta}_2|y-\tilde{s}_2\tau|},\notag\\
&&\big|(V^{S_{2_i}}_{yy},U^{S_{2_i}}_{yy})\big|\leq C\delta_i^3e^{-c_0\delta_1|y-s_{2_i}\tau|},i=1,2,
\ \big|(\tilde{V}^{\tilde{S}_{2}}_{yy},\tilde{U}^{\tilde{S}_{2}}_{yy})\big|\leq C\tilde{\delta}_2^3e^{-c_0\tilde{\delta}_2|y-\tilde{s}_2\tau|},\notag\\
 &&v_-\leq V^{S_{2_1}}\leq v^*,v^*\leq V^{S_{2_2}}\leq v_+,v_m\leq\tilde{V}^{\tilde{S}_{2}}\leq v_+,\tilde{V}^{\tilde{S}_{2}}_{y}>0,V^{S_{2_i}}_{y}>0, i=1,2. \notag
\end{eqnarray}
\end{lemma}

Obviously, if we pull the viscous shock waves and the smooth rarefaction wave back to the coordinate system of $(x,t)$ before the scaling transformation, for any $\epsilon>0$, these corresponding  nonlinear smooth waves  can be expressed as
\begin{equation}\label{smooth-wave}
\left. \begin{array}{lll}
 \displaystyle (V_\epsilon^{S_{2_1}},U_\epsilon^{S_{2_1}})(x,t)= (V^{S_{2_1}},U^{S_{2_1}})\big(\frac{x-x_0-s_{2_1}(t-t_0)}{\epsilon}\big),\\
 \displaystyle (V_\epsilon^{S_{2_2}},U_\epsilon^{S_{2_2}})(x,t)= (V^{S_{2_2}},U^{S_{2_2}})\big(\frac{x-x_0-s_{2_2}(t-t_0)}{\epsilon}\big),\\
 \displaystyle (\tilde{V}_\epsilon^{\tilde{R}_1},\tilde{U}_\epsilon^{\tilde{R}_1})(x,t)= 
 \big(\tilde{V}^{\tilde{R}_1},\tilde{U}^{\tilde{R}_1}\big)\big(\frac{x-x_0}{\epsilon},\frac{t-t_0}{\epsilon}\big),\\
 \displaystyle  (\tilde V_\epsilon^{\tilde S_{2}},\tilde U_\epsilon^{\tilde S_{2}})(x,t)= (\tilde{V}^{\tilde{S}_{2}},\tilde{U}^{\tilde{S}_{2}})\big(\frac{x-x_0-\tilde{s}_{2}(t-t_0)}{\epsilon}\big).
  \end{array}\right.
\end{equation}
 In order to approximate the superposition of two incoming shock waves in the same family, setting
 \begin{equation}\label{the composite viscous wave profile before interaction}
 \left.
  \begin{aligned}
 (\bar{V}^\epsilon,\bar{U}^\epsilon)(x,t)\xlongequal{\mathrm{def}} &\big(V_\epsilon^{S_{2_1}},U_\epsilon^{S_{2_1}}\big)\big(x-x_0-s_{2_1}(t-t_0)\big)
\\
  &\quad+\big(V_\epsilon^{S_{2_2}},U_\epsilon^{S_{2_2}}\big)\big(x-x_0-s_{2_2}(t-t_0)\big)-(v^*,u^*).
   \end{aligned}
\right.
\end{equation}
Without loss of generality, let the above composite viscous shock waves \eqref{the composite viscous wave profile before interaction} at $t=0$ to be the initial data of $v,u$ for the system of \eqref{nsac-lagrange}, that is
\begin{equation}\label{initial data of N-S}
(v,u,\chi)\Big|_{t=0}=\big(\bar{V}^\epsilon,\bar{U}^\epsilon,\chi_0\big)\big(x,0\big)\xrightarrow {x\rightarrow\pm\infty} (v_\pm, u_\pm, \pm1).
\end{equation}
Finally, we define relative entropy as follows:
\begin{gather}
F\big(v \mid V\big)=F(v)-F(V)-F^{\prime}(V)(v-V).\label{relative entropy}
\end{gather}
where $F$ is a differentiable functional defined for $v, V>0$. Letting  
\begin{equation}\label{p-Q}
  Q(v)=\frac {v^{1-\gamma}} {1-\gamma},
\end{equation}
 $Q$ is called as the internal energy. The relative quantity inequalities commonly used in this paper are given as following (see Kang-Vasseur \cite{kv2021} and the references therein):

\begin{lemma}\label{propoties of relative quantity}
For given constants $\gamma>1$, and $M>0$, there exist constants $C, \delta_*>0$, such that the relative quantities associated to $p$ and $Q$ satisfy

\noindent 1) For any $v$, $V$, if $0<v<3 M, 0<V \leq 2 M$, it holds that
 $$|v-V|^2 \leq C Q(v \mid V),\qquad |v-V|^2 \leq C p(v \mid V).$$ 

\noindent 2) For any $v, V>\frac {M}2$, it holds that 
$$|p(v)-p(V)| \leq C|v-V| .$$

\noindent 3) $\forall 0<\delta<\delta_*$, $\forall v, V>0$, if $|p(v)-p(V)|<\delta$, $|p(V)-p(M)| <\delta$, it holds that
\begin{equation}\left. \begin{aligned}
&p\big(v \mid V\big) \leq\big(\frac{\gamma+1}{2 \gamma} \frac{1}{p(V)}+C \delta\big)|p(v)-p(V)|^2,\notag\\
& Q\big(v \mid V\big) \leq\big(\frac{p(V)^{-\frac{1}{\gamma}-1}}{2 \gamma}+C \delta\big)\big|p(v)-p(V)\big|^2,\notag\\
&Q\big(v \mid V\big) \geq \frac{p(V)^{-\frac{1}{\gamma}-1}}{2 \gamma}\mid p(v)-p(V)\big|^2-\frac{1+\gamma}{3 \gamma^2} p(V)^{-\frac{1}{\gamma}-2}\big(p(v)-p(V)\big)^3.\notag
\end{aligned}\right.\end{equation}
\end{lemma}

\section{Estimates Before the Interacting Time}
\setcounter{equation}{0}
\indent\qquad
In this section, the case before the interaction of two homologous shock waves will be considered. For $-\frac{t_0}{\epsilon}\leq \tau\leq0$, the composite wave of two approximate forward shock waves pursuing in the same direction is defined as follows
\begin{equation}\label{the composite wave profile before interaction}
\left.\begin{aligned}
  (\bar{V},\bar{U})(y,\tau)&\xlongequal{\mathrm{def}} (\bar{V}^\epsilon,\bar{U}^\epsilon)(x,t)\\
 &= \big(V^{S_{2_1}},U^{S_{2_1}}\big)\big(y-s_{2_1}\tau\big)
  +\big(V^{S_{2_2}},U^{S_{2_2}}\big)\big(y-s_{2_2}\tau\big)-(v^*,u^*),\end{aligned}\right.
\end{equation}
and it is easy to check that the profile $(\bar{V},\bar{U})$ satisfies
\begin{equation}\label{equation of approximate composite shock wave before interaction}
 \left\{\begin{array}{llll}
\displaystyle \bar{V}_\tau-\bar{U}_y=0, \\
\displaystyle \bar{U}_\tau+p(\bar{V})_y=\Big(\frac{\bar{U}_y}{\bar{V}}\Big)_y+G_y,
\end{array}\right.
\end{equation}
where
\begin{equation}\label{G}
  G=\Big(p(\bar{V})-p(V^{S_{2_1}})-p(V^{S_{2_2}})+p(v^*)\Big)-\Big(\frac{\bar{U}_y}{\bar{V}}-\frac{U^{S_{2_1}}_y}{V^{S_{2_1}}}-\frac{U^{S_{2_2}}_y}{V^{S_{2_2}}}\Big).
\end{equation}
Putting the perturbation around the superposition wave profile $(\bar{V},\bar{U},1)(y,\tau)$ by
\begin{eqnarray}
\big(\phi,\psi,\sigma\big)(y,\tau)=\big(v-\bar V,u-\bar U,\omega-1\big)(y,\tau),\ \ -\frac{t_0}{\epsilon}\leq\tau\leq0,
\end{eqnarray}
and by using the system \eqref{NSAC-scaling},  (\ref{equation of approximate composite shock wave before interaction}), one has
 \begin{equation}\label{system before interaction}
 \left\{\begin{array}{llll}
\displaystyle \phi_\tau-\psi_y=0, \\
\displaystyle \psi_\tau+\big(p(v)-p(\bar{V})\big)_y=\Big(\frac{u_y}{v}\Big)_y-\Big(\frac{\bar{U}_y}{\bar{V}}\Big)_y-G_y-\frac{1}{8}\Big(\frac{\sigma_{y}^{2}}{(\sigma+1) v^{2}}\Big)_{y},\\
\displaystyle\sigma_{\tau}=-2(\phi+\bar V)(\sigma+1)\sigma+(\phi+\bar V)\Big(\frac{\sigma_{\xi}}{(\phi+\bar V)}\Big)_{y}-\frac{\sigma_{y}^{2}}{2(\sigma+1)},\\
\displaystyle (\phi,\psi,\sigma)\Big|_{\tau=-\frac{t_0}{\epsilon}}=(0,0,\chi_0^2-1)\xrightarrow {y\rightarrow\pm\infty} (0, 0, 0).
\end{array}\right.
\end{equation}
Using the antiderivative technique for the above \eqref{system before interaction}$_{1,2}$, that is, letting
\begin{equation}\label{antiferivative technique}
  \Phi(y,\tau)=\int_{-\infty}^y\phi(z,\tau)dz,\ \  \Psi(y,\tau)=\int_{-\infty}^y\psi(z,\tau)dz,
\end{equation}
one gets
 \begin{equation}\label{equation of composite shock wave before interaction}
 \left\{\begin{array}{llll}
\displaystyle \Phi_\tau-\Psi_y=0, \\
\displaystyle \Psi_\tau+\big(p(v)-p(\bar{V})\big)=\frac{u_y}{v}-\frac{\bar{U}_y}{\bar{V}}-G-\frac{\sigma_{y}^{2}}{8(\sigma+1) (\Phi_y+\bar V)^{2}},\\
\displaystyle \sigma_{\tau}=-2(\Phi_y+\bar V)(\sigma+1)\sigma+(\Phi_y+\bar V)\Big(\frac{\sigma_{y}}{\Phi_y+\bar V}\Big)_{y}-\frac{\sigma_{y}^{2}}{2(\sigma+1)},\\
\displaystyle (\Phi,\Psi,\sigma)\big|_{\tau=-\frac{t_0}{\epsilon}}=(0,0,\chi_0^2-1)\xrightarrow {y\rightarrow\pm\infty} (0, 0, 0).
\end{array}\right.
\end{equation}
Linearizing (\ref{equation of composite shock wave before interaction}) around  $(\bar{V},\bar{U},1)$, one has
 \begin{equation}\label{linearization for equation of composite shock wave before interaction}
 \left\{\begin{array}{llll}
\displaystyle \Phi_\tau-\Psi_y=0, \\
\displaystyle \Psi_\tau+p'(\bar{V})\Phi_y-\frac{\Psi_{yy}}{\bar{V}}=-Q-G+H_1,\\
\displaystyle\sigma_{\tau}-\sigma_{yy}=-2(\Phi_y+\bar V)(\sigma+1)\sigma+H_2,\\
\displaystyle (\Phi,\Psi,\sigma)\big|_{\tau=-\frac{t_0}{\epsilon}}=(0,0,\chi_0^2-1)\xrightarrow {y\rightarrow\pm\infty} (0, 0, 0).
\end{array}\right.
\end{equation}
where
\begin{equation}\label{the estimate of Q}  \left.\begin{array}{llll}
\displaystyle Q=p(\bar{V}+\Phi_{y})-p(\bar{V})-p'(\bar{V})\Phi_{y}-\big(\frac{1}{\bar{V}+\Phi_{y}}-\frac{1}{\bar V}\big)(\Psi_{yy}+\bar{U}_y),\\
\displaystyle H_1=-\frac{\sigma_{y}^{2}}{8(\sigma+1) (\Phi_y+\bar V)^{2}},\ \ H_2=-\frac{\sigma_{y}(\Phi_{yy}+\bar V_y)}{\Phi_y+\bar V}-\frac{\sigma_y^{2}}{2(\sigma+1)}.
\end{array}\right.
\end{equation}
We look for the solution of the initial value problem (\ref{linearization for equation of composite shock wave before interaction})
in the following space:
\begin{equation}\label{X-delta}\left.\begin{array}{llll}
&X_{m,M}(I)=\Big\{(\Phi, \Psi,\sigma)\Big|(\Phi,\Psi) \in C(I;H^3(\mathbb{R})),\sigma\in C(I;H^2(\mathbb{R})),\Phi_y\in L^2(I;H^2(\mathbb{R})),\\
&\displaystyle\qquad\qquad\qquad\qquad\Psi_y,\sigma_y\in L^2(I;H^3(\mathbb{R})),\ \sup_{ t \in I}\| (\Phi,\Psi)\|_{H^3(\mathbb{R})}+\| \sigma\|_{H^2(\mathbb{R})}\leq M\Big\},
\end{array}\right.\end{equation}
where $I\subseteq[-\frac{t_0}{\epsilon},0]$ is an interval.
Now we give the existence results of the system (\ref{linearization for equation of composite shock wave before interaction}) before the interaction time, that is
\begin{theorem} \label{existence thm before t0}
(Existence Before the Interaction Time) There are constants $\delta_0$ and $C$  such that if
 \begin{equation}\label{delta}
\| \chi_{0y}\|_{{2}}+\|\chi_0^2-1\|+| v_+-v_-|\leq\delta_0, \end{equation}
 there exists
a unique solution $(\Phi,\Psi,\sigma)\in X_\delta(-\frac{t_{0}}{\epsilon},0)$ to (\ref{linearization for equation of composite shock wave before interaction}). Furthermore, it holds that
\begin{equation}\label{global solution estimation}\left.\begin{array}{llll}
\displaystyle\|(\Phi,\Psi)(\tau)\|_{H^3(\mathbb{R})}^2+\|\sigma(\tau)\|_{H^2(\mathbb{R})}^2+\int_{-\frac{t_0}{\epsilon}}^{\tau}\big(\|\bar V_y^{\frac12}(\Phi,\Psi)(s)\|^2+\|\sigma(s)\|_{H^3(\mathbb{R})}^2\big)ds\\
\displaystyle\qquad+\int_{-\frac{t_0}{\epsilon}}^{\tau}\big(\|\Phi_y(s)\|_{H^2(\mathbb{R})}^2+\|\Psi_y(s)\|_{H^3(\mathbb{R})}^2\big)ds
\leq C\delta_0^{\frac{1}{3}}e^{-\frac{c\delta_0|t-t_0|}{2\epsilon}}.\end{array}\right.
\end{equation}
and
\begin{equation}\label{maximum value}
  -1\leq \chi\leq 1.
\end{equation}
\end{theorem}
For the proof of the Theorem 3.1,  the standard continuation argument is adopted. The basic steps of this argument is to obtain the existence and uniqueness for the local solution of the system (\ref{linearization for equation of composite shock wave before interaction}) through the fixed point theorem, and then find a way to obtain the consistent energy estimation of the local solution with respect to time, and then get the global solution. Because the proof for local solutions is fundamental and conventional, it will not be given in detail for the sake of brevity.  Here we only devote ourselves to obtaining the a priori estimates.
Setting
\begin{eqnarray}\label{a priori assumption}
N(T)\xlongequal{\mathrm{def}} \sup_{\tau\in[-\frac{t_0}{\epsilon},T]}\big(\| (\Phi,\Psi)\|_{H^3(\mathbb{R})}+\|\sigma\|_{H^2(\mathbb{R})}\big)\leq M_0, \,\,\, for\,\,  T\leq0,
\end{eqnarray}
where $[-\frac{t_0}{\epsilon},T]$ is the time interval on which the local solution is assumed to exist, and $M_0$ is a 
positive constant which will be determined later. Obviously, by using the embedding theorem, the  assumption
(\ref{a priori assumption}) with $M_0$ small enough  holds that
\begin{equation}\label{The upper and lower bounds of density}
  \frac{v_-}{2}\leq \bar{V}+\Phi_{y}\leq 2v_+,\qquad -M_0\leq\sigma\leq M_0,
\end{equation}
moreover, noticing that the following equation which be derived from  \eqref{equation of composite shock wave before interaction}$_3$, using the \eqref{a priori assumption} again, for $M_0$ small enough, one has
\begin{equation}\label{The upper and lower bounds of sigma}
  \sigma_{\tau}-\sigma_{yy}+2(\Phi_y+\bar V)\sigma+\frac{\sigma_{y}(\Phi_{yy}+\bar V_y)}{\Phi_y+\bar V}=-\frac{\sigma_y^{2}}{2(\sigma+1)}-2(\Phi_y+\bar V)\sigma^2\leq0,
\end{equation}
and so, using the parabolic equation maximum principle, one obtains  $\sigma\leq0$ and \eqref{maximum value} is  established.
Further, the nonlinear term $Q, H$ in the system \eqref{linearization for equation of composite shock wave before interaction} satisfy the following:
\begin{equation}\label{the estimate of Q and H}\left.\begin{array}{llll}
 \displaystyle
 Q=O(1)\big(\Phi_y^2+|\Phi_y\Psi_{yy}|+|\bar{U}_y\Phi_y|\big),\quad H_1=\frac{\sigma_{y}^{2}}{8(\sigma+1) (\Phi_y+\bar V)^{2}}=O(1)|\sigma_y^2|,\\
 \displaystyle H_2=O(1)\big(|\sigma_y\Phi_{yy}|+|\sigma_y\bar V_y|+|\sigma_y^2|\big).
 \end{array}\right.
\end{equation}
Moreover, by Lemma \ref{estimate of viscous shock}, it is straightforward to obtain
\begin{equation}\label{the estimate of G}
  |G|\leq  C\delta^{2}e^{-c\delta_0|y|-c\delta_0|\tau|}.
\end{equation}
The following is the a priori estimate of the local solution for (\ref{linearization for equation of composite shock wave before interaction}).

\begin{proposition} \label{a Priori Estimate Before the Interaction Time}
(A Priori Estimate Before the Interaction Time)   Suppose that there exists a solution $(\Phi,\Psi,\sigma)\in X[-\frac{t_0}{\epsilon},T]$ with $T\leq0$,
then there exist positive constants $\delta_0$, $M_0$ and $C$ such that, if $\delta\leq\delta_0$ and $N(T)\leq M_0$,
then $(\Phi,\Psi)$ holds that for $\tau\in [-\frac{t_0}{\epsilon},T]$
 \begin{equation}\label{the estimate before interaction time}\left.\begin{array}{llll}
\displaystyle\|(\Phi,\Psi)(\tau)\|_{H^3(\mathbb{R})}^2+\|\sigma(\tau)\|_{H^2(\mathbb{R})}^2+\int_{-\frac{t_0}{\epsilon}}^{\tau}\Big(\|\bar V_y^{\frac12}(\Phi,\Psi)(s)\|^2+\|\sigma(s)\|_{H^3(\mathbb{R})}^2\Big)ds\\
\displaystyle\qquad+\int_{-\frac{t_0}{\epsilon}}^{\tau}\Big(\|\Phi_y\|_{H^2(\mathbb{R})}^2+\|\Psi_y(s)\|_{H^3(\mathbb{R})}^2\Big)ds
\leq C\delta_0^{\frac{1}{3}}e^{-\frac{c\delta_0|t-t_0|}{2\epsilon}}.\end{array}\right.
\end{equation}
\end{proposition}

The  a priori estimate \eqref{the estimate before interaction time} in proposition 3.1 are obtained step by step from the following lemmas.
\begin{lemma} \label{lem of lower estimate}
Under the assumption of Proposition \ref{a Priori Estimate Before the Interaction Time}, if $\delta_0$ and $M_0$  are suitably small, then for
$-\frac{t_0}{\epsilon}\leq\tau\leq T\leq0$,  it holds that
\begin{equation}\label{the lower estimate}\left.\begin{array}{llll}
 \displaystyle
\|(\Phi,\Psi,\sigma)(\tau)\|^2+\int_{-\frac{t_0}{\epsilon}}^{\tau}\big(\|\bar V_y^{\frac12}(\Phi,\Psi)(s)\|^2+\|\sigma(s)\|^2\big)ds+\int_{-\frac{t_0}{\epsilon}}^{\tau}\|(\Psi_y,\sigma_y)(s)\|^2ds\\
 \displaystyle\leq C\delta_0^{\frac{1}{3}}e^{-\frac{c\delta_0|\tau|}{2}}+C(\delta_0+M_0)\int_{-\frac{t_0}{\epsilon}}^{\tau}\|(\Phi_y,\Phi_{yy},\Psi_{yy})(s)\|^2ds.\end{array}\right.
\end{equation}
\end{lemma}
\begin{proof} In order to obtain the low order energy estimate of $\Phi,\Psi,\sigma$,  the weighted integral method is used as following.
 For $1\geq\alpha>0$,
 multiplying (\ref{linearization for equation of composite shock wave before interaction})$_1$ by
$(1+\alpha\bar V)\Phi$, (\ref{linearization for equation of composite shock wave before interaction})$_2$ by $-\frac{(1+\alpha\bar V)}{p'(\bar{V})}\Psi$,  (\ref{linearization for equation of composite shock wave before interaction})$_3$ by $\beta\sigma$ respectively, where $\alpha>0$ and $\beta>0$ are  undetermined constants, adding the three resulting equations together,
 and integrate the result with respect to $y$ over $\mathbb{R}$, note that $-\bar V_\tau=s_{2_1}V^{S_{2_1}}_y+s_{2_2}V^{S_{2_2}}_y>0$, $p''(\bar V)>0>p'(\bar V)$, one has
\begin{eqnarray}\label{lower energy estimate}
  &&\Big(\frac{1}{2}\int_{\mathbb{R}}\big((1+\alpha\bar{V})\Phi^2+\frac{(1+\alpha\bar{V})}{|p'(\bar{V})|}\Psi^2+\beta\sigma^2\big)dy\Big)_{\tau}\notag\\
  &&+\int_{\mathbb{R}}\frac{s_{2_1}V^{S_{2_1}}_y+s_{2_2}V^{S_{2_2}}_y}{2}\Big(\alpha\Phi^{2}+\frac{2\alpha(V^{S_{2_1}}_y+V^{S_{2_2}}_y)}
  {s_{2_1}V^{S_{2_1}}_y\!+\!s_{2_2}V^{S_{2_2}}_y}\Phi\Psi+\big(\frac{p''(\bar V)}{p'^2(\bar V)}+\frac{\alpha}{|p'(\bar V)|}\big)\Psi^{2}\Big)dy\notag\\
&&+\int_{\mathbb{R}}\frac{1+\alpha(V^{S_{2_1}}_y+V^{S_{2_2}}_y)}{\bar V|p'(\bar V)|}\Psi^{2}_ydy+\int_{\mathbb{R}}\beta\big(\sigma_y^2+2(\Phi_y+\bar V)(\sigma+1)\sigma^2\big)dy\\
&&
=\int_{\mathbb{R}}\Big(\frac{1+\alpha\bar V}{p'(\bar{V})}\Psi\big(G-Q\big)+\big(\frac{1+\alpha\bar V}{\bar Vp'(\bar V) }\big)_{\bar V}'\bar{V}_y\Psi\Psi_y-\frac{(1+\alpha\bar V)\Psi H_1}{p'(\bar{V})}+\beta H_2\sigma\Big)dy\notag\\
&&\leq C\int_{\mathbb{R}}\big(|\Psi||G|+|\Psi||Q|+\bar{V}_y|\Psi||\Psi_y|+|\Psi| |H_1|+\beta|H_2||\sigma|\big)dy.\notag
\end{eqnarray}
For the part contained in parentheses in the second term  of (\ref{lower energy estimate}), choosing  $\alpha$ small enough, by direct calculation,  one obtains
\begin{equation}\label{quadratic form}
 \alpha\Phi^{2}+\frac{2\alpha(V^{S_{2_1}}_y+V^{S_{2_2}}_y)}
  {s_{2_1}V^{S_{2_1}}_y\!+\!s_{2_2}V^{S_{2_2}}_y}\Phi\Psi+\big(\frac{p''(\bar V)}{p'^2(\bar V)}+\frac{\alpha}{|p'(\bar V)|}\big)\Psi^{2}\geq c(\Phi^2
+\Psi^2),\end{equation}
 Moreover,  using the Sobolev embedding theorem, combining the H$\mathrm{\ddot{o}}$lder inequality and Young inequality, one gets 
\begin{equation}\left.\begin{array}{llll}
 \displaystyle
\int_{-\frac{t_0}{\epsilon}}^\tau\int_{\mathbb{R}}|G||\Psi|dyds&\leq&\displaystyle\int_{-\frac{t_0}{\epsilon}}^\tau\int_{\mathbb{R}}\delta^{2}|\Psi|e^{-c\delta|y|-c\delta|s|}dyds\\
&\leq&\displaystyle  C\delta\int_{-\frac{t_0}{\epsilon}}^\tau\|\Psi(s)\|^{\frac12}\|\Psi_y(s)\|^{\frac12}e^{-c\delta|s|}ds\\
\displaystyle&\leq& \displaystyle N(T)\int_{-\frac{t_0}{\epsilon}}^\tau\|\Psi_y(s)\|^2ds+C\delta_0^{\frac13}e^{-\frac{c\delta_0|\tau|}{2}},
\end{array}\right.
\end{equation}
\begin{equation}\left.\begin{array}{llll}
 \displaystyle
 \int_{-\frac{t_0}{\epsilon}}^\tau\int_{\mathbb{R}}|Q||\Psi|dyds
&\leq&\displaystyle C(\delta+M_{0})\int_{-\frac{t_0}{\epsilon}}^\tau
\|(\Phi_y,\Psi_{yy})(s)\|^2ds\\
\displaystyle &+&\displaystyle\nu\int_{-\frac{t_0}{\epsilon}}^\tau\|\Phi_{y}\|^2ds+\frac{C\delta_0}{\nu}
\int_{-\frac{t_0}{\epsilon}}^\tau\int_{\mathbb{R}}\bar{V}_y\Psi^{2}dyds,
\end{array}\right.
\end{equation}
\begin{equation}\left.\begin{array}{llll}
 \displaystyle
 \int_{-\frac{t_0}{\epsilon}}^\tau\int_{\mathbb{R}}\bar{V}_y|\Psi||\Psi_y|dyds\leq\nu \int_{-\frac{t_0}{\epsilon}}^\tau\|\Psi_y(s)\|^2ds+\frac{C\delta}{\nu} \int_{-\frac{t_0}{\epsilon}}^\tau\int_{\mathbb{R}}\bar V_y|\Psi|^2dyds,
 \end{array}\right.
\end{equation}

\begin{equation}\left.\begin{array}{llll}
 \displaystyle
 \int_{-\frac{t_0}{\epsilon}}^\tau\int_{\mathbb{R}}|\Psi| |H_1|dyds&\leq& \displaystyle C\|\Psi(s)\|^{\frac12}\|\Psi_y(s)\|^{\frac12}\int_{-\frac{t_0}{\epsilon}}^\tau\int_{\mathbb{R}}|\sigma_y|^2dyds\\
 &\leq& \displaystyle CM_0\int_{-\frac{t_0}{\epsilon}}^\tau\|\sigma_y(s)\|^2ds,
 \end{array}\right.
\end{equation}
and
\begin{equation}\left.\begin{array}{llll}
 \displaystyle
\int_{-\frac{t_0}{\epsilon}}^\tau \int_{\mathbb{R}}|\sigma| |H_2|dyds\leq \displaystyle\int_{-\frac{t_0}{\epsilon}}^\tau \int_{\mathbb{R}}\big(|\sigma||\sigma_y|^2+|\sigma||\sigma_y||\Phi_{yy}|+|\sigma||\sigma_y|\bar V_y\big)dyds\\
 \leq\displaystyle C\|\sigma\|^{\frac12}\|\sigma_y\|^{\frac12}\int_{-\frac{t_0}{\epsilon}}^\tau\big(\|\sigma_y(s)\|^2+\|\Phi_{yy}(s)\|^2\big)ds+\delta^2\int_{-\frac{t_0}{\epsilon}}^\tau\|\sigma_y(s)\|\|\sigma(s)\|ds\\
 \displaystyle\leq C(\delta_0+M_0)\int_{-\frac{t_0}{\epsilon}}^\tau\big(\|\sigma\|^2+\|\sigma_y(s)\|^2+\|\Phi_{yy}(s)\|^2\big)ds,
  \end{array}\right.
\end{equation}
integrating (\ref{lower energy estimate}) over $[-\frac{t_0}{\epsilon},\tau]$, for sufficiently large weighted parameter $\beta$ and sufficiently small weighted parameter $\alpha$, combining with the Lemma \ref{estimate of viscous shock},
(\ref{the lower estimate})  is deduced immediately by selecting the initial disturbance $M_0$  and the shock wave intensity $\delta_0$ are suitable small, and the proof of Lemma \ref{lem of lower estimate} is completed.
\end{proof}
\begin{lemma} \label{lower order estimate before interaction}
Under the assumption of Proposition \ref{a Priori Estimate Before the Interaction Time}, if $\delta_0$ and $M_0$  are suitably small, then for
$-\frac{t_0}{\epsilon}\leq\tau\leq T\leq0$,  it holds that
\begin{equation}\label{the lower order estimate before interaction}\left.\begin{array}{llll}
\displaystyle\|(\Phi,\Psi,\sigma)(\tau)\|_1^2+\int_{-\frac{t_0}{\epsilon}}^{\tau}\Big(\|\bar V_y^{\frac12}(\Phi,\Psi)(s)\|^2+\|\Phi_y\|^2+\|\Psi_y(s)\|_1^2+\|\sigma(s)\|_2^2\Big)ds\\
\displaystyle\leq C\delta_0^{\frac{1}{3}}e^{-\frac{c\delta_0|t-t_0|}{2\epsilon}}+C(\delta_0+M_0)\int_{-\frac{t_0}{\epsilon}}^{\tau}\|\Phi_{yy}(s)\|^2ds.\end{array}\right.
\end{equation}
\end{lemma}
\begin{proof}
Applying $\partial_y$ to (\ref{linearization for equation of composite shock wave before interaction})$_1$, one gets
\begin{equation}\label{The derivative of the mass conservation equation}
   \Phi_{y\tau}-\Psi_{yy}=0,
 \end{equation}
multiplying the above equality \eqref{The derivative of the mass conservation equation}  by
$\Phi_y$, (\ref{linearization for equation of composite shock wave before interaction})$_2$ by $-\bar{V}\Phi_y$ respectively,
adding the two resulting equations together,
 and integrate the result with respect to $y$ over $\mathbb{R}$, one has
\begin{equation}\label{derivative estimation for Phi}\left.\begin{array}{llll}
 \displaystyle
  \Big(\int_{\mathbb{R}}\frac{1}{2}\Phi_y^2dy\Big)_\tau+\int_{\mathbb{R}}\big(\bar{V}|p'(\bar V)|\Phi_y^2-\bar{V}\Psi_{\tau}\Phi_y\big)dy=-\int_{\mathbb{R}}\bar{V}\Phi_y\big(Q-G+H_1\big) dy,\end{array}\right.
\end{equation}
observing that
\begin{equation}\label{derivative estimation for Phi-1}
  \int_{\mathbb{R}}\bar{V}\Psi_{\tau}\Phi_ydy=\Big(\int_{\mathbb{R}}\bar{V}\Psi\Phi_ydy\Big)_\tau+\int_{\mathbb{R}}
  \big((s_{2_1}V_{y}^{S_{2_1}}+s_{2_2}V_{y}^{S_{2_2}})\Psi\Phi_y+\bar V_y\Psi\Psi_y+\bar{V}\Psi_y^2\big)dy,
\end{equation}
plugging  \eqref{derivative estimation for Phi-1} into (\ref{derivative estimation for Phi}), integrating the resulting equation over $[-\frac{t_0}{\epsilon},\tau]$,  one gets
\begin{equation}\label{the inequality of higher order estimate1}
\left.\begin{array}{llll}
\displaystyle\int_{\mathbb{R}}\Big(\frac{1}{2}\Phi_y^2-\bar{V}\Psi\Phi_y\Big)dy+\int_{-\frac{t_0}{\epsilon}}^\tau\int_{\mathbb{R}}\bar{V}|p'(\bar V)|\Phi_y^2dyds\\
\displaystyle\leq\int_{-\frac{t_0}{\epsilon}}^\tau\int_{\mathbb{R}}
  \big((s_{2_1}V_{y}^{S_{2_1}}+s_{2_2}V_{y}^{S_{2_2}})\Psi\Phi_y+\bar V_y\Psi\Psi_y+\bar{V}\Psi_y^2\big)dyds\\
\displaystyle\quad-\int_{-\frac{t_0}{\epsilon}}^\tau\int_{\mathbb{R}}\bar{V}\Phi_y\big(Q-G+H_1\big) dyds,
\end{array}\right.
\end{equation}
multiplying both sides of \eqref{the lower estimate} by a constant large enough, adding it to  (\ref{the inequality of higher order estimate1}), by using the Sobolev embedding theorem, for smallness of  $\delta$ and $M_0$, one obtains
\begin{equation}\label{the higher order estimate of the deriative for Phi}\left.\begin{array}{llll}
\displaystyle\|(\Psi,\sigma)(\tau)\|^2+\|\Phi(\tau)\|_1^2
+\int_{-\frac{t_0}{\epsilon}}^{\tau}\|(\Phi_y,\Psi_y,\sigma_y)(s)\|^2ds\\
\displaystyle\qquad+\int_{-\frac{t_0}{\epsilon}}^{\tau}\big(\|\bar V_y^{\frac12}(\Phi,\Psi)(s)\|^2+\|\sigma(s)\|^2\big)ds\\
\displaystyle\leq C\delta_0^{\frac{1}{3}}e^{-\frac{c\delta_0|t-t_0|}{2\epsilon}}+C(\delta+M_0)\int_{-\frac{t_0}{\epsilon}}^{\tau}\|(\Phi_{yy},\sigma_{yy},\Psi_{yy})(s)\|^2ds.\end{array}\right.
\end{equation}
Multiplying (\ref{linearization for equation of composite shock wave before interaction})$_2$ by $-\Psi_{yy}$, (\ref{linearization for equation of composite shock wave before interaction})$_3$ by $-\sigma_{yy}$, and integrate it over $\mathbb{R}$ with respect to $y$, one has
\begin{equation}\label{the second derivative for Psi}
  \Big(\frac{1}{2}\int_{\mathbb{R}}\Psi_{y}^{2}dy\Big)_\tau+\int_{\mathbb{R}}\frac{1}{\bar{V}}\Psi_{yy}^2dy=\int_{\mathbb{R}}\big(p'(\bar V)\Psi_{yy}\Phi_y-\Psi_{yy}(Q-G+H_1)\big)dy,
\end{equation}
and
\begin{equation}\label{the second deriative for sigma}
 \Big(\frac{1}{2}\int_{\mathbb{R}}\sigma_{y}^{2}dy\Big)_\tau+\int_{\mathbb{R}}\sigma_{yy}^2dy=\int_{\mathbb{R}}\sigma_{yy}\big(2(\Phi_y+\bar V)(\sigma+1)\sigma- H_2\big)dy,
\end{equation}
integrating the above equality over $[\frac{t_0}{\epsilon},\tau]$, by using Sobolev embedding theorem, repeating the previous proof of \eqref{the higher order estimate of the deriative for Phi} steps,  the inequality \eqref{the lower order estimate before interaction} can be derived, and  the proof of Lemma \ref{lower order estimate before interaction} is completed.
\end{proof}

\

\noindent \textbf{Proof of Proposition \ref{a Priori Estimate Before the Interaction Time}.}
In the proof of the Lemma 3.1, we need to use the assumptions \eqref{a priori assumption} for $(\Phi,\Psi,\sigma)$, in fact, $\sup_{\tau\in[-\frac{t_0}{\epsilon},T]}\big\{\| (\Phi,\sigma)\|_2,\|\Psi\|_1\big\}$ is required to be small enough. In order to keep the proof of the energy estimate closed, we need to give the following higher derivative estimate of $(\Phi,\Psi,\sigma)$ for the system \eqref{linearization for equation of composite shock wave before interaction}, for which,
applying $\partial_{yy}$
 to (\ref{linearization for equation of composite shock wave before interaction})$_1$ and $\partial_{y}$ to (\ref{linearization for equation of composite shock wave before interaction})$_2$
one has
 \begin{equation}\label{first derivative of linearization for equation of composite shock wave before interaction}
 \left\{\begin{array}{llll}
\displaystyle \Phi_{yy\tau}-\Psi_{yyy}=0, \\
\displaystyle \Psi_{y\tau}+p'(\bar V)\Phi_{yy}-\frac{\Psi_{yyy}}{\bar{V}}=-p''(\bar V)\bar V_y\Phi_y-\frac{\Psi_{yy}\bar V_y}{\bar V^2}+Q_y-G_y+H_{1y},
\end{array}\right.
\end{equation}
and further,
Applying $\partial_{y}$
 to (\ref{first derivative of linearization for equation of composite shock wave before interaction})$_1$, $\partial_{y}$ to (\ref{first derivative of linearization for equation of composite shock wave before interaction})$_2$
  and  $\partial_{y}$ to (\ref{linearization for equation of composite shock wave before interaction})$_3$, one gets
 \begin{equation}\label{second derivative of linearization for equation of composite shock wave before interaction}
 \left\{\begin{array}{llll}
\displaystyle \Phi_{yyy\tau}-\Psi_{yyyy}=0, \\
\displaystyle \Psi_{yy\tau}+p'(\bar V)\Phi_{yyy}-\frac{\Psi_{yyyy}}{\bar{V}}=-p''(\bar V)\bar V_y\Phi_{yyy}-\frac{\Psi_{yyy}\bar V_y}{\bar{V}^2}-\big(p''(\bar V)\bar V_y\Phi_y\big)\\
\displaystyle\qquad\qquad\qquad\qquad\qquad\qquad\qquad-\Big(\frac{\Psi_{yy}\bar V_y}{\bar V^2}\Big)_y+Q_{yy}-G_{yy}+H_{1yy},\\
\displaystyle \sigma_{\tau y}-\sigma_{yyy}=-2(\Phi_{yy}+\bar V_y)(\sigma+1)\sigma-2(\Phi_y+\bar V)(2\sigma+1)\sigma_y+H_{2y},
\end{array}\right.
\end{equation}
Multiplying (\ref{first derivative of linearization for equation of composite shock wave before interaction})$_1$ by $\Phi_{yy}$, (\ref{first derivative of linearization for equation of composite shock wave before interaction})$_2$ by $ -\bar V\Phi_{yy}$ and (\ref{first derivative of linearization for equation of composite shock wave before interaction})$_2$ by $-\Psi_{yyy}$, repeating the analytic steps in Lemma 3.2 step by step, one obtains
\begin{equation}\label{the second order estimate for Phi}\left.\begin{array}{llll}
\displaystyle\|(\Phi,\Psi)(\tau)\|_{H^2(\mathbb{R})}^2+\|\sigma(\tau)\|_{H^1(\mathbb{R})}^2+\int_{-\frac{t_0}{\epsilon}}^{\tau}\Big(\|\bar V_y^{\frac12}(\Phi,\Psi)(s)\|^2+\|\sigma(s)\|_{H^2(\mathbb{R})}^2\Big)ds\\
\displaystyle+\int_{-\frac{t_0}{\epsilon}}^{\tau}\Big(\|\Phi_y(s)\|_{H^1(\mathbb{R})}^2+\|\Psi_y(s)\|_{H^2(\mathbb{R})}^2\Big)ds\leq C\delta_0^{\frac{1}{3}}e^{-\frac{c\delta_0|t-t_0|}{2\epsilon}}.\end{array}\right.
\end{equation}
Multiplying (\ref{second derivative of linearization for equation of composite shock wave before interaction})$_1$ by $\Phi_{yyy}$, (\ref{second derivative of linearization for equation of composite shock wave before interaction})$_2$ by $ -\bar V\Phi_{yyy}$, (\ref{second derivative of linearization for equation of composite shock wave before interaction})$_2$ by $-\Psi_{yyyy}$ and (\ref{second derivative of linearization for equation of composite shock wave before interaction})$_3$ by $-\sigma_{yyy}$, repeating the proof steps in Lemma 3.2 again, combining with the energy inequality \eqref{the second order estimate for Phi}, the a priori estimate \eqref{the lower order estimate before interaction} can be derived, and thus the proof of Proposition \ref{a Priori Estimate Before the Interaction Time} is completed.


\section{Estimates After the Interacting Time}\setcounter{equation}{0}
\subsection{Time-dependent shift}
\indent\qquad
In this section,  we consider the case after the interaction of two homologous shock waves, for which, as shown in Figure 1, after the interaction, a forward shock wave and a backward rarefaction wave are generated.
 Unlike the system (\ref{NSAC-scaling}) where the initial data at $\tau=-\frac{t_0}{\epsilon}$ is well-chosen such that the perturbation $(v,u)\big| _{\tau=-\frac{t_0}{\epsilon}}=(\bar V,\bar U)$, here the initial data is determined by the solution obtained after the interaction time, thanks for Proposition \ref{a Priori Estimate Before the Interaction Time} and \eqref{Strength after interaction}, one has
 \begin{equation}\label{the estimate at the interaction time-1}
\displaystyle\|(v-\bar V,u-\bar U,\omega-1)(0)\|_{H^2(\mathbb{R})}^2\leq C\delta_0^{\frac{1}{3}}e^{-\frac{c\delta_0|t_0|}{2\epsilon}}.
\end{equation}
and
\begin{equation}\label{strenghth of rarefaction and shock wave after interaction}
\max\big\{|v_m-v_-|,|u_m-u_-|\big\}+\max\big\{|(v^+-v_m|,|u^+-u_m)|\big\}=\tilde{\delta}_1+\tilde{\delta}_2\leq \delta_0,
\end{equation}
On the other hand, due to the difference in the analysis techniques of rarefaction wave and shock wave, the anti-derivative method used in Section 3 is no longer applicable to the analysis of the superposition of rarefaction wave and shock wave. Thanks to the inspiration of Kang-Vasseur-Wang \cite{kvw2021}, we introduce time-dependent shift $\mathbf{X}(\tau)$  for these composite wave, and combine the effective velocity, relative quantities and the weighted estimation method to solve this difficulty.
Letting  
\begin{equation}\label{xi}
  \xi=y-\tilde{s}_{2}\tau,
\end{equation}
under the moving coordinate system $(\xi,\tau)$, the system \eqref{NSAC-scaling} is transformed into the following form
\begin{equation}\label{NSAC-xi}
\left\{\begin{array}{l}
\displaystyle v_\tau-\tilde{s}_2 v_{\xi}-u_{\xi}=0, \\
\displaystyle u_\tau-\tilde{s}_2 u_{\xi}+p(v)_{\xi}=\left(\frac{u_{\xi}}{v}\right)_{\xi}-\frac{1}{8}\Big(\frac{\omega_{\xi}^{2}}{\omega v^{2}}\Big)_{\xi},\\
\displaystyle\omega_{\tau}-\tilde{s}_2\omega_{\xi}=-2v\omega (\omega-1)+v\Big(\frac{\omega_{\xi}}{v}\Big)_{\xi}-\frac{\omega_{\xi}^{2}}{2\omega}.
\end{array}\right.
\end{equation}
The following notations are used throughout in this section, i.e.
\begin{gather}\label{f,X}
f^{\pm \mathbf{X}}(\xi,\tau)\xlongequal{\mathrm{def}} f( \xi \pm \mathbf{X}(\tau),\tau),
\end{gather}with any function $f: \mathbb{R} \times \mathbb{R}^{+} \rightarrow \mathbb{R}$,
moreover, the corresponding smooth approximate superposition outgoing waves is constructed as follows
\begin{gather}\label{35}
\big(\tilde{V}, \tilde{U}\big)(\xi,\tau)\xlongequal{\mathrm{def}} \big(\tilde{V}^{\tilde{R}_1},\tilde{U}^{\tilde{R}_1}\big)(\xi+\tilde{s}_2\tau,\tau)+\big(\tilde{V}^{\tilde{S}_2,-\mathbf{X}},\tilde{U}^{\tilde{S}_2,-\mathbf{X}}\big)(\xi)-(v_m,u_m),
\end{gather}
where $\big(\tilde{V}^{\tilde{R}_{1}}, \tilde{U}^{\tilde{R}_{1}}\big)(\xi,\tau)$ is the approximate 1-rarefaction wave defined in \eqref{1-rarefaction wave}, $\big(\tilde{V}^{\tilde{S}_2}, \tilde{U}^{\tilde{S}_2}\big)(\xi)$ is the 2-viscous shock wave defined in \eqref{2 shock after}, and the shift $\mathbf{X}(\tau)$ of the forward viscous shock wave is determined by the following ordinary differential equation:
\begin{equation}\label{Shift X}
\left\{\begin{aligned}
& \frac{d\mathbf{X}}{d\tau}=-\frac{m_0}{\tilde\delta_2} {\left[\int_{\mathbb{R}} \frac{a(\xi-\mathbf{X})}{\tilde{s}_2} \partial_{\xi} \tilde{h}^{\tilde{S}_2,-\mathbf{X}}\big(p(v)-p (\tilde{V}) \big) d \xi\right.} \\
&\quad \quad \quad \quad\quad\qquad\left.-\int_{\mathbb{R}} a(\xi-\mathbf{X}) \partial_{\xi} p\big(\tilde{V}^{\tilde{S}_2,-\mathbf{X}}\big)\big(v-\tilde{V}\big) d \xi\right], \\
& \mathbf{X}(0)=0,
\end{aligned}\right.
\end{equation}
where $m_0=\frac{5(\gamma+1) (-p^{\prime}(v_m))^\frac32}{8 \gamma p\left(v_m\right)}$, $\tilde{h}^{\tilde{S}_2}=\tilde{U}^{\tilde{S}_2}-\big(\ln \tilde{V}^{\tilde{S}_2}\big)_{\xi}$ and the weight function $a: \mathbb{R} \rightarrow \mathbb{R}$ is defined by
\begin{gather}
a(\xi)\xlongequal{\mathrm{def}} 1+\frac{\lambda}{\tilde\delta_2}\Big(p(v_m)-p\big(\tilde{V}^{\tilde{S}_2}(\xi)\big)\Big),\label{412}
\end{gather}
$\lambda$ is a constant which be chosen  so small but far bigger than $\tilde\delta_2$ such that
\begin{gather}
\tilde\delta_2 \ll \lambda \leq C \sqrt{\tilde\delta_2},\label{413}
\end{gather}
and the following properties hold for the weight function $a(\xi)$
\begin{gather}
1<a(\xi)<1+\lambda,
\quad
a_\xi^{\prime}(\xi)=-\frac{\lambda}{\tilde\delta_2} p^{\prime}\big(\tilde{V}^{\tilde{S}_2}\big) \tilde{V}_{\xi}^{\tilde{S}_2}>0,\quad \big|a_\xi^{\prime}\big| \sim \frac{\lambda}{\tilde\delta_2}\big|\tilde{V}_{\xi}^{\tilde{S}_2}\big| .\label{416}
\end{gather}
The following Lemma \ref{the properities for Shift X} shows the properties of the waves with shift $\mathbf{X}(\tau)$ (see Kang-Vasseur-Wang \cite{kvw2021}).
\begin{lemma}\label{the properities for Shift X}
Let $\mathbf{X}$ be the shift defined by \eqref{Shift X}. There exist constants $\epsilon_1,\delta_0>0$, such that $|\dot{\mathbf{X}}|\leq \epsilon_1$ and $\tilde\delta_1,\tilde{\delta}_2\leq \delta_0$, then for all $t \leq T$,
\begin{equation}
\begin{aligned}
& \big\|\tilde{V}^{\tilde{S}_2,-\mathbf{X}}_{\xi} \big(\tilde{V}^{\tilde{R}_1}-v_m\big)\big\|_{L^1(\mathbb{R})}+\big\|\tilde{V}^{\tilde{R}_1}_{\xi} \tilde{V}^{\tilde{S}_2,-\mathbf{X}}_{\xi} \big\|_{L^1(\mathbb{R})} \leq C \tilde{\delta}_1 \tilde\delta_2 e^{-C \tilde\delta_2 \tau}, \\
& \big\|\tilde{V}^{\tilde{S}_2,-\mathbf{X}}_{\xi} \big(\tilde{V}^{\tilde{R}_1}-v_m\big)\big\|_{L^2(\mathbb{R})}+\big\|\tilde{V}^{\tilde{R}_1}_{\xi} \tilde{V}^{\tilde{S}_2,-\mathbf{X}}_{\xi} \big\|_{L^2(\mathbb{R})} \leq C \tilde{\delta}_1 \tilde\delta_2^{\frac 32} e^{-C \tilde\delta_2 \tau}, \\
& \big\|\tilde{V}^{\tilde{R}_1}_{\xi} \big(\tilde{V}^{{\tilde{S}_2},-\mathbf{X}}-v_m\big)\big\|_{L^2(\mathbb{R})} \leq C \tilde{\delta}_1 \tilde\delta_2 e^{-C \tilde\delta_2 \tau}.
\end{aligned}
\end{equation}
\end{lemma}

It is easy to check that the approximate superposition wave $\big(\tilde{V}, \tilde{U}\big)$ constructed in \eqref{35} satisfies the system
\begin{equation}\label{equation of approximate composite shock wave after interaction}
\left\{\begin{array}{l}
\displaystyle \tilde{V} _\tau-\tilde{s}_2 \tilde{U} _{\xi}- \tilde{U} _{\xi}=-\dot{\mathbf{X}}(\tau)\tilde{V}^{\tilde{S}_2,-\mathbf{X}}_{\xi}  ,\\
\displaystyle \tilde{U} _\tau-\tilde{s}_2 \tilde{U}_{\xi}+p_{\xi} (\tilde{V})=\Big(\frac{\tilde{U}_{\xi}}{\tilde{V}}\Big)_{\xi}-\dot{\mathbf{X}}(\tau)\tilde{U}^{\tilde{S}_2,-\mathbf{X}}_{\xi} +F_1+F_2,
\end{array}\right.
\end{equation}
where
\begin{equation}
F_1=\Big(\frac{\tilde{U}^{\tilde{S}_2,-\mathbf{X}}_{\xi}}{\tilde{V}^{\tilde{S}_2,-\mathbf{X}}}\Big)_{\xi}-\big(\frac{ \tilde{U} _{\xi}}{\tilde{V}}\Big)_{\xi}, \quad F_2=\big(p (\tilde{V}) -p(\tilde{V}^{\tilde{R}_1})-p(\tilde{V}^{\tilde{S}_2,-\mathbf{X}})\big)_{\xi} .\label{310}
\end{equation}
The following Theorem 4.1 is about the existence, uniqueness and asymptotic stability of the solution for system \eqref{NSAC-xi} after the interaction time of the chasing shock waves.

\begin{theorem}\label{thm1}
Let $(v^r,u^r)$ and $(\tilde V^{\tilde{S}_2}, \tilde U^{\tilde{S}_2})$ are the outgoing rarefaction wave and viscous shock wave given by \eqref{23} and  \eqref{smooth-wave} respectively, there exist a positive constant $\tilde{\delta}_0(<\delta_0)$ such that, if the strength satisfies $\delta_1+\delta_2<\tilde{\delta}_0$, then there exists a unique global solution $(v, u,\omega)$ to  \eqref{NSAC-xi}  and \eqref{key Lemma of the estimate for initial data} in  $\tau\in[0,+\infty)$, moreover, it holds that
\begin{eqnarray}
&& v(y,\tau)-\Big( v^{r_1}(\frac{y}{\tau})+\tilde V^{\tilde{S}_2}(y-\tilde{s}_2\tau-X(\tau))-v_m \Big)\in C\big([0,+\infty) ; H^1(\mathbb{R})\big), \notag\\
&& u(y,\tau)-\Big( u^{r_1}(\frac{y}{\tau})+\tilde U^{\tilde{S}_2}(y-\tilde{s}_2\tau-X(\tau))-u_m\Big) \in C\big([0,+\infty) ; H^1(\mathbb{R})\big), \notag\\
&& u_{yy}(y,\tau)-\tilde{U}_{yy}^{\tilde{S}_2}(y-\tilde{s}_2 \tau-\mathbf{X}(\tau)) \in L^2\big([0,+\infty) ; L^2(\mathbb{R})\big),\notag\\
&& \omega-1 \in C\left([0,+\infty ); H^2(\mathbb{R})\right)\cap L^2\big([0,+\infty) ; H^3(\mathbb{R})\big).\notag
\end{eqnarray}
furthermore, 
\begin{equation}\label{sv}
    \lim_{\tau\rightarrow +\infty}\sup _{x \in \mathbb{R}}\left|(v,u,\omega)-(\tilde{V},\tilde{U},1)\right| = 0,
\end{equation}
and
\begin{equation}\label{large time behavior of X(t)}
    \lim_{\tau\rightarrow +\infty}|\dot{\mathbf{X}}(\tau)|= 0.
\end{equation}
\end{theorem}

\begin{remark}
The implication of the theorem 4.1 is that, when two interfaces move in the same direction with different shock wave velocities, if the speed of the rear interface is faster than the speed of the front interface, the rear interface will catch up with the front interface and interact with each other in a finite time, so that it will evolve into a forward interface moving at the shock wave velocity and a backward interface moving at the rarefaction wave, and the tail of the newly generated interface moving backward at the rarefaction wave velocity will produce a shift effect to the interface moving forward at the shock velocity,  the shift is time dependent.
\end{remark}

\begin{remark}
From \eqref{large time behavior of X(t)}, one derives at once that
\begin{equation}\label{sub linearly for X(t)}
 \lim_{\tau\rightarrow +\infty}\frac{\mathbf{X}(\tau)}{\tau}=0,
\end{equation}
which means that the function $\mathbf{X}(\tau)$  grows at most sub-linearly with respect to time $t$.
\end{remark}

Similar to the proof of Theorem 3.1,  the proof of the local solution is omitted,  and the a priori estimate of the solution is given here in detail.
Putting the perturbation around the superposition wave profile $(\tilde{V},\tilde{U},1)$ by
\begin{eqnarray}\label{per1}
\big(\phi,\psi,\sigma\big)(y,\tau)=\big(v-\tilde V,u-\tilde U,\omega-1\big)(y,\tau),\ \ \tau>0,
\end{eqnarray}
subtracting \eqref{equation of approximate composite shock wave after interaction} from \eqref{NSAC-xi}, one has
\begin{equation}\label{nsac-per}
\left\{\begin{array}{l}
\displaystyle \phi _\tau-\tilde{s}_2 \phi _{\xi}- \psi _{\xi}=\dot{\mathbf{X}}(t)\tilde{V}^{\tilde{S}_2,-\mathbf{X}}_{\xi}  ,\\
\displaystyle \psi _\tau-\tilde{s}_2 \psi _{\xi}+\big(p (v)-p (\tilde{V}) \big)_{\xi}=\Big(\frac{ \psi_{\xi}}{v}\Big)_{\xi}-\Big(\frac{ \phi\tilde{U}_{\xi}}{v\tilde{V}}\Big)_{\xi}-\frac{1}{8}\Big(\frac{\sigma_{\xi}^{2}}{(\sigma+1)v^{2}}\Big)_{\xi}\\
\displaystyle\qquad\qquad\qquad\qquad\qquad\qquad\qquad+\dot{\mathbf{X}}(t)\tilde{U}^{\tilde{S}_2,-\mathbf{X}}_{\xi}-F_1-F_2,\\
\displaystyle \sigma_{\tau}-\tilde{s}_2\sigma_{\xi}=-2v(\sigma+1)\sigma+v\Big(\frac{\sigma_{\xi}}{v}\Big)_{\xi}-\frac{\sigma_{\xi}^{2}}{2(\sigma+1)},
\end{array}\right.
\end{equation}
with the initial condition as following
\begin{equation}\label{init-per}
 (\phi,\psi,\sigma)(\xi,0)=(v-\tilde{V},u-\tilde{U},\omega-1)(\xi,0)\xrightarrow {\xi\rightarrow\pm\infty} (0, 0, 0).
\end{equation}
The a prior estimate of the moment of interaction of the chasing shock waves is given below, which means that when the strength of the outgoing waves is suitably small, the initial perturbation in the sense of \eqref{the estimate at the interaction time} is also small after the interacting time.
\begin{lemma}\label{the estimate at the interaction time}
It holds that
\begin{eqnarray}\label{key Lemma of the estimate for initial data}
 \Big\|\big(v-\tilde V,u-\tilde U,\omega-1\big)\Big|_{\tau=0}\Big\|_{H^2(\mathbb{R})}=O(1)\delta_0^{\frac13},
\end{eqnarray}
for $\delta_0$ small enough.
\end{lemma}

\begin{proof}
Since the difference between the viscous shock profiles before and the rarefaction wave and the viscous shock profiles after the interacting time are
 \begin{eqnarray}
\big(\tilde{V}-\bar{V}\big)(y,0)&=&\big(\tilde{V}-\bar{V}\big)(y,0)=\big(\tilde{V}^{\tilde{R}_1}(y,0)+\tilde{V}^{\tilde{S}_2}(y)-v_m\big)-\big(V^{S_{2_1}}(y)+V^{S_{2_2}}(y)-v^*\big)\notag\\
&=&\left\{ \begin{aligned}
(\tilde{V}^{\tilde{R}_1}(y,0)-v_-)&-(V^{S_{2_2}}(y)-v^*)\\
&+(\tilde{V}^{\tilde{S}_2}(y)-v_m)-(V^{S_{2_1}}-v_-),
y\leq0,\notag\\
(\tilde{V}^{\tilde{R}_1}(y,0)-v_m)&-(V^{S_{2_2}}(y)-v_+)\\
&+(\tilde{V}^{\tilde{S}_2}(y)-v_+)-(V^{S_{2_1}}-v^*), y\geq0,
\end{aligned} \right.
\end{eqnarray}
and
 \begin{eqnarray}
\big(\tilde{U}-\bar{U}\big)(y,0)&=&\big(\tilde{U}-\bar{U}\big)(y,0)=\big(\tilde{U}^{\tilde{R}_1}(y,0)+\tilde{U}^{\tilde{S}_2}(y)-v_m\big)-\big(U^{S_{2_1}}(y)+V^{S_{2_2}}(y)-v^*\big)\notag\\
&=&\left\{ \begin{aligned}
(\tilde{U}^{\tilde{R}_1}(y,0)-v_-)&-(U^{S_{2_2}}(y)-v^*)\\
&+(\tilde{U}^{\tilde{S}_2}(y)-v_m)-(U^{S_{2_1}}-v_-),
y\leq0,\notag\\
(\tilde{U}^{\tilde{R}_1}(y,0)-v_m)&-(U^{S_{2_2}}(y)-v_+)\\
&+(\tilde{U}^{\tilde{S}_2}(y)-v_+)-(U^{S_{2_1}}-v^*), y\geq0,
\end{aligned} \right.
\end{eqnarray}
by using Lemma \eqref{the properties of rarefaction wave} and Lemma \eqref{estimate of viscous shock}, one has
\begin{equation*}
  \big|(\tilde V-\bar V)(y,0)\big|=O(1)\delta_0e^{-c\delta_0|y|},\quad \big|(\tilde U-\bar U)(y,0)\big|=O(1)\delta_0e^{-c\delta_0|y|},
\end{equation*}
which gives
\begin{equation}\label{Small perturbations of V after interfaction}
 \Big\|\big(\tilde V-\bar V,\tilde U-\bar U\big)\Big|_{\tau=0}\Big\|_{H^2(\mathbb{R})}=O(1)\delta_0^{\frac12},
\end{equation}
combining with \eqref{the estimate at the interaction time-1}, then the a priori estimate \eqref{the estimate at the interaction time} is obtained immediately, and the proof of Lemma 4.1 is completed.
\end{proof}

Now, similar to the proof procedure  in Section 3,  setting the solution space:
\begin{equation}\label{solution space after interfaction}
\left.\begin{array}{llll}
&X_M([0,T])=\Big\{(\phi, \psi,\sigma)\Big|(\phi,\sigma) \in C([0,T];H^2(\mathbb{R})),\psi\in C([0,T];H^1(\mathbb{R})),\\
&\displaystyle\qquad\qquad\qquad\qquad\sup_{ t \in [0,T]}\| (\phi,\psi,\sigma)\|_{H^1(\mathbb{R})}\leq M,\|(\phi_{yy},\sigma_{yy})\|_{L^{\infty}(0, T ; L^2(\mathbb{R}))}\leq M_1\Big\}.
\end{array}\right.\end{equation}
From \eqref{solution space after interfaction}, by using  the Sobolev embedding, one has
\begin{align}\label{basic1}
\|(\phi,\psi,\sigma)\|_{L^{\infty}((0, T) \times \mathbb{R})} \leq C M,\quad
\|(\phi_{y},\sigma_{y})\|_{L^{\infty}((0, T) \times \mathbb{R})} \leq C \sqrt{MM_1}.
\end{align}
This smallness together with \eqref{Shift X} and Lemma \ref{propoties of relative quantity}  yields that
\begin{align}
&\|p(v)-p(\tilde{V})\|_{L^{\infty}((0, T) \times \mathbb{R})} \leq C\|\phi\|_{L^{\infty}((0, T) \times \mathbb{R})} \leq C M ,\label{49}\\
&|\dot{\mathbf{X}}(\tau)| \leq \frac{C}{\tilde\delta_2}\big\|\,|p(v)|+|p(\tilde{V})|+| \phi \big|\big\|_{L^{\infty}(\mathbb{R})} \int_{\mathbb{R}}\tilde{V}^{\tilde{S}_2}_{\xi}  d \xi \leq C\|\phi\|_{L^{\infty}(\mathbb{R})}.\label{410}
\end{align}

\vskip 0.3cm\begin{proposition}\label{prop}
Suppose that $(\phi,\psi,\sigma)(\xi,\tau)\in X_M([0,T])$ is the solution to the problem \eqref{nsac-per}-\eqref{init-per} on $[0, T]$ for some $T>0$, then, for $\delta_0$ and $M$ small enough, $\forall\tau \leq T$,
it holds that
\begin{equation}\left.\begin{array}{llll}\label{pro1}
&\displaystyle\sup _{[0, T]}\big(\big\|(\phi,\sigma)\big\|^2_{H^2(\mathbb{R})}+\big\|\psi\big\|^2_{H^1(\mathbb{R})}\big)+\tilde\delta_2 \int_0^\tau|\dot{\mathbf{X}}|^2 d s
+\int_0^\tau\int_{\mathbb{R}}(|\tilde{V}^{\tilde{S}_2,-\mathbf{X}}_{\xi}|+|\tilde{U}^{\tilde{R}_1}_{\xi}|)\phi^2 d \xi d s\\
&\displaystyle
+\int_0^\tau\int_{\mathbb{R}}\big|\big(p(v)-p(\tilde{V})\big)_\xi\big|^2 d \xi d s
+\int_0^\tau\int_{\mathbb{R}}\big|\phi_{\xi\xi}\big|^2 d \xi d s+\int_0^\tau\big(\|\psi_\xi\|^2_{H^1(\mathbb{R})}+\|\sigma\|^2_{H^3(\mathbb{R})}\big) d s \\
&\displaystyle\leq C_0\Big(\big\|(\phi,\sigma)(0)\big\|_{H^2(\mathbb{R})}+\big\|\psi(0)\big\|_{H^1(\mathbb{R})}+\delta_0^{\frac13}\Big),\notag
\end{array}\right.\end{equation}
 where  $C_0$  is independent of $T$.
\end{proposition}

\subsection{Energy estimates for weighted relative entropy}
\indent\qquad
In this subsection, the energy estimates for weighted relative entropy will be established. In order to get more information about $v$, thanks to the entropy method provided by Vasseur-Yao \cite{VY-2020}, He-Huang \cite{HH-2020} to overcome  the difficulties caused by the initial perturbations with small energy but possibly large oscillations of shock waves, the following effective velocity is introduced:
\begin{equation}\label{effective velocity}
\left.\begin{array}{l}
\displaystyle h(\xi,\tau)\xlongequal{\mathrm{def}} u-(\ln v)_{\xi}, \qquad\qquad\qquad\tilde{h}^{\tilde{S}_2}(\xi)\xlongequal{\mathrm{def}} \tilde{U}^{\tilde{S}_2}-\big(\ln \tilde{V}^{\tilde{S}_2}\big)_{\xi},\\
\displaystyle\tilde{h}(\xi,\tau)\xlongequal{\mathrm{def}} \tilde{U}^{\tilde{R}_1}+ \tilde{h}^{\tilde{S}_2,-\mathbf{X}}-u_m,
\end{array}\right.
\end{equation}
Further, the system \eqref{NSAC-xi}, the system \eqref{2 shock after} and the system \eqref{equation of approximate composite shock wave after interaction} are transformed into the following new forms respectively
\begin{equation}
\left\{\begin{array}{l}
\displaystyle v_\tau-\tilde{s}_{2} v_{\xi}-h_{\xi}=(\ln v)_{\xi \xi}, \\
\displaystyle h_\tau-\tilde{s}_{2} h_{\xi}+p(v)_{\xi}=-\frac{1}{8}\Big(\frac{\omega_{\xi}^{2}}{\omega v^{2}}\Big)_{\xi},\\
\displaystyle \omega_{\tau}-\tilde{s}_{2}\omega_{\xi}=-2v\omega (\omega-1)+v\Big(\frac{\omega_{\xi}}{v}\Big)_{\xi}-\frac{\omega_{\xi}^{2}}{2\omega}.\label{42}
\end{array}\right.
\end{equation}
and
\begin{equation}
\left\{\begin{array}{l}
\displaystyle -\tilde{s}_2\tilde{V}^{\tilde{S}_2}_{\xi} - \tilde{h}^{\tilde{S}_2}_\xi=(\ln \tilde{V}^{\tilde{S}_2})_{\xi\xi}, \\
\displaystyle -\tilde{s}_2 \tilde{h}^{\tilde{S}_2}_{\xi}+\big(p(\tilde{V}^{\tilde{S}_2})\big)_{\xi}=0, \\
\displaystyle \lim_{\xi\rightarrow-\infty}\big(\tilde{V}^{\tilde{S}_2}, \tilde{h}^{\tilde{S}_2}\big)=(v_m, u_m), \quad\lim_{\xi\rightarrow +\infty}\big(\tilde{V}^{\tilde{S}_2}, \tilde{h}^{\tilde{S}_2}\big)=(v_{+}, u_{+}) .\label{433}
\end{array}\right.
\end{equation}
and
\begin{equation}
\left\{\begin{array}{l}
\tilde{V}_\tau-\tilde{s}_2 \tilde{V}_{\xi}-\tilde{h}_{\xi}=(\ln \tilde{V})_{\xi \xi}-\dot{\mathbf{X}}(\tau)\tilde{V}^{\tilde{S}_2,-\mathbf{X}}_{\xi} +F_3, \\
\tilde{h}_\tau-\tilde{s}_2 \tilde{h}_{\xi}+(p(\tilde{V}))_{\xi}=-\dot{\mathbf{X}}(\tau) \tilde{h}^{\tilde{S}_2,-\mathbf{X}}_{\xi}+F_2,
\end{array}\right.\label{45}
\end{equation}
where $F_2$ is defined in \eqref{310} and
\begin{gather}
F_3=\left(\ln \tilde{V}^{\tilde{S}_2,-\mathbf{X}}-\ln \tilde{V}\right)_{\xi \xi} .\label{46}
\end{gather}
To simplify the system \eqref{42} and the system \eqref{45}, introducing the following vector
\begin{equation}\label{vectors-v-h-w}
\mathbf{w}\xlongequal{\mathrm{def}} \left(\begin{array}{l}
v \\
h \\
\omega
\end{array}\right), \qquad
\tilde{\mathbf{w}}\xlongequal{\mathrm{def}} \left(\begin{array}{l}
\tilde{V} \\
\tilde{h}\\
\tilde{\omega}
\end{array}\right)
=\left(\begin{array}{c}
\tilde{V}^{\tilde{R}_1}(\xi,\tau)+\tilde{V}^{\tilde{S}_2,-\mathbf{X}}(\xi)-v_m \\
\tilde{U}^{\tilde{R}_1}(\xi,\tau)+ \tilde{h}^{\tilde{S}_2,-\mathbf{X}}(\xi)-u_m\\
1
\end{array}\right),
\end{equation}
and the following matrix operators
\begin{equation}\label{matrix-A-M}
A(\mathbf{w})\xlongequal{\mathrm{def}} \left(\begin{array}{c}
\displaystyle-\tilde{s}_2 v-h \\
\displaystyle-\tilde{s}_2 h+p(v)\\
\displaystyle-\tilde{s}_2\omega
\end{array}\right),
\   M(\mathbf{w})\xlongequal{\mathrm{def}} \left(\begin{array}{ccc}
\displaystyle \frac{1}{\gamma p(v)} & 0 & 0\\
\displaystyle 0 & 0& \displaystyle-\frac{\omega_{\xi}}{8v^2\omega}\\
\displaystyle 0 & 0& 1
\end{array}\right),
\end{equation}
\begin{equation}\label{matrix-H}
  H(\mathbf{w})\xlongequal{\mathrm{def}} \left(
\begin{array}{c}
0 \\ 0\\
\displaystyle-2v\omega (\omega-1)-\frac{\omega_{\xi}v_{\xi}}{v}\!-\!\frac{\omega_{\xi}^{2}}{2\omega}
\end{array}
\right),
\end{equation}
further, we define the entropy function
\begin{equation}\label{entropy}
  \eta(\mathbf{w})\xlongequal{\mathrm{def}} Q(v)+\frac{h^2}{2}+\frac{\omega^2}{2},
\end{equation}
where
  \begin{equation}\label{Q}
   Q(v)=\frac{v^{-\gamma+1}}{\gamma-1}.
 \end{equation}
Moreover, the relative entropy of $\eta$, the relative flux of $A$ and the relative pressure of $p$ are as follows
\begin{eqnarray}
&&\eta\big(\mathbf{w}  \big| \tilde{\mathbf{w}}\big)=\eta(\mathbf{w})-\eta(\tilde{\mathbf{w}})-\nabla \eta(\tilde{\mathbf{w}})(\mathbf{w}-\tilde{\mathbf{w}})
=Q(v \big| \tilde{V})+\frac{(h-\tilde{h})^2}{2}+\frac{\sigma^2}{2},\label{422}\\
&&A\big(\mathbf{w}  \big| \tilde{\mathbf{w}}\big)=A(\mathbf{w})-A(\tilde{\mathbf{w}})-\nabla A(\tilde{\mathbf{w}})(\mathbf{w}-\tilde{\mathbf{w}})= \left(\begin{array}{c}
0 \\
p\big(v \big| \tilde{V}\big)\\ 0
\end{array}\right).\\
&&p(v\big| \tilde{V})=p(v)-p(\tilde{V})-p'(\tilde{V})(v-\tilde{V}). \label{423-1}
\end{eqnarray}
Thus, the system \eqref{42} and \eqref{45} are transformed into the following vector forms respectively
\begin{gather}
\partial_\tau \mathbf{w}+\partial_{\xi} A(\mathbf{w})=\partial_{\xi}\Big(M(\mathbf{w}) \partial_{\xi} \nabla \eta(\mathbf{w})\Big)+H(\mathbf{w}).\label{419}
\end{gather}
  \begin{equation}
\partial_\tau \tilde{\mathbf{w}}+\partial_{\xi} A(\tilde{\mathbf{w}})=\partial_{\xi}\Big(M(\tilde{\mathbf{w}}) \partial_{\xi} \nabla \eta(\tilde{\mathbf{w}})\Big)-\dot{\mathbf{X}} \partial_{\xi}\tilde{\mathbf{w}}^{\tilde{S}_2,-\mathbf{X}}+\left(\begin{array}{l}
F_3 \\ F_2 \\ 0
\end{array}\right),\label{421}
\end{equation}
where
\begin{equation}\label{tilde-vector-w-matrix-A-M}
 \tilde{\mathbf{w}}^{\tilde{S}_2,-\mathbf{X}}=\left(\begin{array}{c}
\displaystyle\tilde{V}^{\tilde{S}_2,-\mathbf{X}} \\
\displaystyle\tilde{h}^{\tilde{S}_2,-\mathbf{X}}\\
\displaystyle1\notag
\end{array}\right),
 \ A(\tilde{\mathbf{w}})=\left(\begin{array}{c}
\displaystyle-\tilde{s}_2 \tilde{V}-\tilde{h} \\
\displaystyle-\tilde{s}_2 h+p(\tilde{V})\\
\displaystyle-\tilde{s}_2\notag
\end{array}\right),
\quad   M(\tilde{\mathbf{w}})=\left(\begin{array}{ccc}
\displaystyle \frac{1}{\gamma p(\tilde{V})} & 0 & 0\\
\displaystyle 0 & 0& \displaystyle0\\
\displaystyle 0 & 0& 1\notag
\end{array}\right),
\end{equation}
and $F_2, F_3$ are defined in \eqref{310}, \eqref{46} respectively.
By using \eqref{Shift X}, \eqref{412}, \eqref{419}-\eqref{421}, the expression of the rate of change of the weighted relative entropy integral with respect to time can be obtained by direct calculation
\begin{gather}\label{wre1}
\frac{d}{d \tau} \int_{\mathbb{R}} a^{-\mathbf{X}} \eta\big(\mathbf{w}\big| \tilde{\mathbf{w}}\big) d \xi=\dot{\mathbf{X}}(\tau) \mathbf{Y}(\mathbf{w})+\mathcal{B}(\mathbf{w})-\mathcal{G}(\mathbf{w}),
\end{gather}
where
\begin{equation}
  \mathbf{Y}(\mathbf{w}) =\int_{\mathbb{R}} a^{-\mathbf{X}} \tilde{h}_{\xi}^{\tilde{S}_2,{-\mathbf{X}}}\hbar d \xi-\int_{\mathbb{R}} a^{-\mathbf{X}} p^{\prime}(\tilde{V}) \tilde{V}_{\xi}^{\tilde{S}_2,{-\mathbf{X}}}\phi d \xi-\int_{\mathbb{R}} a_{\xi}^{-\mathbf{X}}\big(\frac{| h-\tilde h |^2}{2}+Q(v \big| \tilde{V})\big) d \xi ,\notag
\end{equation}
\begin{equation}
\begin{aligned}
\mathcal{B}(\mathbf{w}) & =\frac{1}{2 \tilde{s}_2} \int_{\mathbb{R}} a_{\xi}^{-\mathbf{X}}|p(v)-p(\tilde{V})|^2 d \xi+\tilde{s}_2 \int_{\mathbb{R}} a^{-\mathbf{X}}\tilde{V}^{\tilde{S}_2,-\mathbf{X}}_{\xi}  p(v \big| \tilde{V}) d \xi \\
& -\int_{\mathbb{R}} a_{\xi}^{-\mathbf{X}} \frac{p(v)-p(\tilde{V})}{\gamma p(v)} \partial_{\xi}(p(v)-p(\tilde{V})) d \xi+\int_{\mathbb{R}} a_{\xi}^{-\mathbf{X}}(p(v)-p(\tilde{V}))^2 \frac{\partial_{\xi} p(\tilde{V})}{\gamma p(v) p(\tilde{V})} d \xi \\
& -\int_{\mathbb{R}} a^{-\mathbf{X}} \partial_{\xi}(p(v)-p(\tilde{V})) \frac{p(\tilde{V})-p(v)}{\gamma p(v) p(\tilde{V})} \partial_{\xi} p(\tilde{V}) d \xi+\int_{\mathbb{R}} a^{-\mathbf{X}}(p(v)-p(\tilde{V})) F_3 d \xi \\
& -\int_{\mathbb{R}} a^{-\mathbf{X}} (h-\tilde h)  F_2 d \xi-\frac{1}{8}\int_{\mathbb{R}} a^{-\mathbf{X}} (h-\tilde h) \Big(\frac{\sigma_{\xi}^{2}}{(\sigma+1) v^{2}}\Big)_{\xi}  d \xi-\int_{\mathbb{R}} a^{-\mathbf{X}}_{\xi}\sigma\sigma_{\xi}d \xi \\
& -\int_{\mathbb{R}} a^{-\mathbf{X}}\Big(2v\sigma^3+\frac{(\phi_{\xi}+\tilde{V}_{\xi})\sigma_{\xi}\sigma}{v}+\frac{\omega_{\xi}^{2}\sigma}{\sigma+1}\Big)  d \xi ,\notag
\end{aligned}
\end{equation}
and
\begin{equation}
    \begin{aligned}
\mathcal{G}(\mathbf{w}) & =\frac{\tilde{s}_2}{2} \int_{\mathbb{R}} a_{\xi}^{-\mathbf{X}}\big| (h-\tilde h) -\frac{p(v)-p(\tilde{V})}{\tilde{s}_2}\big|^2 d \xi+\tilde{s}_2 \int_{\mathbb{R}} a_{\xi}^{-\mathbf{X}} Q(v | \tilde{V}) d \xi+\int_{\mathbb{R}} a^{-\mathbf{X}} \tilde{U}_{\xi}^{\tilde{R}_1} p\big(v | \tilde{V}\big) d \xi\\
& +\int_{(\mathbb{R}} (\frac{(\tilde{s}_2 a_{\xi}^{-\mathbf{X}}+4va^{-\mathbf{X}})\sigma^2}{2}  d \xi+\int_{\mathbb{R}} a^{-\mathbf{X}} |\sigma_\xi|^2 d \xi +\int_{\mathbb{R}} \frac{a^{-\mathbf{X}}\big|\partial_{\xi}(p(v)-p(\tilde{V}))\big|^2}{\gamma p(v)} d \xi.\notag
\end{aligned}
\end{equation}
At this point of preparation, the relative entropy energy inequality is given as follows:
\begin{lemma}\label{lem-vh}
 Under the hypotheses of Proposition \ref{prop}, there exists $C>0$ (independent of $\left.\delta_0, M, T\right)$ such that for all $t \in(0, T]$,
\begin{eqnarray}\label{47}
&&\displaystyle \int_{\mathbb{R}}Q(v\big|\tilde{V})d \xi+\| h-\tilde h \|_{L^{2}(\mathbb{R})}^{2}+\|\sigma\|_{H^{1}(\mathbb{R})}^{2}+\tilde{\delta}_2 \int_0^\tau|\dot{\mathbf{X}}|^2 d s\notag\\
&&\quad
+\int_0^\tau(G_1+G^{\tilde{S}_2}+D_1+\|\sigma\|_{H^{2}(\mathbb{R})}^{2}) d s\\
&&\displaystyle \leq C \int_{\mathbb{R}}Q(v\big|\tilde{V})(\xi,0)d \xi+\|(h-\tilde{h})(0)\|_{L^{2}(\mathbb{R})}^{2}+\|\omega(0)-1\|_{H^{1}(\mathbb{R})}^{2}+C \tilde{\delta}_1^{\frac13},\notag
\end{eqnarray}
where $h(\xi,0)=u(\xi,0)-\big(\ln v(\xi,0)\big)_{\xi}$ and
\begin{equation}\label{48}
    \begin{aligned}
 G_1=&\frac{\lambda}{\tilde{\delta}_2} \int_{\mathbb{R}}|\tilde{V}^{\tilde{S}_2,-\mathbf{X}}_{\xi} |\cdot\Big| h-\tilde{h} -\frac{p(v)-p(\tilde{V})}{\tilde{s}_2}\Big|^2 d \xi,\\
G^{\tilde{S}_2}=&\int_{\mathbb{R}}|\tilde{V}^{{\tilde{S}_2},-\mathbf{X}}_{\xi}|\big|p(v)-p(\tilde{V})\big|^2 d \xi,\qquad
D_1=\int_{\mathbb{R}}\big|(p(v)-p(\tilde{V}))_\xi\big|^2 d \xi.
    \end{aligned}
\end{equation}
\end{lemma}
\begin{proof}
By using the translation coordinate $\xi \mapsto \xi+\mathbf{X}(\tau)$  for \eqref{wre1}, one has
\begin{gather}\label{4333}
\frac{d}{d \tau} \int_{\mathbb{R}} a(\xi) \eta(\mathbf{w}^{\mathbf{x}} \big| \tilde{\mathbf{w}}^{\mathbf{x}}) d \xi=\dot{\mathbf{X}}(\tau) \mathbf{Y}\left(\mathbf{w}^{\mathbf{x}}\right)+\mathcal{B}\left(\mathbf{w}^{\mathbf{x}}\right)-\mathcal{G}\left(\mathbf{w}^{\mathbf{x}}\right),
\end{gather}
Without causing ambiguity, for the sake of brevity, we omit the dependence of the solution on the shift, that is, $(v,h,\omega)=(v^{\mathbf{X}},h^{\mathbf{X}},\omega^{\mathbf{X}}), \mathbf{w}=\mathbf{w}^{\mathbf{X}}$. Moreover, the following substitute symbols are introduced
\begin{equation}
\begin{aligned}
&\Delta P\xlongequal{\mathrm{def}} p(v)-p\big(\tilde{V}^{\mathbf{X}}\big),
\mathcal{B}\xlongequal{\mathrm{def}} \sum_{i=1}^3 \mathbf{B}_i+\sum_{i=1}^3\mathbf{P}_i+\sum_{i=1}^2\mathbf{S}_i,\\
&\mathcal{G}\xlongequal{\mathrm{def}} \sum_{i=1}^2\mathbf{G}_i+\mathbf{G}^{\tilde{R}_1}+\mathbf{D}_1+\mathbf{A}_1,\label{4343}
\end{aligned}
\end{equation}
where
\begin{equation}
\begin{aligned}
& \mathbf{B}_1=\frac{1}{2 \tilde{s}_2} \int_{\mathbb{R}} a_{\xi}(\xi) (\Delta P)^2 d \xi+\tilde{s}_2 \int_{\mathbb{R}} a(\xi)\tilde{V}^{\tilde{S}_2}_{\xi} p\big(v \big| \tilde{V}^{\mathbf{X}}\big) d \xi, \\ &\mathbf{B}_2=\int_{\mathbb{R}} \frac{a(\xi) (\Delta P)^2_{\xi} p_{\xi}\big(\tilde{V}^{\mathbf{X}}\big)}{2\gamma p(v) p\big(\tilde{V}^{\mathbf{X}}\big)}  d \xi-\int_{\mathbb{R}} \frac{a_{\xi} (\Delta P)^2_{\xi}}{2\gamma p(v)} d \xi, \\
& \mathbf{B}_3=\int_{\mathbb{R}}\frac{a_{\xi} (\Delta P)^2  p_{\xi}(\tilde{V}^{\mathbf{X}})}{\gamma p(v) p(\tilde{V}^{\mathbf{X}})} d \xi, \ \mathbf{P}_1=- \frac{1}{8}\int_{\mathbb{R}} a(h-\tilde h^\mathbf{X})\big(\frac{2\sigma_{\xi}\sigma_{\xi\xi}}{(\sigma+1)v^2}-\frac{\sigma_{\xi}^3}{(\sigma+1)^2v^2} -\frac{2\sigma_{\xi}^2\tilde{V}^{\mathbf{X}}_{\xi}}{(\sigma+1)v^3}\big)d\xi, \\
& \mathbf{P}_2=-\int_{\mathbb{R}} a_{\xi}\sigma\sigma_{\xi}d \xi -\int_{\mathbb{R}} a\Big(2v\sigma^3-\frac{\tilde{V}^{\mathbf{X}}_{\xi}\sigma_{\xi}\sigma}{v}-\frac{\sigma_{\xi}^{2}\omega}{\sigma+1}\Big)  d \xi, \\
& \mathbf{P}_3=\frac{1}{4}\int_{\mathbb{R}} a(h-\tilde h^\mathbf{X})\frac{\sigma_{\xi}^2\phi_{\xi}}{(\sigma+1)v^3}d\xi+\int_{\mathbb{R}} a\frac{\phi_{\xi}\sigma_{\xi}\sigma}{v}d\xi,
\notag
\end{aligned}
\end{equation}
and
\begin{equation}
\begin{aligned}
& \mathbf{S}_1=\int_{\mathbb{R}} a (\Delta P)\big(\ln \tilde{V}^{\tilde{S}_2}-\ln \tilde{V}^{\mathbf{X}}\big)_{\xi \xi} d \xi, \ \  \mathbf{S}_2=-\int_{\mathbb{R}} a(h-\tilde h^\mathbf{X})\big(p(\tilde{V}^{\mathbf{X}})-p(\tilde{V}^{\tilde{R}_1,\mathbf{X}})-p(\tilde{V}^{\tilde{S}_2})\big)_{\xi} d \xi,\\
&\mathbf{G}_1=\frac{\tilde{s}_2}{2} \int_{\mathbb{R}} a_{\xi}\big| (h-\tilde h^{\mathbf{X}} )-\frac{ \Delta P}{\tilde{s}_2}\big|^2 d \xi, \ \mathbf{G}_2=\tilde{s}_2 \int_{\mathbb{R}} a_{\xi} Q\big(v \big| \tilde{V}^{\mathbf{X}}\big) d \xi, \ \mathbf{G}^{\tilde{R}_1}=\int_{\mathbb{R}} a\tilde{U}_{\xi}^{\tilde{R}_1,\mathbf{X}} p\big(v \big| \tilde{V}^{\mathbf{X}}\big) d \xi, \\
& \mathbf{D}_1=\int_{\mathbb{R}} \frac{a}{\gamma p(v)}\big| (\Delta P)_{\xi}\big|^2 d \xi,\quad\mathbf{A}_1=\int_{\mathbb{R}} \frac{\tilde{s}_2 a_{\xi}^{-\mathbf{X}}+4va^{-\mathbf{X}}}{2}  \sigma^2 d \xi+\int_{\mathbb{R}} a|\sigma_{\xi}|^2 d \xi.\notag
\end{aligned}
\end{equation}
Further, from \eqref{Shift X} and \eqref{wre1},  the functional $\mathbf{Y}$ can be decomposed into the following form
\begin{equation}\label{yyy}
    \begin{aligned}
\mathbf{Y}&= -\int_{\mathbb{R}} a_{\xi}(\xi)\Big(\frac{| (h-\tilde h^{\mathbf{X}} ) |^2}{2}+Q(v \big| \tilde{V}^{\mathbf{X}})\Big) d \xi +\int_{\mathbb{R}} a(\xi) \tilde{h}_{\xi}^{\tilde{S}_2}(h-\tilde h^{\mathbf{X}} ) d \xi\\
&\quad-\int_{\mathbb{R}} a(\xi) p^{\prime}(\tilde{V}^{\mathbf{X}}) \tilde{V}_\xi^{\tilde S_2}(v-\tilde V^\mathbf{X} )d \xi
=\sum_{i=1}^5 \mathbf{Y}_i,
\end{aligned}
\end{equation}
where $\mathbf{Y}_i,(i=1,\cdots,5)$ satisfy
\begin{equation}
\begin{aligned}
\mathbf{Y}_1=&\int a(\xi)\Big(\frac{\tilde{h}_{\xi}^{\tilde{S}_2} \Delta P}{\tilde{s}_2}  - p'(\tilde{V}^{\tilde{S}_2})\tilde{V}^{\tilde{S}_2}_\xi(v-\tilde V^\mathbf{X})\Big) d \xi,\quad
\mathbf{Y}_2=\int  a(\xi) \tilde{h}_{\xi}^{\tilde{S}_2}\Big( (h-\tilde h^{\mathbf{X}} )-\frac{ \Delta P}{{\tilde{s}_2}}\Big) d \xi,\\
\mathbf{Y}_3=&-\int  a(\xi)\big(p^{\prime}(\tilde{V}^{\mathbf{X}})-p^{\prime}(\tilde{V}^{\tilde{S}_2})\big) \tilde{V}_{\xi}^{\tilde{S}_2}\phi d \xi,\\
\mathbf{Y}_4=&-\frac{1}{2} \int_{\mathbb{R}}  a_{\xi}(\xi)\big(h-\tilde h^{\mathbf{X}} -\frac{ \Delta P}{{\tilde{s}_2}}\big)\big( h-\tilde h^{\mathbf{X}} +\frac{ \Delta P}{{\tilde{s}_2}}\big) d \xi,\\
\mathbf{Y}_5=&-\int  a_{\xi}(\xi) \Big(Q(v \big| \tilde{V}^{\mathbf{X}})+\frac{ (\Delta P)^2}{2 \tilde{s}_2^2}\Big) d \xi.\notag
\end{aligned}
\end{equation}
Notice from \eqref{Shift X} that $\dot{\mathbf{X}}(\tau)=-\frac{m_0}{\tilde{\delta}_2}\mathbf{Y}_1$
which together with \eqref{yyy} yields
\begin{gather}
\dot{\mathbf{X}}(\tau) \mathbf{Y}=-\frac{\tilde{\delta}_2}{m_0}|\dot{\mathbf{X}}(\tau)|^2+\dot{\mathbf{X}}(\tau) \sum_{i=2}^5 \mathbf{Y}_i,\label{439}
\end{gather}
and then, by using the Poincar$\acute{e}$ inequality below (Kang-Vasseur-Wang \cite{kvw2021})
\begin{equation}\label{Poincare inequality}
 \int_0^1\Big|g-\int_0^1gdz\Big|^2dz\leq\frac12\int_0^1z(1-z)\big|g'\big|dz,
\end{equation}
where $\forall g:[0,1]\rightarrow\mathbb{R}$ and  satisfiying $\int_0^1z(1-z)\big|g'\big|dz<\infty$,  through the tedious but not difficult calculations, one immediately knows, there exist $C, C_1>0$ such that:
\begin{equation}\label{key11}
\begin{aligned}
& -\frac{\tilde{\delta}_2}{2 M}|\dot{\mathbf{X}}|^2+\mathbf{B}_1-\mathbf{G}_2-\frac{3}{4} \mathbf{D}_1\leq-C_1 \mathbf{G}^{\tilde{S}_2}+C (K_1+K_2),
\end{aligned}
\end{equation}
where
\begin{equation}
    \begin{aligned}
&\mathbf{G}^{\tilde{S}_2}=\int |\tilde { V } ^ { \tilde{S}_2 } _ { \xi } | (\Delta P)^2 d \xi, \\
&K_1=\int | a_{\xi}(\xi)|  (\Delta P)^3 d \xi,\qquad
K_2=\int | a_{\xi}(\xi)||\tilde{V}^{\tilde{R}_1,\mathbf{X}}-v_m|  (\Delta P)^2 d \xi.
    \end{aligned}
\end{equation}
Next, by using \eqref{4333}, \eqref{439} and \eqref{key11}, we have
\begin{equation}
\begin{aligned}
\frac{d}{d \tau} \int_{\mathbb{R}}  a(\xi) \eta(\mathbf{w} \big| \tilde{\mathbf{w}}^{\mathbf{X}}) d \xi\leq & -\frac{\tilde\delta_2}{2 M}|\dot{\mathbf{X}}|^2+\dot{\mathbf{X}} \sum_{i=2}^5 \mathbf{Y}_i+\sum_{i=2}^3 \mathbf{B}_i+\sum_{i=1}^3 \mathbf{P}_i+\sum_{i=1}^2\mathbf{S}_i\\
&-\mathbf{G}_1-\mathbf{G}^{\tilde{R}_1}-\frac{1}{4} \mathbf{D}_1-\mathbf{A}_1 -C_1 \mathbf{G}^{\tilde{S}_2}+C (K_1+K_2).\notag
\end{aligned}
\end{equation}
By using the Young's inequality, one obtains
\begin{equation}\label{wre2}
\begin{aligned}
&\frac{d}{d \tau} \int_{\mathbb{R}}  a(\xi) \eta(\mathbf{w} \big| \tilde{\mathbf{w}}^{\mathbf{X}}) d \xi+\frac{\tilde{\delta}_2}{4 M}|\dot{\mathbf{X}}|^2+\mathbf{G}_1+\mathbf{G}^{\tilde{R}_1}+\frac{1}{4} \mathbf{D}_1+\mathbf{A}_1+C_1 \mathbf{G}^{\tilde{S}_2} \\
&\leq C(K_1+K_2)+\frac{C}{\tilde{\delta}_2} \sum_{i=2}^5\left|\mathbf{Y}_i\right|^2+\sum_{i=2}^3 \mathbf{B}_i+\sum_{i=1}^3 \mathbf{P}_i+\sum_{i=1}^2\mathbf{S}_i .
\end{aligned}
\end{equation}
Note that from \eqref{49} and \eqref{416}, it holds
\begin{gather}\label{ggg}
G_1 \sim \mathbf{G}_1, \quad G^{\tilde{S}_2}=\mathbf{G}^{\tilde{S}_2}, \quad D_1 \sim \mathbf{D}_1, \quad \|\sigma\|_{H^{1}}^{2} \sim \mathbf{A}_1.
\end{gather}
 we will use the above good terms $\mathbf{G}_1, \mathbf{G}^{\tilde{R}_1}, \mathbf{G}^{\tilde{S}_2}, \mathbf{D}_1$ and $\mathbf{A}_1$ in \eqref{wre2} to control the above bad terms on the right-hand side of \eqref{wre2},
and the following are appropriate estimates for the these bad terms in turn. The first step is to consider the estimation of  $K_{1}$ and $K_2$. By using interpolation inequality, \eqref{412}, \eqref{49} and \eqref{416} with Lemma  \ref{estimate of viscous shock}, one has
\begin{equation}\label{e-k1}
\begin{aligned}
K_1 & \leq C \frac{\lambda}{\tilde{\delta}_2} \| \Delta P\|_{L^{\infty}(\mathbb{R})}^2\int |\tilde{V}^{\tilde{S}_2}_{\xi}|\cdot|  \Delta P| d \xi \\
& \leq C \frac{\lambda}{\tilde{\delta}_2}\big\| (\Delta P)_{\xi}\big\|_{L^2(\mathbb{R})}\| \Delta P\|_{L^2(\mathbb{R})} \big(\int|\tilde{V}^{\tilde{S}_2}_{\xi}| (\Delta P)^2 d \xi\big)^{1/2}\big(\int|\tilde{V}^{\tilde{S}_2}_{\xi}| d \xi\big)^{1/2}\\
& \leq C M\left\| (\Delta P)_{\xi}\right\|_{L^2(\mathbb{R})}^2+C M \int|\tilde{V}^{\tilde{S}_2}_{\xi}| (\Delta P)^2 d \xi \leq \frac{1}{40}\big(\mathbf{D}_1+C_1 \mathbf{G}^{\tilde{S}_2}\big) .
\end{aligned}
\end{equation}
Similarly, from Lemma \ref{the properities for Shift X}, combining with \eqref{solution space after interfaction}, \eqref{413}, \eqref{416},  one gets
\begin{equation}\label{e-k2}
\begin{aligned}
K_2 & \leq C \frac{\lambda}{\tilde{\delta}_2}\| \Delta P\|_{L^4(\mathbb{R})}^2 \|\tilde{V}^{\tilde{S}_2}_{\xi}(\tilde{V}^{{\tilde{R}_1},\mathbf{X}}-v_m) \|_{L^2(\mathbb{R})} \\
& \leq C \frac{\lambda}{\tilde{\delta}_2}\left\| (\Delta P)_{\xi}\right\|_{L^2(\mathbb{R})}^{\frac12}\| \Delta P\|_{L^2(\mathbb{R})}^{\frac32}\tilde{\delta}_2^{\frac32} \tilde{\delta}_1 e^{-C \tilde{\delta}_2 \tau}\\
& \leq C M\big\| (\Delta P)_{\xi}\big\|_{L^2(\mathbb{R})}^2+C M \tilde{\delta}_2^{\frac43} \tilde{\delta}_1^{\frac43} e^{-C \tilde{\delta}_2 \tau} \leq \frac{1}{40} \mathbf{D}_1+C M \tilde{\delta}_2^{\frac43} \tilde{\delta}_1^{\frac43} e^{-C \tilde{\delta}_2 \tau} .
\end{aligned}
\end{equation}
The second step is to give the estimates of the terms $\mathbf{Y}_{i}$ $(i=2,\cdots,5)$. Since $\tilde{h}^{\tilde{S}_2}_{\xi}=\frac{\partial_{\xi}p(\tilde{V}^{\tilde{S}_2})}{\tilde{s}_2}=-\frac{\tilde{\delta}_2}{\lambda \tilde{s}_2}a_{\xi}$, one derives
\begin{align}
&\left|\mathbf{Y}_2\right|\leq C\frac{\tilde{\delta}_2}{\lambda}\int \big|a_{\xi}\big|\cdot\big|h-\tilde{h}-\frac{ \Delta P}{\tilde{s}_2}\big|d \xi \leq C \frac{\tilde{\delta}_2}{\sqrt{\lambda}} |\mathbf{G}_1|^{\frac12}.\notag
\end{align}
From Lemma \ref{estimate of viscous shock}, one gets
\begin{gather}
\left|\mathbf{Y}_3\right| \leq C \int\big|\tilde{V}^{{\tilde{R}_1},\mathbf{X}}-v_m\big||\tilde{V}_{\xi}^{\tilde{S}_2}|\big| \phi\big| d \xi \leq C \tilde{\delta}_1 \int|\tilde{V}_{\xi}^{\tilde{S}_2}|\Delta P d \xi \leq C \tilde\delta_1 \sqrt{\tilde{\delta}_2} |\mathbf{G}^S|^{\frac12}.\notag
\end{gather}
 Observing that
\begin{equation}
\begin{aligned}
| h-\tilde{h}| & \leq\big|u-\tilde{U}\big|+\big|(\ln v)_{\xi}-\big(\ln \tilde{V}^{\tilde{S}_2}\big)_{\xi}\big| \\
& \leq\left|\psi\right|+C\Big(\big| \phi_{\xi}\big|+\big|\tilde{V}_{\xi}^{\mathbf{X}}\big|\big|v-\tilde{V}^{\tilde{S}_2}\big|+\big|\tilde{V}_{\xi}^{{\tilde{R}_1},\mathbf{X}}\big|\Big) \\
& \leq\left|\psi\right|+C\Big(\left|\phi_{\xi}\right|+\big|\tilde{V}_{\xi}^{\mathbf{X}}\big||\phi|+\big|\tilde{V}_\xi^{\tilde S_2}\big|\big|\tilde{V}^{\tilde{R}_1,\mathbf{X}}-v_m\big|
+\big|\tilde{V}_{\xi}^{\tilde{R}_1,\mathbf{X}}\big|\Big),\label{455}
\end{aligned}
\end{equation}
together with  Lemma \ref{the properties of rarefaction wave}, Lemma \ref{estimate of viscous shock} and \eqref{solution space after interfaction}, one gets
\begin{gather}
\| h-\tilde{h}\|_{L^2(\mathbb{R})} \leq C(\|\psi\|_{L^2(\mathbb{R})}+\| \phi\|_{H^1(\mathbb{R})}+\tilde{\delta}_1)\leq C(M+\tilde{\delta}_1) ,\label{456}
\end{gather}
combining with \eqref{456}, \eqref{solution space after interfaction} and $\left\|a_{\xi}\right\|_{L^{\infty}} \leq C \lambda \tilde{\delta}_2$, one obtains
\begin{equation}
\begin{aligned}
& \left|\mathbf{Y}_4\right| \leq C\|a_{\xi}(\xi)\|_{L^{\infty}}^{\frac{1}{2}}\left|\mathbf{G}_1\right|^{\frac{1}{2}}(\| h-\tilde{h}\|_{L^2(\mathbb{R})}+\| \phi\|_{L^2(\mathbb{R}})\leq C\big(M+\tilde{\delta}_1\big)\big(\lambda \tilde{\delta}_2\big)^{\frac{1}{2}} \mathbf{G}_1^{\frac{1}{2}}.\notag
\end{aligned}
\end{equation}
Using Lemma \ref{estimate of viscous shock} with \eqref{49}, one derives
\begin{gather}
\left|\mathbf{Y}_5\right| \leq C\int\left|a_{\xi}(\xi)\right|  (\Delta P)^2 d \xi\leq C\| \Delta P \|_{L^2(\mathbb{R})}|\mathbf{G}_S|^{\frac{1}{2}}\leq CM|\mathbf{G}_S|^{\frac{1}{2}}.\notag
\end{gather}
Thus, combining the estimates for $\mathbf{Y}_{i}$ $(i=2,\cdots,5)$, the following estimates is obtained
\begin{gather}\label{e-yy}
\frac{C}{\tilde{\delta}_2}\sum_{i=2}^5|\mathbf{Y}_i|^2 \leq C\tilde{\delta}_2\mathbf{G}_1+C\tilde{\delta}_1\mathbf{G}^{\tilde{S}_2}+C(M+\tilde{\delta}_1)\mathbf{G}_1 +CM^2\mathbf{G}_S\leq \frac{1}{2}\mathbf{G}_1+\frac{C_1}{20}\mathbf{G}^{\tilde{S}_2}.
\end{gather}
Next,  the terms $\mathbf{B}_i$ $(i=2,3)$ are considered as following. By using Young's inequality, one has
\begin{equation}\label{e-b2}
\begin{aligned}
|\mathbf{B}_2| &\leq C \mathbf{D}_1^{\frac{1}{2}} \big(\int_{\mathbb{R}}|a_{\xi}|^2  (\Delta P)^2 d \xi\big)^{\frac{1}{2}}+ C \mathbf{D}_1^{\frac{1}{2}}\delta_0 \big(\int_{\mathbb{R}}|\tilde{V}^{\mathbf{X}}_{\xi}|^2
 (\Delta P)^2 d \xi\big)^{\frac{1}{2}}\\
&\leq C \mathbf{D}_1^{\frac{1}{2}} \tilde\delta_2|\mathbf{G}^{\tilde{S}_2}|^{\frac{1}{2}} +C \mathbf{D}_1^{\frac{1}{2}}\delta_0 \left(\mathbf{G}^S+\mathbf{G}^R\right)^{\frac{1}{2}}\\
&\leq \frac{1}{40}\mathbf{D}_1+\frac18(\mathbf{G}^{{\tilde{R}}_1}+C_1 \mathbf{G}^{\tilde{S}_2}),
\end{aligned}
\end{equation}
It should be noted that  $|\partial_{\xi} p(\tilde{V}^{\mathbf{X}})| \leq C(|\tilde{V}_{\xi}^{\tilde{S}_2}|+|\tilde{U}_{\xi}^{\tilde{R}_1,\mathbf{X}}|)$
and Lemma \ref{propoties of relative quantity} are used in the derivation of the inequality \eqref{e-b2}.
Similarly, it holds
\begin{gather}\label{4577}
\left|\mathbf{B}_3\right| \leq C \lambda \tilde{\delta}_2 \int_{\mathbb{R}}\big(|\tilde{V}_{\xi}^{\tilde{S}_2}|+|\tilde{U}_{\xi}^{\tilde{R}_1,\mathbf{X}}|\big) (\Delta P)^2 d \xi \leq \frac{1}{8}(C_1 \mathbf{G}^{\tilde{S}_2}+\mathbf{G}^{\tilde{R}_1}).
\end{gather}
The next step is to give the estimates for the terms $\mathbf{P}_i$ $(i=1,2,3)$ as following.
For $\mathbf{P}_1$ and $\mathbf{P}_2$, using \eqref{456} and Young's inequality, one has
\begin{equation}\label{458}
\begin{split}
  \left|\mathbf{P}_1\right|&\leq C\| h-\tilde h\|_{L^2}\big(\|\sigma_{\xi\xi}\|_{L^2}\|\sigma_{\xi}\|_{L^\infty}+\|\sigma_{\xi}\|_{L^2}\|\sigma_{\xi}\|_{L^\infty}^2+\| \tilde{V}^{\mathbf{X}}_\xi\|_{L^\infty}\|\sigma_{\xi}\|_{L^4}^2\big)\\
  &\leq C(\epsilon_{1}+\tilde{\delta}_{1}) \big(\|\sigma_{\xi\xi}\|_{L^2}^{3/2}\|\sigma_{\xi}\|_{L^2}^{1/2}+\|\sigma_{\xi}\|_{L^2}^2\|\sigma_{\xi\xi}\|_{L^2}+(\tilde{\delta}_2+\tilde{\delta}_1)\|\sigma_{\xi}\|_{L^2}^{3/2}\|\sigma_{\xi\xi}\|_{L^2}^{1/2}\big)\\
  &\leq C(\epsilon_{1}+\tilde{\delta}_{1})(\|\sigma_{\xi}\|_{L^2}^2+\|\sigma_{\xi\xi}\|_{L^2}^2)\leq \frac18(\mathbf{A}_1 +\|\sigma_{\xi\xi}\|_{L^2}^2).\\
  \left|\mathbf{P}_2\right|&\leq C\left\| a_\xi\right\|_{L^\infty}\left(\|\sigma\|_{L^2}^2+\|\sigma_{\xi}\|_{L^2}^2\right)+C\|\sigma\|_{L^2}^2\|\sigma\|_{L^\infty}+C\left\| \tilde{V}^{\mathbf{X}}_\xi\right\|_{L^\infty}\|\sigma\|_{L^2}\|\sigma_{\xi}\|_{L^2}\\
  &\leq C(\tilde{\delta}_2+\tilde{\delta}_1)(\|\sigma_{\xi}\|_{L^2}^2+\|\sigma\|_{L^2}^2)+C\epsilon_{1}\|\sigma\|_{L^2}^2\leq \frac14\mathbf{A}_1.
  \end{split}
\end{equation}
To give the estimate of  $\mathbf{P}_3$, we start with the following inequality of  $\phi_\xi$ in terms of $\mathbf{D}_1, \mathbf{G}_R$ and $\mathbf{G}_S$. Observing that there exists a constant $\bar{v}$ between $v$ and $\tilde{V}^{\mathbf{X}}$ such that
\begin{equation}
\big(p(v)-p(\tilde{V}^{\mathbf{X}})\big)_\xi=p^{\prime}(v)\phi_{\xi}
   -\big( p^{\prime}(v)-p^{\prime}(\tilde{V}^{\mathbf{X}}) \big)\tilde{V}^{\mathbf{X}}_{\xi}=p^{\prime}(v)\phi_{\xi}
   -p^{\prime\prime}(\bar{v})\tilde{V}^{\mathbf{X}}_{\xi}\phi,\notag
\end{equation}
which together with the algebraic inequality $\frac{p^2}{2}-q^2\leq(p-q)^2$ implies
\begin{equation}
\begin{split}
 \mathbf{D}_1&=\int_{\mathbb{R}}\frac{a}{\gamma p(v)}\left( p^{\prime}(v)\phi_{\xi}
   -\big( p^{\prime}(v)-p^{\prime}(\tilde{V}^{\mathbf{X}}) \big)\tilde{V}^{\mathbf{X}}_{\xi}\right)^{2}d\xi\\
   & \geq \frac{1}{2}\int_{\mathbb{R}}\frac{a(p^{\prime}(v))^{2}}{\gamma p(v)}|\phi_{\xi}|^{2}d\xi-\int_{\mathbb{R}}\frac{a(p^{\prime\prime}(\bar{v}))^2}{\gamma p(v)}(\tilde{V}^{\mathbf{X}}_{\xi})^{2}\phi^{2}d\xi\\
   & \geq \frac{1}{C_2}\int_{\mathbb{R}}\phi_{\xi}^{2}d\xi-C(\tilde\delta_2+\tilde\delta_1)\int_{\mathbb{R}}\left(\tilde{V}_{\xi}^{\tilde{S}_2}+\tilde{V}_{\xi}^{\tilde{R}_1,\mathbf{X}}\right)\phi^{2} d\xi\\
   & \geq \frac{1}{C_2}\int_{\mathbb{R}}\phi_{\xi}^{2}d\xi-C(\tilde\delta_2+\tilde\delta_1)\left( \mathbf{G}^{\tilde{S}_2}+\mathbf{G}^{\tilde{R}_1}\right),\notag
   \end{split}
\end{equation}
where the facts that $\phi^{2}\leq Cp(v\big| \tilde{V}^{\mathbf{X}})\leq C\big(p(v)-p(\tilde{V}^{\mathbf{X}})\big)$ are used, and it implies
\begin{equation}\label{phi-xi}
    \|\phi_{\xi}\|_{L^2}^2\leq C_2\mathbf{D}_1+C(\tilde\delta_2+\tilde\delta_1)\big( \mathbf{G}^{\tilde{S}_2}+\mathbf{G}^{\tilde{R}_1}\big),
\end{equation}
which together with \eqref{458}, one gets
\begin{equation}
\begin{aligned}
 |\mathbf{P}_3|&\leq C\| h-\tilde h\|_{L^2}\|\phi_{\xi}\|_{L^2}\|\sigma_{\xi}\|_{L^\infty}^2+C\|\sigma\|_{L^\infty}\|\sigma_{\xi}\|_{L^2}\|\phi_{\xi}\|_{L^2}\\
  &\leq C(\epsilon_{1}+\tilde\delta_1) \left(\|\phi_{\xi}\|_{L^2}^{2}+\|\sigma_{\xi\xi}\|_{L^2}^{2}\right)+C\epsilon_{1}\left(\|\phi_{\xi}\|_{L^2}^{2}+\|\sigma_{\xi}\|_{L^2}^{2}\right)\\
  &\leq \frac18(\mathbf{A}_1+\|\sigma_{\xi\xi}\|_{L^2}^2)+\frac{1}{40}(\mathbf{D}_1+C_1\mathbf{G}^{\tilde{S}_2}+\mathbf{G}^{\tilde{R}_1}).\\
\end{aligned}
\end{equation}
And next,  $\mathbf{S}_1$ and $\mathbf{S}_2$  are estimated as follows. Since $\tilde{V}^{\tilde{S}_2}, \tilde{V}^{\mathbf{X}},\tilde{V}^{\tilde{R}_1,\mathbf{X}} \in\left(\frac{v_{+}}2,2 v_{+}\right)$,
one has
\begin{eqnarray}
&&\Big|(\ln \tilde{V}^{\tilde S_2}-\ln \tilde{V}^{\mathbf{X}})_{\xi \xi}\Big|\\
&&=\Big|\tilde{V}_{\xi \xi}^{\tilde S_2}(\frac{1}{\tilde{V}^{\tilde S_2}}-\frac{1}{\tilde{V}^{\mathbf{X}}})+\frac{1}{\tilde{V}^{\mathbf{X}}}(\tilde{V}_{\xi \xi}^{\tilde S_2}-\tilde{V}_{\xi \xi}^{\mathbf{X}})-\frac{1}
{{(\tilde{V}^{\tilde S_2})}^2}\big((\tilde{V}_\xi^{\tilde S_2})^2-(\tilde{V}_{\xi}^{\mathbf{X}})^2\big)\notag\\
&&\ \qquad\qquad-\big(\tilde{V}_{\xi}^{\mathbf{X}}\big)^2\Big(\frac{1}{(\tilde{V}^{\tilde S_2})^2}-\frac{1}{(\tilde{V}^{\mathbf{X}})^2}\Big)\Big| \label{459}\\
&&\leq C\left(|\tilde{V}_{\xi \xi}^{\tilde{S}_2}||\tilde{V}^{\tilde R_1,\mathbf{X}}-v_m|+|\tilde{V}^{\tilde R_1,\mathbf{X}}_{\xi \xi}|+|\tilde{V}^{\tilde R_1,\mathbf{X}}_{\xi}||\tilde{V}_{\xi}^{\tilde{S}_2}|+|\tilde{V}^{\tilde R_1,\mathbf{X}}_{\xi}|^2+|\tilde{V}_{\xi}^{\tilde{S}_2}|^2|\tilde{V}^{\tilde R_1,\mathbf{X}}-v_m|\right)\notag\\
&&\leq C\left(|\tilde{V}_\xi^{\tilde S_2}||\tilde{V}^{\tilde R_1,\mathbf{X}}-v_m|+|\tilde{V}^{\tilde R_1,\mathbf{X}}_{\xi}||\tilde{V}_\xi^{\tilde S_2}|+|\tilde{V}^{\tilde R_1,\mathbf{X}}_{\xi \xi}|+|\tilde{V}^{\tilde R_1,\mathbf{X}}_{\xi}|^2\right).\notag
\end{eqnarray}
Then, one obtains
\begin{eqnarray}
\left|\mathbf{S}_1\right|&\leq & C \int_{\mathbb{R}}| (\Delta P)|\left(|\tilde{V}_\xi^{\tilde S_2}||\tilde{V}^{\tilde R_1,\mathbf{X}}\!-\!v_m|+|\tilde{V}^{\tilde R_1,\mathbf{X}}_{\xi}||\tilde{V}_\xi^{\tilde S_2}|+|\tilde{V}^{\tilde R_1,\mathbf{X}}_{\xi \xi}|+|\tilde{V}^{\tilde R_1,\mathbf{X}}_{\xi}|^2\right)d \xi \notag\\
& \leq& C\| (\Delta P)\|_{L^2}\big\||\tilde{V}_\xi^{\tilde S_2}|\cdot|\tilde{V}^{\tilde R_1,\mathbf{X}}-v_m|+|\tilde{V}^{\tilde R_1,\mathbf{X}}_{\xi}|\cdot|\tilde{V}_\xi^{\tilde S_2}|\big\|_{L^2}+C\| \Delta P\|_{L^{\infty}}\|\tilde{V}^{\tilde R_1,\mathbf{X}}_{\xi \xi}\|_{L^1} \notag\\
&\ &+C\| \Delta P\|_{L^{2}}\|\tilde{V}^{\tilde R_1,\mathbf{X}}_{\xi}\|_{L^4}^2\label{e-s1}\\
& \leq& C(M+\tilde\delta_1)\tilde\delta_1\tilde\delta_2e^{-C\tilde\delta_2\tau}+CM^{1 / 2}|\mathbf{D}_1|^{1 / 4}\|\tilde{V}^{\tilde R_1,\mathbf{X}}_{\xi \xi}\|_{L^1}\!+\!CM\|\tilde{V}^{\tilde R_1,\mathbf{X}}_{\xi}\|_{L^4}^2\notag \\
& \leq& \frac{1}{40} \mathbf{D}_1+C(M+\tilde\delta_1)\tilde\delta_1\tilde\delta_2e^{-C\tilde\delta_2\tau}+CM^{2/3}|\|\tilde{V}^{\tilde R_1,\mathbf{X}}_{\xi \xi}\|_{L^1}^{4/3}\!+\!CM\|\tilde{V}^{\tilde R_1,\mathbf{X}}_{\xi}\|_{L^4}^2.\notag
\end{eqnarray}
Similarly, observing that
\begin{gather}
\Big|\big(p(\tilde{V}^{\mathbf{X}})-p(\tilde{V}^{\tilde R_1,\mathbf{X}})-p(\tilde{V}^{\tilde S_2})\big)_{\xi}\Big| \leq C\Big(\big|\tilde{V}^{\tilde R_1,\mathbf{X}}_{\xi} \big||\tilde{V}^{\tilde S_2}-v_m|+\big|\tilde{V}_\xi^{\tilde S_2}\big|\big|\tilde{V}^{\tilde R_1,\mathbf{X}}-v_m\big|\Big),\label{460}
\end{gather}
then, one has
\begin{equation}\label{e-s2}
\begin{aligned}
\left|\mathbf{S}_2\right|&\leq  C \int_{\mathbb{R}}|h-\tilde h|\big(|\tilde{V}^{\tilde R_1,\mathbf{X}}_{\xi} ||\tilde{V}^{\tilde S_2}-v_m|+|\tilde{V}_\xi^{\tilde S_2}||\tilde{V}^{\tilde R_1,\mathbf{X}}-v_m|\big)d \xi \\
& \leq C\|h-\tilde h\|_{L^2}\Big\||\tilde{V}^{\tilde R_1,\mathbf{X}}_{\xi} ||\tilde{V}^{\tilde S_2}-v_m|+|\tilde{V}_\xi^{\tilde S_2}||\tilde{V}^{\tilde R_1,\mathbf{X}}-v_m|\Big\|_{L^2} \\
& \leq C\big(M+\tilde\delta_1\big)\tilde\delta_1\tilde\delta_2e^{-C\tilde\delta_2\tau}.
\end{aligned}
\end{equation}
and so far, one obtains
\begin{equation}\label{wre3}
\begin{aligned}
&\frac{d}{d \tau} \int_{\mathbb{R}} a \eta(\mathbf{w} | \tilde{\mathbf{w}}^{\mathbf{X}}) d \xi+\frac{\tilde\delta_2}{4 M}|\dot{\mathbf{X}}|^2+\frac{\mathbf{G}_1}{2}+\frac{\mathbf{G}^{\tilde{R}_1}}{2}+\frac{\mathbf{D}_1}{8} +\frac{\mathbf{A}_1}{2}+\frac{C_1}{2} \mathbf{G}^{\tilde{S}_2} \\
\leq& \frac14\|\sigma_{\xi\xi}\|_{L^2}^2
+C\big(M+\tilde\delta_1\big)\tilde\delta_1\tilde\delta_2e^{-C\tilde\delta_2\tau}
+CM^{\frac23}\|\tilde{V}^{\tilde R_1,\mathbf{X}}_{\xi \xi}\|_{L^1}^{\frac43}+CM\|\tilde{V}^{\tilde R_1,\mathbf{X}}_{\xi}\|_{L^4}^2.
\end{aligned}
\end{equation}
Finally, multiplying $\eqref{nsac-per}_{3}$ by $-\sigma_{\xi\xi}$,  and integrating over $\mathbb{R}$, one has
\begin{equation}
\begin{split}
    &\frac{1}{2}\frac{d}{dt}\int_{\mathbb{R}}|\sigma_\xi|^{2}d\xi+\|\sigma_{\xi\xi}\|_{L^2}^2\\
   = &\int_{\mathbb{R}}2v(\sigma+\sigma^2)\sigma_{\xi\xi}d\xi+\int_{\mathbb{R}} \frac{\sigma_\xi v_\xi}{v}\sigma_{\xi\xi}d\xi +\int_{\mathbb{R}}\frac{\sigma_\xi^2\sigma_{\xi\xi}}{\sigma+1}d\xi\\
\leq&\!-\!\int_{\mathbb{R}}2v|\sigma_\xi|^{2}d\xi\!-\!\int_{\mathbb{R}}2(\phi_{\xi}\!+\!\tilde{V}^{\mathbf{X}}_\xi)\sigma\sigma_\xi d\xi\!+\!\int_{\mathbb{R}}2v\sigma^2\sigma_{\xi\xi}d\xi\!+\!\int_{\mathbb{R}} \frac{\sigma_\xi (\phi_{\xi}\!+\!\tilde{V}^{\mathbf{X}}_\xi)}{v}\sigma_{\xi\xi}d\xi\!+\!\int_{\mathbb{R}}\frac{\sigma_\xi^2\sigma_{\xi\xi}}{\sigma+1}d\xi\\
\leq&\!-\!\int_{\mathbb{R}}2v|\sigma_\xi|^{2}d\xi\!+\!C\|\sigma\|_{L^\infty}(\|\sigma_\xi\|_{L^2}\|\phi_\xi\|_{L^2}\!+\!\|\sigma\|_{L^2}\|\sigma_{\xi\xi}\|_{L^2})\!+\!C\|\tilde{V}^{\mathbf{X}}_\xi\|_{L^\infty}\|\sigma_\xi\|_{L^2}\|\sigma\|_{L^2})\\
&+\!C\|\sigma_\xi\|_{L^\infty}(\|\sigma_{\xi\xi}\|_{L^2}\|\phi_\xi\|_{L^2}\!+\!\|\sigma_\xi\|_{L^2}\|\sigma_{\xi\xi}\|_{L^2})\!+\!C\|\tilde{V}^{\mathbf{X}}_\xi\|_{L^\infty}\|\sigma_\xi\|_{L^2}\|\sigma_{\xi\xi}\|_{L^2}\\
\leq&\!-\!\int_{\mathbb{R}}2v|\sigma_\xi|^{2}d\xi\!+\!CM(\|\phi_\xi\|_{L^2}^2\!+\!\mathbf{A}_1)\!+\!C\sqrt{M}(\|\phi_\xi\|_{L^2}^2\!+\!\|\sigma_{\xi\xi}\|_{L^2}^2)\!+\!C(\tilde\delta_2\!+\!\tilde\delta_1)(\mathbf{A}_1\!+\!\|\sigma_{\xi\xi}\|_{L^2}^2),\notag
\end{split}
\end{equation}
which together with \eqref{phi-xi} yields
\begin{equation}\label{d22}
\begin{split}
&\frac{1}{2}\frac{d}{dt}\int_{\mathbb{R}}|\sigma_\xi|^{2}d\xi+\|\sigma_{\xi\xi}\|_{L^2}^2+\int_{\mathbb{R}}2v|\sigma_\xi|^{2}d\xi\\
&\leq\frac14\|\sigma_{\xi\xi}\|_{L^2}^2+\frac14(\mathbf{A}_1+\mathbf{G}^R+C_1\mathbf{G}^S)+\frac{1}{16}\mathbf{D}_1,
\end{split}
\end{equation}
combining with \eqref{wre3} \eqref{d22}, one gets
\begin{equation}\label{wre4}
\begin{aligned}
&\frac{d}{d \tau} \int_{\mathbb{R}} a \eta(\mathbf{w} \big| \tilde{\mathbf{w}}^{\mathbf{X}}) d \xi+\frac{1}{2}\frac{d}{d\tau}\int_{\mathbb{R}}|\sigma_\xi|^{2}d\xi+\frac{\tilde\delta_2}{4 M}|\dot{\mathbf{X}}|^2\\
&+\frac{\mathbf{G}_1}{2}+\frac{\mathbf{D}_1}{16} +\frac{\mathbf{A}_1}{4}+\frac{C_1}{4} \mathbf{G}^{\tilde{S}_2}+\frac{1}{2}\|\sigma_{\xi\xi}\|_{L^2}^2 \\
\leq& C\left(M+\tilde\delta_1\right)\tilde\delta_1\tilde\delta_2e^{-C\tilde\delta_2\tau}
+CM^{\frac23}|\|\tilde{V}^{\tilde R_1,\mathbf{X}}_{\xi \xi}\|_{L^1}^{\frac43}\!+\!CM\|\tilde{V}^{\tilde R_1,\mathbf{X}}_{\xi}\|_{L^4}^2.
\end{aligned}
\end{equation}
Noticing that by Lemma \ref{the properties of rarefaction wave},
\begin{gather}
\left\|\tilde{V}^{\tilde{R}_1}_{\xi \xi}\right\|_{L^1} \leq \begin{cases}\tilde\delta_1 & \text { if } 1+\tau \leq \tilde\delta_1^{-1}, \\ \frac{1}{1+\tau} & \text { if } 1+\tau \geq \tilde\delta_1^{-1},\end{cases},\quad
\left\|\tilde{V}^{\tilde{R}_1}_{\xi} \right\|_{L^4} \leq \begin{cases}\tilde\delta_1 & \text { if } 1+\tau \leq \tilde\delta_1^{-1}, \\ \tilde\delta_1^{1 / 4} \frac{1}{(1+\tau)^{3 / 4}} & \text { if } 1+\tau \geq \tilde\delta_1^{-1},\end{cases}\notag
\end{gather}
which implies
\begin{gather}
\int_0^{\infty}\big\|\tilde{V}^{\tilde{R}_1,{\mathbf{X}}}_{\xi \xi}\big\|_{L^1}^{4 / 3} d s \leq C \tilde\delta_1^{\frac13}, \quad \int_0^{\infty}\big\|\tilde{V}^{\tilde{R}_1,{\mathbf{X}}}_{\xi} \big\|_{L^4}^2 d s \leq C \tilde\delta_1 .\label{464}
\end{gather}
Integrating \eqref{wre4} over $[0, \tau]$ for any $\tau \leq T$ and using \eqref{464}, we have
\begin{equation}
\begin{aligned}
& \sup _{t \in[0, T]} \int_{\mathbb{R}} \left(\eta(\mathbf{w} | \tilde{\mathbf{w}}^{\mathbf{X}})\!+\!|\sigma_\xi|^{2}\right) d \xi+\tilde\delta_2 \int_0^t|\dot{\mathbf{X}}|^2 d s
+\int_0^t\left(\mathbf{G}_1\!+\!\mathbf{G}^{\tilde{S}_2}\!+\!\mathbf{D}_1\!+\!\mathbf{A}_1\!+\!\|\sigma_{\xi\xi}\|_{L^2}^2\right) d s \\
&\leq  C \int_{\mathbb{R}} \left(\eta(\mathbf{w}_0 | \tilde{\mathbf{w}}(\xi,0))+|\sigma_{0\xi}|^{2}\right) d \xi+C \tilde\delta_1^{\frac13},\notag
\end{aligned}
\end{equation}
which together with \eqref{ggg}, the desired estimate \eqref{47} is obtained, and the proof of Lemma \ref{lem-vh} is completed.
\end{proof}

\subsection{Lower order energy estimates}
\qquad
The purpose of this section is to derive the low order energy estimation of the solution for the system \eqref{nsac-per}-\eqref{init-per} after the shock waves interaction in the same group.  Since the main difficulty is to give the lower derivative estimate of the velocity for the flow, we will not use the effective velocity technique considered in the subsection 4.2.

\begin{lemma}\label{lem-low}
 Under the hypotheses of Proposition \ref{prop}, there exists $C>0$ (independent of $\left.\delta_0, M, T\right)$ such that for all $\tau \in(0, T]$,
\begin{equation}
\begin{aligned}
\|  &\phi , \sigma \|_{H^1(\mathbb{R})}^2+\| \psi\|_{L^2(\mathbb{R})}^2\\
&+\tilde\delta_2 \int_0^t|\dot{\mathbf{X}}|^2 d s +\int_0^\tau\left(G^{\tilde{S}_2}+G^{\tilde{R}_1}+D_1+\|\psi_{\xi}\|_{L^2(\mathbb{R})}^2+\|\sigma\|_{H^{2}(\mathbb{R})}^{2}\right) d s \\
\leq & C\left(\left\|\phi_{0} , \sigma_{0}\right\|_{H^1(\mathbb{R})}^2+\left\|\psi_{0}\right\|_{L^2(\mathbb{R})}^2\right)+C \tilde\delta_1^{\frac13},\label{51}
\end{aligned}
\end{equation}
where $G^{\tilde{S}_2}, D_1$ are as in \eqref{48}, and
\begin{gather}
\begin{aligned}
G^{\tilde{R}_1}  & \xlongequal{\mathrm{def}} \int_{\mathbb{R}} \tilde{U}_{\xi}^{\tilde{R}_1} p(v \big| \tilde{V}) d \xi.\label{52}
\end{aligned}
\end{gather}
\end{lemma}

\begin{proof} Here we turn to the system \eqref{nsac-per}-\eqref{init-per}. As in the proof step of Lemma 4.2,
 the following vectors and matrix operators are introduced
\begin{equation}\label{z}
\displaystyle\mathbf{z}\xlongequal{\mathrm{def}} \left(\begin{array}{l}
\displaystyle v \\
\displaystyle u \\
\omega
\end{array}\right), \
\displaystyle\mathbf{\tilde{z}}\xlongequal{\mathrm{def}} \left(\begin{array}{l}
\displaystyle \tilde{V} \\
\displaystyle \tilde{U} \\
1
\end{array}\right),\
 M_1(\mathbf{z})\xlongequal{\mathrm{def}} \left(\begin{array}{ccc}
\displaystyle 0 & \displaystyle 0 &\displaystyle  0\\
\displaystyle 0 & \displaystyle\frac{1}{v}&\displaystyle -\frac{\omega_{\xi}}{8v^2\omega}\\ \displaystyle 0 & \displaystyle 0&\displaystyle 1
\end{array}\right),
\end{equation}
moreover, similarly, the following entropy function is introduced
\begin{equation}\label{entropy function of z}
  \eta(\mathbf{z})=Q(v)+\frac{u^2}{2}+\frac{\omega^2}{2},
\end{equation}
where $Q$ satisfies \eqref{Q}, and the relative entropy of $\eta$ is as follows
\begin{equation}\label{425-2}
\begin{aligned}
\eta(\mathbf{z}\big| \tilde{\mathbf{z}})  =Q(v \big| \tilde{V})+\frac{ \psi ^2}{2}+\frac{\sigma^2}{2}.
\end{aligned}
\end{equation}
From the definition \eqref{z} and \eqref{entropy function of z},   the system $\eqref{NSAC-xi}$ and the system \eqref{equation of approximate composite shock wave after interaction} are transformed into the following  vector form equations respectively
\begin{gather}
\partial_\tau \mathbf{z}+\partial_{\xi} A(\mathbf{z})=\partial_{\xi}\big(M_1(\mathbf{z}) \partial_{\xi} \nabla \eta(\mathbf{z})\big)+H(\mathbf{z}).\label{53}
\end{gather}
and
\begin{equation}
\partial_t \tilde{\mathbf{z}}+\partial_{\xi} A(\tilde{\mathbf{z}})=\partial_{\xi}\left(M_1(\tilde{\mathbf{z}}) \partial_{\xi} \nabla \eta(\tilde{\mathbf{z}})\right)-\dot{\mathbf{X}}\tilde{\mathbf{z}}^{\tilde S_2,{-\mathbf{X}}}_{\xi}+\left(\begin{array}{c}
0 \\
F_1+F_2\\0
\end{array}\right),\label{54}
\end{equation}
where $F_1, F_2$ are as in \eqref{310}. Adopting the same method for the estimates of the weighted relative entropy in \eqref{wre1} with $a\equiv 1$, one has
\begin{gather}\label{re1}
\frac{d}{d \tau} \int_{\mathbb{R}} \eta(\mathbf{z}\big| \tilde{\mathbf{z}}) d \xi=\dot{\mathbf{X}}\mathcal{Y}(\mathbf{z})+\mathcal{B}(\mathbf{z})-\mathcal{G}(\mathbf{z}),
\end{gather}
where
\begin{equation}
\begin{aligned}
\mathcal{Y}(\mathbf{z}) & =\underbrace{-\int_{\mathbb{R}} p^{\prime}(\tilde{V}) \tilde{V}_{\xi}^{\tilde S_2,{-\mathbf{X}}}\phi d \xi}_{\mathcal{Y}_1}+\underbrace{\int_{\mathbb{R}} \tilde{U}_{\xi}^{\tilde S_2,{-\mathbf{X}}}\psi d \xi}_{\mathcal{Y}_2},\\
\mathcal{B}(\mathbf{z}) & =\underbrace{-\int_{\mathbb{R}} \tilde{U}^{\tilde S_2,{-\mathbf{X}}}_{\xi}  p(v \big| \tilde{V}) d \xi}_{\mathbb{B}_1} \underbrace{-\int_{\mathbb{R}} (\frac{1}{v}-\frac{1}{\tilde{V}})\tilde{U}_\xi\psi_\xi d \xi}_{\mathbb{B}_2} +\underbrace{\frac{1}{8}\int_{\mathbb{R}} \frac{\sigma_{\xi}^{2}\psi_\xi}{(\sigma+1) v^{2}} d \xi}_{\mathbb{B}_3}\\
& \underbrace{-\int_{\mathbb{R}} \psi(F_1+F_2)   d \xi}_{\mathbb{B}_4}\underbrace{-\int_{\mathbb{R}} \Big(2v\sigma^3+\frac{(\phi_{\xi}+\tilde{V}_{\xi})\sigma_{\xi}\sigma}{v}+\frac{\sigma_{\xi}^{2}\sigma}{\sigma+1}\Big)  d \xi}_{\mathbb{B}_5} ,\notag
\end{aligned}
\end{equation}
and
\begin{equation}
    \begin{aligned}
\mathcal{G}(\mathbf{z}) & =\underbrace{\int_{\mathbb{R}} \frac{1}{v}\left|\psi_\xi\right|^2 d \xi}_{\mathbb{E}_1}+\underbrace{\int_{\mathbb{R}} \tilde{U}_{\xi}^{\tilde{R}_1} p(v | \tilde{V}) d \xi}_{G^{\tilde{R}_1}}+\underbrace{\int_{\mathbb{R}} (2v\sigma^2+|\sigma_\xi|^2) d \xi}_{\mathbb{A}_1}.\notag
\end{aligned}
\end{equation}
Noticing that $\mathcal{G}(\mathbf{z})$ consists of good terms, while $\mathcal{B}(\mathbf{z})$ consists bad terms related to the perturbation $\phi,\psi$ and $\sigma$.
The good terms $G^{\tilde{R}_1}, \mathbb{A}_1$ and $\mathbb{E}_1$ will be used to control the bad terms.
The proof start  with the energy estimation of $\mathcal{Y}_1$ and $\mathcal{Y}_2$. By using Lemma \ref{lem-vh}, one deduces from \eqref{phi-xi} that
\begin{equation}
\begin{aligned}
\left|\mathcal{Y}_1\right| \leq & C\left(\int_{\mathbb{R}} |\tilde{V}_{\xi}^{\tilde S_2,{-\mathbf{X}}} | d \xi\right)^{1/2} \left(\int_{\mathbb{R}} |\tilde{V}_{\xi}^{\tilde S_2,{-\mathbf{X}}} | \phi^2 d \xi\right)^{1/2} \leq C\sqrt{\tilde\delta_2} \sqrt{G^{\tilde{S}_2}},\\
\left|\mathcal{Y}_2\right| \leq & C \int |\tilde{U}_{\xi}^{\tilde S_2,{-\mathbf{X}}} |\Big(\Big| h-\tilde h -\frac{p(v)-p(\tilde{V})}{\tilde{s}_2}\Big|+|p(v)-p(\tilde{V})|\\
& +\left|\phi_{\xi}\right|+\left|\tilde{V}_{\xi}\right|\left|\phi\right|+|\tilde{V}_{\xi}^{\tilde S_2,{-\mathbf{X}}}|\left|\tilde{V}^{\tilde{R}_1}-v_m\right|+\left|\tilde{V}_{\xi}^{\tilde{R}_1}\right|\Big) d \xi \\
\leq & C\Big(\frac{\tilde\delta_2}{\sqrt{\lambda}} \sqrt{G_1}+\sqrt{\tilde\delta_2} \sqrt{G^{\tilde{S}_2}}+\tilde\delta_2 \sqrt{D_1}+\tilde\delta_2 \sqrt{G^{\tilde{R}_1}}+\tilde\delta_2 \tilde\delta_1 e^{-C \tilde\delta_2 \tau}\Big) ,\notag
\end{aligned}
\end{equation}
where in the estimates of $\mathcal{Y}_2$, we have use the fact (as done in \eqref{455}):
\begin{gather}
|\psi| \leq \big|h-\tilde h\big|+C\Big(|\phi_{\xi}|+|\tilde{V}_{\xi}|\cdot|\phi|+|\tilde{V}_{\xi}^{\tilde S_2,{-\mathbf{X}}}||\tilde{V}^{\tilde{R}_1}-v_m|+|\tilde{V}_{\xi}^{\tilde{R}_1}|\Big).\notag
\end{gather}
Thus,
\begin{gather}
|\dot{\mathbf{X}}\mathcal{Y}(\mathbf{z})|\leq \frac{\tilde\delta_2}{2}|\dot{\mathbf{X}}|^2+C \frac{\tilde\delta_2}{\lambda} G_1+C G^{\tilde{S}_2}+C \tilde\delta_2 D_1+C \tilde\delta_2 G^{\tilde{R}_1}+C \tilde\delta_2 \tilde\delta_1^2 e^{-C \tilde\delta_2 \tau} .\label{y1y2}
\end{gather}
Next, we will give the energy estimates of bad terms $\mathcal{B}(\mathbf{z})$. For $\mathbb{B}_1,\mathbb{B}_2$ and $\mathbb{B}_3$, using Lemma \ref{the properties of rarefaction wave}, Lemma \ref{estimate of viscous shock} and Young's inequality, we have
\begin{align*}
\left|\mathbb{B}_1\right| &\leq C G^{\tilde{S}_2},\\
\left|\mathbb{B}_2\right| &\leq \int_{\mathbb{R}}\left|\psi_{\xi}\right||\phi|\big(|\tilde{U}_{\xi}^{\tilde{R}_1}|+|\tilde{U}^{\tilde S_2,{-\mathbf{X}}}_{\xi}|\big) d \xi \leq \frac{1}{8} \mathbb{E}_1+C \tilde\delta_1 G^{\tilde{R}_1}+C \tilde\delta_2 G^{\tilde{S}_2}, \\
\left|\mathbb{B}_3\right| &\leq C\int_{\mathbb{R}}\left|\psi_{\xi}\right||\sigma_\xi|^2d \xi\leq C\|\sigma_\xi\|_{L^\infty}\|\sigma_\xi\|_{L^2}\mathbb{E}_1^{1/2}\leq \frac{1}{4} \mathbb{E}_1+\frac{1}{4} \mathbb{A}_1.
\end{align*}
For $\mathbb{B}_4$, using Lemma \ref{the properties of rarefaction wave}, Lemma \ref{estimate of viscous shock} again, as done in \eqref{459}, one gets
\begin{eqnarray}
|F_1| &\leq& C\big(\big|\tilde{U}^{\tilde{R}_1}_{\xi \xi}\big|+|\tilde{U}^{\tilde{R}_1}_{\xi}||\tilde{V}^{\tilde{R}_1}_{\xi}|+(| \tilde{U}^{\tilde S_2,{-\mathbf{X}}} _{\xi \xi} |+|\tilde{U}^{\tilde S_2,{-\mathbf{X}}}_{\xi}||\tilde{V}^{\tilde S_2,{-\mathbf{X}}}_{\xi}|)|\tilde{V}^{\tilde{R}_1}-v_m|\notag\\
&&+|\tilde{U}^{\tilde{R}_1}_{\xi}||\tilde{V}^{\tilde S_2,{-\mathbf{X}}}_{\xi}|+|\tilde{V}^{\tilde{R}_1}_{\xi}||\tilde{U}^{\tilde S_2,{-\mathbf{X}}}_{\xi}|\big)\label{55} \\
&\leq& C\big(\big|\tilde{U}^{\tilde{R}_1}_{\xi \xi}\big|+\big|\tilde{U}^{\tilde{R}_1}_{\xi}\big|^2+(|\tilde{V}^{\tilde S_2,{-\mathbf{X}}}_{\xi \xi}|+|\tilde{V}^{\tilde S_2,{-\mathbf{X}}}_{\xi}|^2)|\tilde{V}^{\tilde{R}_1}-v_m|+|\tilde{V}^{\tilde{R}_1}_{\xi}||\tilde{V}^{\tilde S_2,{-\mathbf{X}}}_{\xi}|\big),\notag
\end{eqnarray}

which together with \eqref{460}, using the same estimates as in \eqref{e-s1} with \eqref{key Lemma of the estimate for initial data} yields
\begin{equation}
\begin{aligned}
|\mathbb{B}_4|\leq & C \int_{\mathbb{R}}|\psi|\big(\big|\tilde{U}^{\tilde{R}_1}_{\xi \xi}\big|+\big|\tilde{U}^{\tilde{R}_1}_{\xi}\big|^2\big) d \xi \\
& +C \int_{\mathbb{R}}|\psi|\big(|\tilde{V}^{\tilde S_2,{-\mathbf{X}}}_{\xi}||\tilde{V}^{\tilde{R}_1}-v_m|+|\tilde{V}^{\tilde{R}_1}_{\xi}||\tilde{V}^{\tilde S_2}-v_m|+|\tilde{V}^{\tilde{R}_1}_{\xi}||\tilde{V}^{\tilde S_2,{-\mathbf{X}}}_{\xi}|\big) d \xi \\
\leq& C\|\psi\|_{L^2}^{1 / 2}\left\|\psi_{\xi}\right\|_{L^2}^{1 / 2}\big\|\tilde{U}^{\tilde{R}_1}_{\xi \xi}\big\|_{L^1}+C\|\psi\|_{L^2}\big\|\tilde{U}^{\tilde{R}_1}_{\xi}\big\|_{L^4}^2+C M\tilde\delta_1\tilde\delta_2e^{-C\tilde\delta_2\tau} \\
\leq& \frac{1}{8} \mathbb{E}_1+M^{2/3}\big\|\tilde{U}^{\tilde{R}_1}_{\xi \xi}\big\|_{L^1}^{\frac43}+CM\big\|\tilde{U}^{\tilde{R}_1}_{\xi}\big\|_{L^4}^2+C M\tilde\delta_1\tilde\delta_2e^{-C\tilde\delta_2\tau}.\notag
\end{aligned}
\end{equation}
For $\mathbb{B}_5$, by using Young's inequality, one obtains
\begin{equation}
\begin{aligned}
|\mathbb{B}_5|&\leq C\|\sigma\|_{L^\infty}\left(\|\sigma\|_{L^2}^2+\|\sigma_{\xi}\|_{L^2}^2\right)+\|\sigma\|_{L^\infty}\|\sigma_{\xi}\|_{L^2}\|\phi_\xi\|_{L^2}+C\big\| \tilde{V}^{\mathbf{X}}_\xi\big\|_{L^\infty}\|\sigma\|_{L^2}\|\sigma_{\xi}\|_{L^2}\\
  &\leq C(M+\tilde\delta_2+\tilde\delta_1)\left(\|\sigma\|_{L^2}^2+\|\sigma_{\xi}\|_{L^2}^2\right)+CM\|\phi_\xi\|_{L^2}^2\\
&\leq \frac14\mathbb{A}_1+CM(D_1+G^{\tilde{S}_2}+G^{\tilde{R}_2}).\notag
\end{aligned}
\end{equation}
Therefore, from the above estimates, we find that for some constant $c_1>0$,
\begin{eqnarray}
&&\frac{d}{d \tau} \int_{\mathbb{R}} \eta(\mathbf{z} | \tilde{\mathbf{z}})  d \xi+\frac{1}{2} G^R+\frac{1}{2} \mathbb{E}_1+\frac{1}{2} \mathbb{A}_1 \notag\\
&& \leq \frac{\tilde\delta_2}{2}|\dot{\mathbf{X}}|^2+C \frac{\tilde\delta_2}{\lambda} G_1+c_1 G^S+C (\tilde\delta_2+M) D_1\\
&&\quad+M^{\frac23}\big\|\tilde{U}^{\tilde{R}_1}_{\xi \xi}\big\|_{L^1}^{\frac43}+CM\big\|\tilde{U}^{\tilde{R}_1}_{\xi}\big\|_{L^4}^2+C M\tilde\delta_1\tilde\delta_2e^{-C\tilde\delta_2\tau}.\notag
\end{eqnarray}
Integrating the above inequality over $[0, \tau]$ for any $\tau \leq T$, and using \eqref{464}, one derives
\begin{equation}
\begin{aligned}
\int_{\mathbb{R}} & \left(\frac{\psi^2}{2}+\frac{\sigma^2}{2}+Q(v | \tilde{V})\right) d \xi+\frac{1}{2} \int_0^\tau\left(G^R+\mathbb{E}_1+\mathbb{A}_1\right) d s \\
\leq & \int_{\mathbb{R}}\Big(\frac{\psi_{0}^2}{2}+\frac{\sigma_{0}^2}{2}+Q\big(v_0 \big| \tilde{V}(\xi,0)\big)\Big) d \xi +C \tilde\delta_1^{\frac13}\\
& +\int_0^t\Big(\frac{\tilde\delta_2}{2}|\dot{\mathbf{X}}|^2+C \frac{\tilde\delta_2}{\lambda} G_1+c_1 G^S+C (\tilde\delta_2+M) D_1\Big) d s.\label{57}
\end{aligned}
\end{equation}
Therefore, multiplying \eqref{57} by the constant $\frac {1}{2 \max \{1, c_1\}}$, and then adding the result to \eqref{47}, together with the smallness of $\tilde\delta_2 / \lambda, \tilde\delta_2, \tilde\delta_1, M$, we have
\begin{equation}
\begin{aligned}
& \|\phi,\psi,h-\tilde h\|_{L^2(\mathbb{R})}^2+\|\sigma\|_{H^{1}(\mathbb{R})}^{2}+\tilde\delta_2 \int_0^\tau|\dot{\mathbf{X}}|^2 d s\\
&\qquad+\int_0^\tau\left(G^{\tilde{R}_1}+G^{\tilde{S}_2}+D_1+\|\psi_{\xi}\|_{L^2(\mathbb{R})}^2+\|\sigma\|_{H^{2}(\mathbb{R})}^{2}\right) d s \\
& \leq C( \|\phi_{0},\psi_{0},(h-\tilde h)_{0}\|_{L^2(\mathbb{R})}^2+\|\sigma_{0}\|_{H^{1}(\mathbb{R})}^{2})+C \tilde\delta_1^{\frac13},\label{58}
\end{aligned}
\end{equation}
where we have used the fact that $ Q(v | \tilde{V}) \sim \phi^2$ and $\|\psi_{\xi}\|_{L^2}^2 \sim \mathbb{E}_1$.
At the end of the proof,  the energy estimate of $\left\|\phi_{\xi}\right\|_{L^2(\mathbb{R})}$ is established.
By using the definition of $h$ and $\tilde{h}$ in \eqref{effective velocity}, one has
\begin{eqnarray}
\psi-(h-\tilde h)&=&(u-\tilde{U})-(h-\tilde{h})=\left(\ln v-\ln \tilde{V}^{\tilde S_2,{-\mathbf{X}}}\right)_{\xi}\notag\\
&=&\frac{\big(v-\tilde{V}^{\tilde S_2,{-\mathbf{X}}}\big)_{\xi}}{v}+\frac{\tilde{V}^{\tilde S_2,{-\mathbf{X}}}_{\xi} \big(\tilde{V}^{\tilde S_2,{-\mathbf{X}}}-v\big)}{v\tilde{V}^{\tilde S_2,{-\mathbf{X}}}},\notag
\end{eqnarray}
which yields
\begin{eqnarray}
\phi_{\xi} & =&\big(v-\tilde{V}^{\tilde S_2,{-\mathbf{X}}}\big)_{\xi}-\left(\tilde{V}-\tilde{V}^{\tilde S_2,{-\mathbf{X}}}\right)_{\xi}\notag \\
&=&v\psi-v(h-\tilde h)+\frac{\tilde{V}^{\tilde S_2,{-\mathbf{X}}}_{\xi} \big(\phi+(\tilde{V}^{\tilde{R}_1}-v_m)\big)}{\tilde{V}^{\tilde S_2,{-\mathbf{X}}}}-\tilde{V}_{\xi}^{\tilde R_1},\notag
\end{eqnarray}
together with Lemma \ref{estimate of viscous shock} and Lemma \ref{the properities for Shift X} implies
\begin{gather} \left\|\phi_{\xi}\right\|_{L^2(\mathbb{R})}^2 \leq C\left(\|\phi,\psi,h-\tilde h\|_{L^2(\mathbb{R})}^2+\tilde\delta_1^2\right).\label{59}
\end{gather}
As in \eqref{455}, one deduces from Lemma \ref{the properties of rarefaction wave} and Lemma \ref{estimate of viscous shock} that
\begin{eqnarray}
\|(h-\tilde h)(0)\|_{L^2(\mathbb{R})}^2 &\leq& C\Big(\big\|\phi_{0}\big\|_{H^1(\mathbb{R})}^2+\big\|\psi_{0}\big\|_{L^2(\mathbb{R})}^2+\tilde\delta_1^2\big\|\tilde{V}_\xi^{\tilde S_2}(0)\big\|_{L^2(\mathbb{R})}^2+\big\|\tilde{V}_{\xi}^{\tilde R_1}(0)\big\|_{L^2(\mathbb{R})}^2\Big)\notag\\
&\leq& C\left(\left\|\phi_{0}\right\|_{H^1(\mathbb{R})}^2+\left\|\psi_{0}\right\|_{L^2(\mathbb{R})}^2+\tilde\delta_1^2\right).\label{510}
\end{eqnarray}
Hence, the combination of \eqref{58}, \eqref{59} and \eqref{510} implies the desired estimate \eqref{51} and the proof of Lemma \ref{lem-low} is completed.
\end{proof}

\subsection{Higher order energy estimates}
\indent\qquad
 The purpose of this section is to derive the higher order energy estimation of the solution for the system \eqref{nsac-per}-\eqref{init-per} after the shock waves interaction in the same group.
\begin{lemma}\label{lem-uh}
  Under the hypotheses of Proposition \ref{prop}, there exist $C_1, C>0$ (independent of $\left.\delta_0, M, T\right)$ such that for all $t \in(0, T]$,
\begin{equation}
\begin{aligned}
& \|\phi,\psi,\sigma\|_{H^1(\mathbb{R})}^2+\tilde\delta_2 \int_0^\tau|\dot{\mathbf{X}}|^2 d s +\int_0^\tau\left(G^{\tilde{S}_2}+G^{\tilde{R}_1}+D_1+\|\psi_{\xi}\|_{H^1}^2+\|\sigma\|_{H^{2}}^{2}\right) d s \\
& \leq C\| \phi_{0},\psi_{0},\sigma_{0}\|_{H^1(\mathbb{R})}^2+C \tilde\delta_1^{\frac13},\label{5111}
\end{aligned}
\end{equation}
where $G^{\tilde{S}_2}, D_1$ are as in \eqref{48}, and $G^{\tilde{R}_1}$ are as in \eqref{52}.
\end{lemma}
\begin{proof}
Multiplying the system $\eqref{nsac-per}_{2}$ by $-\psi_{\xi \xi}$ and integrating the result with respect to $\xi$ over $\mathbb{R}$,  one has
\begin{equation}
\begin{aligned}
& \frac{d}{d \tau} \Big(\int_{\mathbb{R}} \frac{\left|\psi_{\xi}\right|^2}{2} d \xi\Big)+\int_{\mathbb{R}} \frac{1}{v}\left|\psi_{\xi \xi}\right|^2 d \xi \\
&=-\dot{\mathbf{X}}(t) \int_{\mathbb{R}}\tilde{U}^{\tilde S_2,{-\mathbf{X}}}_{\xi}  \psi_{\xi \xi} d \xi+\int_{\mathbb{R}}(p(v)-p(\tilde{V}))_{\xi} \psi_{\xi \xi} d \xi+\int_{\mathbb{R}}\frac{\psi_\xi (\phi_\xi+\tilde{V}_\xi)}{v^2}\psi_{\xi \xi} d \xi \\
&\quad -\int_{\mathbb{R}}\Big(\frac{\tilde{U}_{\xi}\phi}{v\tilde{V}}\Big)_{\xi} \psi_{\xi \xi} d \xi+\int_{\mathbb{R}}\left(F_1+F_2\right) \psi_{\xi \xi} d \xi-\frac{1}{8}\int_{\mathbb{R}}\Big(\frac{\sigma_{\xi}^{2}}{(\sigma+1)v}\Big)_{\xi}\psi_{\xi \xi} d \xi \xlongequal{\mathrm{def}}  \sum_{i=1}^6 J_i.
\end{aligned}
\end{equation}
First, by using Young's inequality, Lemma \ref{the properties of rarefaction wave}, Lemma \ref{estimate of viscous shock} and \eqref{key Lemma of the estimate for initial data}, one has
\begin{equation}\begin{aligned}
&\left|J_1\right| \leq |\dot{\mathbf{X}}| \tilde\delta_2^2 \int_{\mathbb{R}}e^{-C\tilde\delta_{2}|\xi-\mathbf{X}|}\left|\psi_{\xi \xi}\right| d \xi, \\
& \leq \frac{\tilde\delta_2}{2}|\dot{\mathbf{X}}|^2+C \tilde\delta_2^3 \int_{\mathbb{R}} \frac{1}{v}\left|\psi_{\xi \xi}\right|^2 d \xi \leq \frac{\tilde\delta_2}{2}|\dot{\mathbf{X}}|^2+\frac{1}{12} \int_{\mathbb{R}} \frac{1}{v}\left|\psi_{\xi \xi}\right|^2 d \xi, \\
&  \left|J_2\right| \leq \frac{1}{12} \int_{\mathbb{R}} \frac{1}{v}\left|\psi_{\xi \xi}\right|^2 d \xi+C D_1, \notag\\
&\left|J_{3}\right| \leq\left\|\phi_{\xi}\right\|_{L^2(\mathbb{R})}\left\|\psi_{\xi}\right\|_{L^{\infty}(\mathbb{R})}\big\|\psi_{\xi \xi}\big\|_{L^2(\mathbb{R})}+\big\|\tilde{V}_{\xi}\big\|_{L^{\infty}(\mathbb{R})}\left\|\psi_{\xi}\right\|_{L^2(\mathbb{R})}\left\|\psi_{\xi \xi}\right\|_{L^2(\mathbb{R})}\\
& \qquad \leq C M\left\|\psi_{\xi}\right\|_{L^2(\mathbb{R})}^{\frac12}\left\|\psi_{\xi \xi}\right\|_{L^2(\mathbb{R})}^{\frac12}\left\|\psi_{\xi \xi}\right\|_{L^2(\mathbb{R})}+C\big(\tilde\delta_2+\tilde\delta_1\big)\left\|\psi_{\xi}\right\|_{L^2(\mathbb{R})}\left\|\psi_{\xi \xi}\right\|_{L^2(\mathbb{R})} \\
&\qquad\leq \frac{1}{12} \int_{\mathbb{R}} \frac{1}{v}\left|\psi_{\xi \xi}\right|^2 d \xi+C\big(M+\tilde\delta_2+\tilde\delta_1\big) \left\|\psi_{\xi}\right\|_{L^2(\mathbb{R})}^2,
\end{aligned}
\end{equation}
\begin{equation}
\begin{aligned}
\left|J_{4}\right| &\leq C \int_{\mathbb{R}}\big(|\tilde{U}_{\xi}^{\tilde{S}_2}|+|\tilde{U}_{\xi}^{\tilde{R}_1}|\big)\left(|\phi|+\left|\phi_{\xi}\right|\right)\left|\psi_{\xi \xi}\right| d \xi \\
& \leq \frac{1}{12} \int_{\mathbb{R}} \frac{1}{v}\left|\psi_{\xi \xi}\right|^2 d \xi+C\left(\tilde\delta_2+\tilde\delta_1\right)\left(G^{\tilde{S}_2}+G^{\tilde{R}_1}+D_1\right) .\notag
\end{aligned}
\end{equation}
By using \eqref{460}, \eqref{55} and Lemma \ref{lem-vh}, a straightforward computation yields
\begin{equation}
\begin{aligned}
\left|J_5\right| & \leq C\left\|\psi_{\xi \xi}\right\|_{L^2(\mathbb{R})}\cdot\Big\|\big|\tilde{U}^{\tilde{R}_1}_{\xi \xi}\big|+\big|\tilde{U}^{\tilde{R}_1}_{\xi}\big|^2+\big(|\tilde{V}^{\tilde S_2,{-\mathbf{X}}}_{\xi\xi}|+|\tilde{V}^{\tilde S_2,{-\mathbf{X}}}_{\xi}|^2\big)\big|\tilde{V}^{\tilde{R}_1}-v_m\big|\\
&\qquad\qquad\qquad\qquad\qquad+|\tilde{V}^{\tilde{R}_1}_{\xi} ||\tilde{V}^{\tilde S_2,{-\mathbf{X}}}_{\xi}  |\Big\|_{L^2(\mathbb{R})} \\
& \leq \frac{1}{12}\int_{\mathbb{R}} \frac{1}{v}\left|\psi_{\xi \xi}\right|^2 d \xi+C\big\|\tilde{U}^{\tilde{R}_1}_{\xi \xi}\big\|_{L^2(\mathbb{R})}^2+C\big\|\tilde{U}^{\tilde{R}_1}_{\xi}\big\|_{L^4(\mathbb{R})}^4+C\tilde\delta_1\tilde\delta_2^{\frac32}e^{-c\tilde\delta_2\tau} .\notag
\end{aligned}
\end{equation}
For $J_{6}$, by using Young's inequality, one  gets
\begin{equation}
    \begin{split}
 |J_{6}|&\leq C\int_{\mathbb{R}}\big(|\sigma_{\xi}\sigma_{\xi \xi}|+|\sigma_{\xi}|^3+|\sigma_{\xi}|^2|\phi_{\xi}|+|\sigma_{\xi}|^2|\tilde{V}_{\xi}|\big)\left|\psi_{\xi \xi}\right| d \xi \\
 &\leq \frac{1}{12}\int_{\mathbb{R}} \frac{1}{v}\left|\psi_{\xi \xi}\right|^2 d \xi+CM\left(\|\sigma_{\xi}\|_{H^1(\mathbb{R})}^2+G^{\tilde{S}_2}+G^{\tilde{R}_1}+D_1\right).\notag
    \end{split}
\end{equation}
Therefore,  for some $c_2>0$, one obtains
\begin{equation}
\begin{aligned}
&\frac{1}{2}\frac{d}{d t} \int_{\mathbb{R}} \psi_{\xi}^{2} d \xi+\frac{1}{2} \int_{\mathbb{R}} \frac{1}{v}\left|\psi_{\xi \xi}\right|^2 d \xi\\
 &\leq  \frac{\tilde\delta_2}{2}|\dot{\mathbf{X}}|^2+c_2 D_1+C\left(M+\tilde\delta_2+\tilde\delta_1\right)\left(G^{\tilde{S}_2}+G^{\tilde{R}_1}+\|\psi_{\xi}\|_{L^2(\mathbb{R})}^2+\|\sigma_{\xi}\|_{H^1(\mathbb{R})}^2\right) \\
& \quad +C\big\|\tilde{U}^{\tilde{R}_1}_{\xi \xi}\big\|_{L^2(\mathbb{R})}^2+C\big\|\tilde{U}^{\tilde{R}_1}_{\xi}\big\|_{L^4(\mathbb{R})}^4+C\tilde\delta_1\tilde\delta_2^{\frac32}e^{-c\tilde\delta_2\tau} .\notag
\end{aligned}
\end{equation}
Integrating the above estimate over $[0, \tau]$ for any $\tau \leq T$, one has
\begin{equation}\label{u-h}
\begin{aligned}
&\frac{1}{2}\int_{\mathbb{R}} \psi_{\xi}^{2} d \xi +\frac{1}{2}\int_{0}^{\tau} \int_{\mathbb{R}} \frac{1}{v}\left|\psi_{\xi \xi}\right|^2 d \xi ds\\&
\leq  \frac{1}{2}\int_{\mathbb{R}} \psi_{0\xi}^{2} d \xi+C \tilde\delta_1 +\int_0^\tau\big(\frac{\tilde\delta_2}{2}|\dot{\mathbf{X}}|^2+c_2D)ds\\
& \qquad+C(M+\tilde\delta_2+\tilde\delta_1)\int_0^\tau(G^{\tilde{S}_2}+G^{\tilde{R}_1}+\|\psi_{\xi}\|_{L^2(\mathbb{R})}^2+\|\sigma_{\xi}\|_{H^1(\mathbb{R})}^2) d s,
\end{aligned}
\end{equation}
where the Lemma 3.2 is used to get
\begin{gather}
\int_0^{\infty}\big\|\tilde{U}^{\tilde{R}_1}_{\xi \xi}\big\|_{L^2(\mathbb{R})}^2 d s \leq C \tilde\delta_1, \quad \int_0^{\infty}\big\|\tilde{U}^{\tilde{R}_1}_{\xi}\big\|_{L^4(\mathbb{R})}^4 d s \leq C \tilde\delta_1^3.\label{l2l4}
\end{gather}
Multiplying \eqref{u-h} by the constant $\frac{1}{2 \max \left(1, c_2\right)}$, and then adding the result to \eqref{51}, together with the smallness of $M, \tilde\delta_2, \tilde\delta_1$, we have
\begin{equation}
\begin{aligned}
&\|  \phi, \psi,\sigma\|_{H^1(\mathbb{R})}^2+\tilde\delta_2 \int_0^\tau|\dot{\mathbf{X}}|^2 d s+\int_0^\tau\left(G^{\tilde{R}_1}+G^{\tilde{S}_2}+D_1+\|\psi_\xi\|_{H^{1}(\mathbb{R})}^{2}+\|\sigma\|_{H^{2}(\mathbb{R})}^{2}\right) d s \\
& \leq C\| \phi_{0}, \psi_{0},\sigma_{0}\|_{H^1(\mathbb{R})}^2+C \tilde\delta_1^{\frac13},\notag
\end{aligned}
\end{equation}
and thus the proof of Lemma \ref{lem-uh} is finished.
\end{proof}
\vskip 0.3cm
\begin{lemma}\label{lem-vp}
Under the hypotheses of Proposition \ref{prop}, there exists $ C>0$ (independent of $\left.\delta_0, M, T\right)$ such that for all $\tau \in(0, T]$,
\begin{equation}
    \begin{split}
       & \|(\phi,\sigma)\|_{H^{2}(\mathbb{R})}^{2}+\|\psi\|_{H^{1}(\mathbb{R})}^{2}+\tilde\delta_{2}\int_{0}^{\tau}|\dot{\mathbf{X}}|^{2}ds\\        &\quad
       +\int_{0}^{t}\left(G^{\tilde{R}_1}+G^{\tilde{S}_2}+D_1+\|\phi_{\xi\xi}\|_{L^{2}(\mathbb{R})}^{2}+\|\psi_\xi\|_{H^{1}}^{2}+\|\sigma\|_{H^{3}(\mathbb{R})}^{2}\right)ds\\
        &\leq C(\|\phi_{0},\sigma_{0}\|_{H^{2}}^{2}+\|\psi_{0}\|_{H^{1}}^{2})+C\tilde\delta_1^{\frac{1}{3}}.\label{5144}
    \end{split}
\end{equation}
\end{lemma}
\begin{proof}
Differentiating the system $\eqref{nsac-per}_{2}$ with respect to $\xi$ yields
\begin{equation}\begin{aligned}
\big(\frac{ \psi_{\xi}}{v}\big)_{\xi\xi}&=
\psi _{\tau\xi}-\tilde{s}_2 \psi _{\xi\xi}+\left(p (v)-p (\tilde{V}) \right)_{\xi\xi}+\big(\frac{ \phi\tilde{U} _{\xi}}{v\tilde{V}}\big)_{\xi\xi}\\
&\quad+\frac{1}{8}\Big(\frac{\sigma_{\xi}^{2}}{(\sigma\!+\!1)v^{2}}\Big)_{\xi\xi}-\dot{\mathbf{X}}(\tau)\tilde{U}^{\tilde S_2,{-\mathbf{X}}}_{\xi\xi}+{F_1}_{\xi}+{F_2}_{\xi}.\label{hh1}
\end{aligned}\end{equation}
By using \eqref{nsac-per}$_{1}$, one has the following fact
\begin{equation}
    \begin{split}
\left(\frac{\psi_{\xi}}{v}\right)_{\xi\xi}
&=\left( \frac{\phi_{\tau}}{v} \right)_{\xi\xi}-\tilde{s}_2 \left( \frac{\phi_{\xi}}{v} \right)_{\xi\xi}-\dot{\mathbf{X}}(\tau)\Big(
\frac{\tilde{V}_{\xi}^{\tilde S_2,-{\mathbf{X}}}}{v} \Big)_{\xi\xi}\\
&=\left( \frac{\phi_{\xi\xi}}{v} \right)_{\tau}+\frac{\phi_{\xi\xi}v_{\tau}-2\phi_{\xi \tau}v_{\xi}-\phi_{\tau}v_{\xi\xi}  }{v^{2}}+\frac{2\phi_{\tau}v_{\xi}^{2}}{v^{3}}-\tilde{s}_2 \left( \frac{\phi_{\xi}}{v} \right)_{\xi\xi}-\dot{\mathbf{X}}(\tau)\Big(
\frac{\tilde{V}_{\xi}^{\tilde S_2,-{\mathbf{X}}}}{v} \Big)_{\xi\xi},\notag
    \end{split}
\end{equation}
which together with \eqref{hh1} implies
\begin{equation}\begin{aligned}
\left( \frac{\phi_{\xi\xi}}{v} \right)_{\tau}=&-\frac{\phi_{\xi\xi}\tilde{V}_{\tau}-2\phi_{\xi \tau}v_{\xi}-\phi_{\tau}\tilde{V}_{\xi\xi}  }{v^{2}}-\frac{2\phi_{\tau}v_{\xi}^{2}}{v^{3}}+\tilde{s}_2 \left( \frac{\phi_{\xi}}{v} \right)_{\xi\xi}+\dot{\mathbf{X}}(\tau)\Big(
\frac{\tilde{V}_{\xi}^{\tilde S_2,-{\mathbf{X}}}}{v}-\tilde{U}^{\tilde S_2,{-\mathbf{X}}} \Big)_{\xi\xi}\\
&+\psi _{\tau\xi}-\tilde{s}_2  \psi _{\xi\xi}+\left(p (v)-p (\tilde{V}) \right)_{\xi\xi}+\big(\frac{ \phi\tilde{U} _{\xi}}{v\tilde{V}}\big)_{\xi\xi}+\frac{1}{8}\Big(\frac{\sigma_{\xi}^{2}}{(\sigma\!+\!1)v^{2}}\Big)_{\xi\xi}+{F_1}_{\xi}+{F_2}_{\xi}.\notag
\end{aligned}\end{equation}
Multiplying the above equation by $\frac{\phi_{\xi\xi}}{v}$, integrating the resulting equality by parts, one obtains
\begin{equation}
    \begin{split}
    \frac{d}{d\tau}\int_{\mathbb{R}}\big|\frac{\phi_{\xi\xi}}{v}\big|^{2}d\xi
        &=\dot{\mathbf{X}}(\tau)I_0+\sum_{i=1}^6I_i,\label{514}
    \end{split}
\end{equation}
where
\begin{equation}
\begin{aligned}
I_0&=\int_{\mathbb{R}}\Big(
\frac{\tilde{V}_{\xi}^{\tilde S_2,-{\mathbf{X}}}}{v}-\tilde{U}^{\tilde S_2,{-\mathbf{X}}} \Big)_{\xi\xi}\frac{\phi_{\xi\xi}}{v}d\xi,\quad
I_1=-\int_{\mathbb{R}}\Big(\frac{\phi_{\xi\xi}\tilde{V}_{\tau}-2\phi_{\xi \tau}v_{\xi}-\phi_{\tau}\tilde{V}_{\xi\xi}  }{v^{2}}+\frac{2\phi_{\tau}v_{\xi}^{2}}{v^{3}}\Big)\frac{\phi_{\xi\xi}}{v}d\xi,\\
I_2&=\int_{\mathbb{R}}\sigma\left( \frac{\phi_{\xi}}{v} \right)_{\xi\xi}\frac{\phi_{\xi\xi}}{v}d\xi,\quad
I_3=\int_{\mathbb{R}}\left(\psi _{\tau\xi}-\sigma \psi _{\xi\xi}\right)\frac{\phi_{\xi\xi}}{v}d\xi,\quad
I_4=\int_{\mathbb{R}}\left(p (v)\!-\!p (\tilde{V}) \right)_{\xi\xi}\frac{\phi_{\xi\xi}}{v}d\xi,\\
I_5&=\int_{\mathbb{R}}\big(\frac{ \phi\tilde{U} _{\xi}}{v\tilde{V}}\big)_{\xi\xi}\frac{\phi_{\xi\xi}}{v}d\xi,\quad
I_6=\int_{\mathbb{R}}\Big(\frac{\sigma_{\xi}^{2}}{8(\sigma\!+\!1)v^{2}}\Big)_{\xi\xi}\frac{\phi_{\xi\xi}}{v}d\xi,\quad
I_7=\int_{\mathbb{R}}\Big({F_1}_{\xi}+{F_2}_{\xi}\Big)\frac{\phi_{\xi\xi}}{v}d\xi.\notag
\end{aligned}
\end{equation}
Observing that
\begin{equation}
 \big(p (v)-p (\tilde{V}) \big)_{\xi\xi}= p^{\prime} (v)\phi_{\xi\xi}+(p^{\prime}(v)-p^{\prime}(\tilde{V}))\tilde{V}_{\xi\xi}+p^{\prime\prime}(v)v_{\xi}^{2}-p^{\prime\prime}(\tilde{V})\tilde{V}_{\xi}^{2},\notag
\end{equation}
one gets
\begin{equation}
 I_4=-\underbrace{\int_{\mathbb{R}}\frac{-p^{\prime} (v)}{v}|\phi_{\xi\xi}|^2d\xi}_{\mathbb{D}_2}+\underbrace{\int_{\mathbb{R}}\big((p^{\prime}(v)-p^{\prime}(\tilde{V}))\tilde{V}_{\xi\xi}+
 p^{\prime\prime}(v)v_{\xi}^{2}-p^{\prime\prime}(\tilde{V})\tilde{V}_{\xi}^{2}\big)\frac{\phi_{\xi\xi}}{v}d\xi}_{\mathbb{I}_4},\notag
\end{equation}
where $p^{\prime} (v)<0$ and $\mathbb{D}_2$ is the good term which can used to control the other terms. Now we start with the estimate for $I_0$. Noting that
\begin{equation}
   \Big(\frac{\tilde{V}_{\xi}^{\tilde S_2,-{\mathbf{X}}}}{v}\Big)_{\xi\xi}\frac{\phi_{\xi\xi}}{v}=\Big(\frac{\tilde{V}_{\xi\xi\xi}^{\tilde S_2,-{\mathbf{X}}}}{v}-\frac{2\tilde{V}_{\xi\xi}^{\tilde S_2,-{\mathbf{X}}}v_\xi}{v^{2}}-
   \frac{\tilde{V}_{\xi}^{\tilde S_2,-{\mathbf{X}}}(\phi_{\xi\xi}+\tilde{V}_{\xi\xi}^{\tilde S_2,-{\mathbf{X}}})}{v^{2}}+\frac{2\tilde{V}_{\xi}^{\tilde S_2,-{\mathbf{X}}}v_{\xi}^2}{v^{3}}\Big)\frac{\phi_{\xi\xi}}{v},\notag
\end{equation}
by using \eqref{basic1}, Lemma \ref{the properties of rarefaction wave} and Lemma \ref{estimate of viscous shock}, one has
\begin{equation}\begin{aligned}
 I_0&=-\int_{\mathbb{R}}
\frac{\tilde{V}_{\xi}^{\tilde S_2,-{\mathbf{X}}}}{v^2}|\phi_{\xi\xi}|^2d\xi\\
&+\int_{\mathbb{R}}\Big(\frac{\tilde{V}_{\xi\xi\xi}^{\tilde S_2,-{\mathbf{X}}}}{v}-\frac{2\tilde{V}_{\xi\xi}^{\tilde{S}_2,-{\mathbf{X}}}v_\xi}{v^{2}}-\frac{\tilde{V}_{\xi}^{\tilde{S}_2,-{\mathbf{X}}}\tilde{V}_{\xi\xi}^{\tilde{S}_2,-{\mathbf{X}}}}{v^{2}}+
\frac{2\tilde{V}_{\xi}^{\tilde S_2,-{\mathbf{X}}}v_{\xi}^2}{v^{3}}-\tilde{U}^{\tilde{S}_2,-\mathbf{X}}_{\xi\xi} \Big)\frac{\phi_{\xi\xi}}{v}d\xi\\
&\displaystyle \leq C\tilde\delta_2^2\mathbb{D}_2+C\tilde\delta_2\big(\tilde\delta_2+\tilde\delta_1+M^{\frac12}\big)\mathbb{D}_2^{\frac12},\notag
\end{aligned}
\end{equation}
which together with \eqref{410} yields
\begin{equation}\label{e-i0}
   \dot{\mathbf{X}}(\tau)I_0 \leq C\big(\tilde\delta_2+\tilde\delta_1+M^{\frac12}\big)\mathbb{D}_2+\frac{\tilde\delta_2}{2}|\dot{\mathbf{X}}|^2\leq \frac{1}{16}\mathbb{D}_2+\frac{\tilde\delta_2}{6}|\dot{\mathbf{X}}|^2.
\end{equation}
From \eqref{equation of approximate composite shock wave after interaction}$_1$, \eqref{nsac-per}$_1$, Lemma \ref{the properties of rarefaction wave}, Lemma \ref{estimate of viscous shock} and \eqref{basic1}, one obtains
\begin{equation}\begin{aligned}\label{e-i1}
 I_1&\leq C\int_{\mathbb{R}}\Big((|\tilde{V}_{t}|+|v_\xi|)|\phi_{\xi\xi}|^2+|v_\xi||\psi_{\xi\xi}||\phi_{\xi\xi}|+(|\tilde{V}_\xi|\!+\!|v_\xi|^2)(|\phi_{\xi}|\!+\!|\psi_{\xi}|)|\phi_{\xi\xi}|\Big)d\xi\\
&\qquad +C\int_{\mathbb{R}}|\dot{\mathbf{X}}|\Big(|\tilde{V}_{\xi\xi}^{\tilde S_2,-{\mathbf{X}}}|+|\tilde{V}_{\xi}^{\tilde S_2,-{\mathbf{X}}}||\tilde{V}_{\xi\xi}|+|\tilde{V}_{\xi}^{\tilde S_2,-{\mathbf{X}}}||v_{\xi}|^2 \Big)|\phi_{\xi\xi}|d\xi\\
&\leq C\big(\tilde\delta_2+\tilde\delta_1+M^{\frac12}\big)\big(\mathbb{D}_2+\|\psi_{\xi\xi}\|_{L^2(\mathbb{R})}^2+\|\phi_{\xi}\|_{L^2(\mathbb{R})}^2\big)+\frac{\tilde\delta_2}{4}|\dot{\mathbf{X}}|^2\\
&\leq \frac{1}{16}\mathbb{D}_2+\frac{\tilde\delta_2}{6}|\dot{\mathbf{X}}|^2 +C\big(\tilde\delta_2+\tilde\delta_1+M^{\frac12}\big)\big(\|\psi_{\xi\xi}\|_{L^2(\mathbb{R})}^2+G^{\tilde{R}_1}+G^{\tilde{S}_2}+D_1\big),
\end{aligned}
\end{equation}
where the following facts are used here
\begin{equation}\begin{aligned}
&|\tilde{V}_{\tau}|\leq C\big(|\tilde{V}_{\xi}|+|\tilde{U}_{\xi}|+|\dot{\mathbf{X}}\tilde{V}_{\xi}^{\tilde S_2,-{\mathbf{X}}}|\big)\leq C(\tilde\delta_2+\tilde\delta_1),\\
& \displaystyle \psi_{\tau\xi}=\tilde{s}_2 \phi_{\xi\xi}+\psi_{\xi\xi}+\dot{\mathbf{X}}\tilde{V}_{\xi\xi}^{\tilde S_2,-{\mathbf{X}}}.\label{vvv}
\end{aligned}\end{equation}
Similarly, a straightforward computation yields
\begin{equation}\begin{aligned}\label{e-i2}
 I_2&=\underbrace{\int_{\mathbb{R}}\tilde{s}_2 \Big( \frac{\phi_{\xi\xi}}{v} \Big)_{\xi}\frac{\phi_{\xi\xi}}{v}d\xi}_{=0}-\int_{\mathbb{R}}\sigma\Big( \frac{\phi_{\xi}v_\xi}{v^2} \Big)_{\xi}\frac{\phi_{\xi\xi}}{v}d\xi\\
 &\leq\int_{\mathbb{R}}|v_\xi||\phi_{\xi\xi}|^2d\xi+
\int_{\mathbb{R}}(|\tilde{V}_{\xi\xi}|+|v_\xi|^2)|\phi_{\xi}||\phi_{\xi\xi}|d\xi\\
&\leq C\big(\tilde\delta_2+\tilde\delta_1+M^{\frac12}\big)(\mathbb{D}_2+\|\phi_{\xi}\|_{L^2}^2)\\
&\leq \frac{1}{16}\mathbb{D}_2+C(\tilde\delta_2+\tilde\delta_1+M^{\frac12})\big(G^{\tilde{R}_1}+G^{\tilde{S}_2}+D_1\big).
\end{aligned}
\end{equation}
For $I_3$, by using \eqref{vvv} (as done in \eqref{e-i1}), one gets
\begin{equation}\label{e-i3}
    \begin{split}
I_3&=\int_{\mathbb{R}}\left(\psi _{\tau\xi}-\tilde{s}_2 \psi _{\xi\xi}\right)\frac{\phi_{\xi\xi}}{v}d\xi,\\
        &=\frac{d}{d \tau} \int_\mathbb{R} \frac{\psi_\xi \phi_{\xi\xi}}
        {v} d \xi-\int_\mathbb{R} \psi_\xi\left(\frac{\phi_{\xi\xi}}{v}\right)_\tau d \xi -\int_{\mathbb{R}}\tilde{s}_2 \psi _{\xi\xi}\frac{\phi_{\xi\xi}}{v}d\xi,\\
 & =\frac{d}{d \tau} \int_\mathbb{R} \frac{\psi_\xi \phi_{\xi\xi}}
        {v} d \xi+\int_\mathbb{R}\left(\frac{\psi_{\xi}}{v}\right)_\xi \phi_{\tau \xi} d \xi+\int_{\mathbb{R}}\frac{\psi_{\xi} \phi_{\xi\xi}}{v}\left(\phi_\tau+\tilde{V}_\tau\right) d \xi-\int_{\mathbb{R}}\tilde{s}_2  \psi _{\xi\xi}\frac{\phi_{\xi\xi}}{v}d\xi,\\
        &\leq \frac{d}{d \tau} \int_\mathbb{R} \frac{\psi_\xi \phi_{\xi\xi}}
        {v} d \xi+\frac{1}{16}\mathbb{D}_2+\frac{\tilde\delta_2}{6}|\dot{\mathbf{X}}|^2 +C\|\psi_{\xi}\|_{H^1(\mathbb{R})}^2.
    \end{split}
\end{equation}
and for $\mathbb{I}_4$, observing that
\begin{equation}\begin{aligned}
&\big(p^{\prime}(v)-p^{\prime}(\tilde{V})\big)\tilde{V}_{\xi\xi}+p^{\prime\prime}(v)v_{\xi}^{2}-p^{\prime\prime}(\tilde{V})\tilde{V}_{\xi}^{2}
=p^{\prime\prime}(v_1)\phi\tilde{V}_{\xi\xi}+p^{\prime\prime}(v)\phi_{\xi}^{2}+p^{\prime\prime\prime}(v_2)\phi\tilde{V}_{\xi}^{2},\notag
\end{aligned}
\end{equation}
where $v_1,v_2$ are constants between $v$ and $\tilde{V}$.
This together with Lemma \ref{the properties of rarefaction wave} and Lemma \ref{estimate of viscous shock} yields
\begin{equation}\begin{aligned}\label{e-i4}
 \mathbb{I}_4 &\leq C\int_{\mathbb{R}}|\phi_{\xi}|^2|\phi_{\xi\xi}|d\xi+C
\int_{\mathbb{R}}|\phi|\big(|\tilde{V}_{\xi\xi}|+|\tilde{V}_\xi|^2\big)|\phi_{\xi\xi}|d\xi\\
&\leq \frac{1}{16}\mathbb{D}_2+C(\tilde\delta_2+\tilde\delta_1+M^{\frac12})(G^{\tilde{R}_1}+G^{\tilde{S}_2}+D_1).
\end{aligned}
\end{equation}
Similarly, one gets
\begin{equation}\begin{aligned}\label{e-i5}
 I_5 &\leq C\int_{\mathbb{R}}\Big(|\tilde{U}_{\xi}||\phi_{\xi\xi}|^2+|\phi_\xi|\big(|\tilde{U}_{\xi\xi}|+|\tilde{V}_\xi||\tilde{U}_\xi|\big)|\phi_{\xi\xi}|+|\phi||\tilde{U}_{\xi}||\tilde{V}_{\xi\xi}||\phi_{\xi\xi}|\Big)d\xi\\
&\leq \frac{1}{16}\mathbb{D}_2+C(\tilde\delta_2+\tilde\delta_1+M^{\frac12})\big(G^{\tilde{R}_1}+G^{\tilde{S}_2}+D_1\big).
\end{aligned}
\end{equation}
For $I_6$, by \eqref{basic1}, Lemma \ref{the properties of rarefaction wave} and Lemma \ref{estimate of viscous shock},  a straightforward computation yields
\begin{equation}\begin{aligned}\label{e-i6}
 I_6 &\leq C\int_{\mathbb{R}}|\sigma_{\xi}|^2|\phi_{\xi\xi}|^2d\xi+C
\int_{\mathbb{R}}|\sigma_\xi|(|\sigma_\xi|+|\sigma_{\xi\xi\xi}|)|\phi_{\xi\xi}|d\xi+C
\int_{\mathbb{R}}|\sigma_{\xi\xi}|^2|\phi_{\xi}|^2d\xi\\
&\quad +C\int_{\mathbb{R}}|\sigma_{\xi\xi}||\sigma_{\xi\xi\xi}||\phi_{\xi}|d\xi+C\int_{\mathbb{R}}|\sigma_{\xi\xi}|^2|\sigma_{\xi}||\phi_{\xi}|d\xi\\
&\leq \frac{1}{16}\mathbb{D}_2+CM^{\frac12}\|\sigma_{\xi\xi}\|_{H^1(\mathbf{}\mathbb{R})}^2.
\end{aligned}
\end{equation}
Finally, as the definition of \eqref{310} and straight calculation, one obtains
\begin{equation}
    \begin{split}
F_{1\xi}&=\Big(\big(\frac{\tilde{U}^{\tilde S_2,{-\mathbf{X}}}_{\xi} }{\tilde{V}^{\tilde S_2,{-\mathbf{X}}}}\big)_{\xi}-\big(\frac{\tilde{U}_{\xi}}{\tilde{V}}\big)_{\xi}\Big)_{\xi}\\
& =\Big(\frac{\tilde{U}^{\tilde S_2,{-\mathbf{X}}}_{\xi\xi}-\tilde{U}_{\xi\xi}}{\tilde{V}^{\tilde S_2,{-\mathbf{X}}}}+\Big(\frac{1}{\tilde{V}^{\tilde S_2,{-\mathbf{X}}}}-\frac{1}{\tilde{V}}\Big) \tilde{U}_{\xi\xi} \\
& \qquad +\tilde{U}_\xi \tilde{V}_\xi\Big(\frac{1}{\tilde{V}^2}-\frac{1}{\big(\tilde{V}^{\tilde S_2,{-\mathbf{X}}}\big)^2}\Big)+\frac{\tilde{U}_\xi \tilde{V}_\xi-\tilde{U}^{\tilde S_2,{-\mathbf{X}}}_\xi\tilde{V}^{\tilde S_2,{-\mathbf{X}}}_\xi}{\big(\tilde{V}^{\tilde S_2,{-\mathbf{X}}}\big)^2} \Big)_{\xi} \\
&=\Big(\frac{\tilde{U}^{\tilde{R}_1}_{\xi\xi}}{\tilde{V}^{\tilde S_2,{-\mathbf{X}}}}+\frac{(\tilde{V}^{\tilde{R}_1}-v_{m})\tilde{U}_{\xi\xi}}{\tilde{V}\tilde{V}^{\tilde S_2,{-\mathbf{X}}}} \\
&\qquad -\frac{\tilde{U}_\xi \tilde{V}_\xi\big(\tilde{V}+\tilde{V}^{\tilde S_2,{-\mathbf{X}}}\big)\big(\tilde{V}^{\tilde{R}_1}-v_m\big)}{\tilde{V}^2\big(\tilde{V}^{\tilde S_2,{-\mathbf{X}}}\big)^2}+\frac{\tilde{U}_\xi \tilde{V}_\xi^{\tilde{R}_1}+\tilde{U}_{\xi}^{\tilde{R}_1}\tilde{V}^{\tilde S_2,{-\mathbf{X}}}_{\xi} }{\big(\tilde{V}^{\tilde S_2,{-\mathbf{X}}}\big)^2}\Big)_{\xi}\\
&\leq C\Big(|\tilde{U}^{\tilde{R}_1}_{\xi\xi}|+|\tilde{V}^{\tilde{R}_1}_{\xi\xi}|+|\tilde{U}^{\tilde{R}_1}_{\xi\xi\xi}|+|\tilde{V}^{\tilde{R}_1}-v_m||\tilde{V}^{\tilde S_2,{-\mathbf{X}}}_{\xi}|+|\tilde{V}_{\xi}^{\tilde{R}_1}|^{2}+|\tilde{U}_{\xi}^{\tilde{R}_1}|^{2}
\Big),\notag
    \end{split}
\end{equation}
and
\begin{equation}
    \begin{split}
        F_{2\xi}&=\Big(p (\tilde{V}) -p\big(\tilde{V}^{\tilde{R}_1}\big)-p\big(\tilde{V}^{\tilde S_2,{-\mathbf{X}}}\big)\Big)_{\xi\xi} \\
        &=\Big(\big(p^{\prime}(\tilde{V})-p^{\prime}\big(\tilde{V}^{\tilde{R}_1}\big)\big) \tilde{V}_\xi^{\tilde{R}_1}+\big(p^{\prime}(\tilde{V})-p^{\prime}\big(\tilde{V}^{\tilde S_2,{-\mathbf{X}}}\big)\big)\tilde{V}^{\tilde S_2,{-\mathbf{X}}}_{\xi}\Big)_{\xi}\\
        &\leq C\left( |\tilde{V}_{\xi\xi
        }^{\tilde{R}_1}|+|\tilde{V}_{\xi}^{\tilde{R}_1}|^{2}+|\tilde{V}_{\xi}^{\tilde{R}_1}||\tilde{V}^{\tilde S_2,{-\mathbf{X}}}_{\xi}|+|\tilde{V}^{\tilde{R}_1}-v_{m}||\tilde{V}^{\tilde S_2,{-\mathbf{X}}}_{\xi}| \right),\notag
    \end{split}
\end{equation}
which yields that
\begin{equation}\begin{aligned}\label{e-i7}
 I_7 &\leq C\int_{\mathbb{\tilde{R}_1}}\left( |\tilde{V}_{\xi\xi
        }^{\tilde{R}_1}|+|\tilde{V}_{\xi}^{\tilde{R}_1}|^{2}+|\tilde{V}_{\xi}^{\tilde{R}_1}||\tilde{V}^{\tilde S_2,{-\mathbf{X}}}_{\xi}|+|\tilde{V}^{\tilde{R}_1}-v_{m}||\tilde{V}^{\tilde S_2,{-\mathbf{X}}}_{\xi}| \right)\phi_{\xi\xi}d\xi\\
&\leq \frac{1}{16}\mathbb{D}_2+C\left( \|\tilde{V}_{\xi\xi
        }^{\tilde{R}_1}\|_{L^2(\mathbb{R})}^2+\|\tilde{V}_{\xi}^{\tilde{R}_1}\|_{L^4(\mathbb{R})}^4\right)+C\tilde\delta_1\tilde\delta_2^{\frac32}e^{-C\tilde\delta_2\tau}.
\end{aligned}
\end{equation}
Combining \eqref{e-i0}-\eqref{e-i7} with \eqref{514}, one has
\begin{equation}
    \begin{aligned}
    \frac{d}{d\tau}\int_{\mathbb{R}}\big|\frac{\phi_{\xi\xi}}{v}\big|^{2}d\xi+\frac12\mathbb{D}_2
         &\leq \frac{d}{d t} \int_\mathbb{R} \frac{\psi_\xi \phi_{\xi\xi}}
        {v} d \xi+ \frac{\tilde\delta_2}{2}|\dot{\mathbf{X}}|^2\!+\!C \|\psi_{\xi}\|_{H^1(\mathbb{R})}^2\\
        &\quad+\!C\big(\tilde\delta_2+\tilde\delta_1+M^{\frac12}\big)\big(G^{\tilde{S}_2}+G^{\tilde{R}_1}+D_1+\|\sigma_{\xi\xi}\|_{H^1}^2\big) \\
& \quad +C\big\|\tilde{U}^{\tilde{R}_1}_{\xi \xi}\big\|_{L^2(\mathbb{R})}^2+C\big\|\tilde{U}^{\tilde{R}_1}_{\xi}\big\|_{L^4(\mathbb{R})}^4+C\tilde\delta_1\tilde\delta_2^{\frac32}e^{-c\tilde\delta_2\tau} .\label{pho-h}
\end{aligned}
\end{equation}
Differentiating \eqref{nsac-per}$_3$ with respect to $\xi$ twice, multiplying the resulting by $\sigma_{\xi\xi}$, and integrating that inequality by part, one gets
 \begin{equation}
    \begin{split}
       & \frac{1}{2}\frac{d}{d\tau}\Big(\int_{\mathbb{R}}|\sigma_{\xi\xi}|^{2}d\xi\Big)+\int_{\mathbb{R}}|\sigma_{\xi\xi\xi}|^{2}d\xi+\int_{\mathbb{R}}2v|\sigma_{\xi\xi}|^{2}d\xi\\
       &=\int_{\mathbb{R}}\Big( 2v\sigma^{2}+\frac{\sigma_{\xi}v_{\xi}}{v}+\frac{\sigma_{\xi}^{2}}{2(\sigma+1)} \Big)_{\xi}\sigma_{\xi\xi\xi}d\xi-2\int_{\mathbb{R}}(v_{\xi\xi}\sigma\sigma_{\xi\xi}+2v_\xi\sigma_\xi\sigma_{\xi\xi})d\xi\\
  &\leq \frac{1}{2}\|\sigma_{\xi\xi\xi}\|^{2}+C\epsilon_{1}^{\frac12}\big(\mathbb{D}_2+\|\sigma\|_{H^3(\mathbb{R})}^2\big),\notag
    \end{split}
\end{equation}
which implies
 \begin{equation}\label{phi-hh}
           \frac{1}{2}\frac{d}{d\tau}\Big(\int_{\mathbb{R}}|\sigma_{\xi\xi}|^{2}d\xi\Big)+\frac{1}{2}\|\sigma_{\xi\xi\xi}\|^{2}+\int_{\mathbb{R}}2v|\sigma_{\xi\xi}|^{2}d\xi\leq C\epsilon_{1}^{\frac12}\left(\mathbb{D}_2+\|\sigma\|_{H^3(\mathbb{R})}^2\right).
    \end{equation}
Adding \eqref{pho-h} to \eqref{phi-hh} and integrating the resulting over $[0, \tau]$ for any $\tau \leq T$, using \eqref{l2l4}, we have that for some $c_3>0$,
 \begin{equation}
    \begin{split}
       & \int_{\mathbb{R}}\Big(\big|\frac{\phi_{\xi\xi}}{v}\big|^{2}+\frac{1}{2}|\sigma_{\xi\xi}|^{2}\Big)d\xi+\frac12\int_0^\tau(\mathbb{D}_2+\|\sigma_{\xi\xi}\|_{H^1(\mathbb{R})}^{2})ds\\
       &\leq C\int_{\mathbb{R}}\Big(\big|\frac{\phi_{0\xi\xi}}{v_0}\big|^{2}+\frac{1}{2}|\sigma_{0\xi\xi}|^{2}+|\psi_{0\xi}|^2\Big)d\xi+\int_0^\tau\Big(\frac{\tilde\delta_2}{2}|\dot{\mathbf{X}}|^2+c_3 \|\psi_{\xi}\|_{H^1(\mathbb{R})}^2\\
        &\quad+\!C(\tilde\delta_2+\tilde\delta_1+M^{\frac12})(G^{\tilde{S}_2}+G^{\tilde{R}_1}+D_1+\mathbb{D}_2+\|\sigma\|_{H^3(\mathbb{R})}^2)\Big)ds.
    \end{split}
\end{equation}
Multiplying the above estimate by the constant $\frac{1}{2 \max \left(1, c_3\right)}$, and then adding the result to \eqref{5111}, together with the smallness of $M, \tilde\delta_2, \tilde\delta_1$, one obtains
\begin{equation}
    \begin{split}
       & \|(\phi,\sigma)\|_{H^{2}(\mathbb{R})}^{2}\!+\!\|\psi\|_{H^{1}(\mathbb{R})}^{2}+\tilde\delta_{2}\int_{0}^{\tau}|\dot{\mathbf{X}}|^{2}ds\\
       &\qquad+\int_{0}^{\tau}\left(G^{\tilde{R}_1}+G^{\tilde{S}_2}+D_1+\|\phi_{\xi\xi}\|_{L^{2}(\mathbb{R})}^{2}+\|\psi_\xi\|_{H^{1}(\mathbb{R})}^{2}+\|\sigma\|_{H^{3}(\mathbb{R})}^{2}\right)ds\\
        &\leq C(\|\phi_{0},\sigma_{0}\|_{H^{2}(\mathbb{R})}^{2}+\|\psi_{0}\|_{H^{1}(\mathbb{R})}^{2})+C\tilde\delta_1^{\frac{1}{3}},\notag
    \end{split}
\end{equation}
and the proof of Lemma \ref{lem-vp} is completed.
\end{proof}

{\it Proof of Proposition \ref{prop}.}
Combining with Lemmas \ref{lem-low}-\ref{lem-vp}, Lemma \ref{the properities for Shift X} and the following estimates
\begin{equation}
    \int_{\mathbb{R}}\big(|\tilde{V}^{\tilde S_2,{-\mathbf{X}}}_{\xi}|+|\tilde{U}^{\tilde{R}_1}_{\xi}|\big)\phi^2 d \xi\leq C\big(G^{\tilde{S}_2}+G^{\tilde{R}_1}\big),\notag
\end{equation}
 the desired estimate \eqref{pro1} is obtained. Thus the proof of Proposition \ref{prop} is completed. \endproof

\subsection{The proof of Theorem 4.1}\label{sec-g}

\indent\qquad
With all the a priori estimates in Subsection 2-4 at hand, we are going to  prove the main theorems in this subsection. By \eqref{per1}, Theorem \ref{thm1} can be deduced immediately from  Lemma \ref{the properties of rarefaction wave}. Moreover, using a natural scaling argument \eqref{scaling method}, we reduce the proof of Theorem 1.1 from the
nonlinear asymptotic stability analysis of composite wave in Theorem \ref{thm1}.

\vskip 0.3cm
{\it Proof of Theorem 4.1.}
In order to prove the global existence of the solutions to the problem \eqref{nsac-per}-\eqref{init-per}, we employ the standard continuation argument based on a local existence theorem and the a priori estimates.
Once Proposition \ref{prop} is obtained, similar to \cite{kvw2021}, we employ the standard continuation argument to extend the local solution to a global one.
Therefore, to complete the proof of Theorem 4.1, we need only to investigate the large-time behavior for the solution $(\phi,\psi,\sigma)$ of system \eqref{nsac-per}. First, we define a function $g$ on $(0,+\infty)$ by
 \begin{gather}
     g(\tau)\xlongequal{\mathrm{def}} \|(\phi,\psi,\sigma)_{\xi}\|_{L^{2}(\mathbb{R})}^{2}.\notag
 \end{gather}
The main aim is to show the following estimate
\begin{equation}\label{527}
    \int_{0}^{\infty}( |g(\tau)|+|g^{\prime}(\tau)| )d\tau<\infty.
\end{equation}
Using \eqref{5144}, one can easily get
\begin{gather}
    \int_{0}^{\infty}|g(\tau)| d\tau<\infty.\label{5288}
\end{gather}
Using \eqref{5111}, \eqref{5144} and \eqref{nsac-per}, one has
\begin{equation}
    \begin{split}
 & \int_{0}^{\infty}|g^{\prime}(\tau)|dt=\int_{0}^{\infty}2\Big|\int_{\mathbb{R}}\big( \phi_{\xi}\phi_{\tau\xi}+\psi_{\xi}\psi_{\tau\xi}+\sigma_{\xi}\sigma_{\tau\xi}\big)d\xi \Big|  d\tau    \\
  &\leq 2\int_{0}^{\infty}\int_{\mathbb{R}}\left| \phi_{\xi\xi}\phi_{\tau}+\psi_{\xi\xi}\psi_{\tau}+\sigma_{\xi\xi}\sigma_{\tau}\right|d\xi   d\tau\\
  &\leq C\int_{0}^{\infty}\Big(\|\phi,\psi,\sigma\|_{H^{2}(\mathbb{R})}^{2}+\tilde\delta_{2}\dot{\mathbf{X}}^{2}+|D_1|\Big)d\tau
  +C\int_{0}^{\infty}\int_{\mathbb{R}}\Big(F_{1}^{2}+F_{2}^{2}+\big(\frac{u_{\xi}}{v}-\frac{\tilde{U}_{\xi}}{\tilde{V}}\big)_{\xi}^{2} \Big)d\xi d\tau\\
  &\leq C\int_{0}^{\infty}\left(\|\phi,\psi,\sigma\|_{H^{2}(\mathbb{R})}^{2}+\tilde\delta_{2}\dot{\mathbf{X}}^{2}+|D_1|+|G^{\tilde{S}_2}|+|G^{\tilde{R}_1}|\right)d\tau\\
  &\quad +C\int_{0}^{\infty}\Big\| (\tilde{U}^{\tilde{R}_1}_{\xi})^{2}+ |\tilde{U}^{\tilde{R}_1}_{\xi \xi}|+ |\tilde{V}^{\tilde S_2,{-\mathbf{X}}}_{\xi} |\cdot |\tilde{V}^{\tilde{R}_1}-v_m|+|\tilde{V}^{\tilde{R}_1}_{\xi} ||\tilde{V}^{\tilde S_2,{-\mathbf{X}}}_{\xi}|\Big\|_{L^2(\mathbb{R})}^2 d\tau\\
  &\leq C.\notag
    \end{split}
\end{equation}
Thus, the above estimates with \eqref{5288}, we get \eqref{527}. Therefore, one obtains
\begin{gather}
    \lim_{\tau\to+\infty}\| \left(\phi, \psi, \sigma \right)_{\xi} \|_{L^{2}(\mathbb{R})}=0,
\end{gather}
this along with interpolation inequality implies
\begin{gather}
    \lim_{\tau\to+\infty}\| (\phi, \psi, \sigma)  \|_{L^{\infty}(\mathbb{R})}=0,
\end{gather}
which together with \eqref{410} yields
\begin{equation}\label{thm3-lim}
    \lim_{\tau\rightarrow +\infty}\sup _{\xi \in \mathbb{R}}\left|(\phi,\psi,\sigma)(\xi,\tau)\right| = 0,\quad \lim_{\tau\rightarrow +\infty}|\dot{\mathbf{X}}(\tau)|= 0.
\end{equation}
 and the proof of Theorem 4.1 is completed. \endproof

\

Finally,Theorems 1.1 and Theorem 1.2 are directly derived from Theorem 3.1 and Theorem 4.2, and the proof of the main conclusions in this paper is completed.

\section*{Acknowledgments:} 
Qiaolin He and Yazhou Chen acknowledge support from  National Natural Sciences Foundation of China  (NSFC) (No. 12371434). Qiaolin He acknowledges support from the National key R \& D Program of China (No.2022YFE03040002).
Xiaoping Wang acknowledges support from National Natural Science Foundation of China (NSFC) (No. 12271461), the key project of NSFC (No. 12131010),   the University Development Fund from The Chinese University of Hong Kong, Shenzhen (UDF01002028), and Hetao Shenzhen-Hong Kong Science and Technology  Innovation Cooperation Zone Project (No.HZQSWS-KCCYB-2024016).  Xiaoding Shi acknowledges support from  National Natural Sciences Foundation of China  (NSFC) (No. 12171024).

\section*{ Conflict of Interest} The authors declare that they have no conflict of interest.

\end{document}